\documentclass[12pt, twoside, a4paper]{article}
\usepackage[english]{babel}  
\usepackage[utf8]{inputenc}  
\usepackage[T2A]{fontenc}      
\usepackage{wrapfig}
\usepackage{amsmath, amssymb, amsthm, amsfonts, amsthm, dsfont,mathrsfs}
\usepackage{graphicx} 
\usepackage{accents,epigraph}
\usepackage{abstract}
\DeclareMathOperator{\supp}{supp}
\DeclareMathOperator{\dist}{dist}

\DeclareMathOperator{\const}{const}
\DeclareMathOperator{\clos}{clos}
\DeclareMathOperator{\loc}{{loc}}
\DeclareMathOperator{\Harm}{Harm}

\DeclareMathOperator{\Innt}{Int}
\DeclareMathOperator{\Per}{{Per}}
\DeclareMathOperator{\Cap2}{{Cap_2}}
\DeclareMathOperator{\diam}{{diam}}
\DeclareMathOperator{\Cylr}{{Cyl_r}}

\DeclareMathOperator{\arctanh}{{arctanh}}

\DeclareMathOperator{\card}{card}
\DeclareMathOperator{\divv}{div}
\DeclareMathOperator{\BS}{BS}
\DeclareMathOperator{\rot}{curl}

\newcommand{\eps}{\varepsilon}
\newcommand{\R}{{\mathbb R}}
\newcommand{\lo}{{{L^{2,1}_c}(\Omega)}}
\newcommand{\dl}{\,d\lambda_2}
\newcommand{\mao}{{\mathcal Adm(\Omega)}}
\newcommand{\sob}{{\overset{\circ}{W}{}^{1,2}(\mathbb D)}}

\newtheorem{theorem}{Theorem}[section]
\newtheorem{lemma}[theorem]{Lemma}
\newtheorem{define}[theorem]{Definition}
\newtheorem{predl}[theorem]{Proposition}
\newtheorem{sled}[theorem]{Corollary}

\usepackage{fancyhdr}
\makeatletter
\def\blindfootnote{\gdef\@thefnmark{}\@footnotetext}
\makeatother

\begin{document}
\thispagestyle{empty}

\title{Interpolation by periods in planar domain}
\author{M.B. Dubashinskiy} 
\date{\today}
\maketitle

\setlength{\epigraphrule}{0pt}
\epigraphwidth=5cm
\epigraph{Dedicated to the memory of Victor Petrovich Havin.}{}

\blindfootnote{Chebyshev Laboratory, St.~Petersburg State University, 14th Line 29b, Vasilyevsky Island, Saint~Petersburg 199178, Russia.}
\blindfootnote{\hspace{0.2mm}e-mail: \texttt{mikhail.dubashinskiy@gmail.com}}
\blindfootnote{The work is supported by the Russian Science Foundation grant 14-21-00035.}
\blindfootnote{{Keywords:} \emph{infinite-connected domain, periods of forms, interpolation, Riesz basis, harmonic functions}.}
\blindfootnote{\hspace{0.0mm}{MSC 2010 Primary}: 30C85; {Secondary}: 31A15, 30E05, 30H20, 58A14, 26D15.}

\thispagestyle{empty}

\renewcommand{\abstractname}{}

\begin{abstract}\footnotesize
\noindent {\bf Abstract.} Let $\Omega\subset\R^2$ be a countably connected domain. To any closed differential form of degree $1$ in $\Omega$ with components in $L^2(\Omega)$ one associates the sequence of its periods around holes in $\Omega$, that is around bounded connected components of $\R^2\setminus \Omega$. \emph{For which $\Omega$ the collection of such period sequences coincides with $\ell^2$?} We give the answer in terms of metric properties of holes in~$\Omega$.
\end{abstract}

{\footnotesize
	\tableofcontents
}

\pagestyle{fancy}

\section*{Introduction}\markboth{Introduction}{Introduction}
\addcontentsline{toc}{section}{Introduction}

Let $\mathbb D$ be the unit disk in the plane $\mathbb C \simeq \R^2$ and connected compact sets $B_1, B_2, \dots\subset \mathbb D$ be disjoint. Consider a planar countably connected domain $\Omega:=\mathbb D\setminus\bigcup_{j=1}^\infty B_j$; sets  $B_1, B_2, \dots$ are called \emph{holes} in $\Omega$ (unbounded connected component of $\R^2\setminus \Omega$ will have a special status). Denote by $\lo$ the Hilbert space of all closed real differential forms $\omega=\omega_x dx+\omega_y dy$ of degree $1$ defined in $\Omega$ and such that components $\omega_x$ and $\omega_y$ are square-integrable over Lebesgue measure in $\Omega$ (one can also deal with vector fields which satisfy analogous conditions). For $\omega\in\lo$ and $j=1,2,\dots$, let $\Per_j\omega$ be the period of $\omega$ around hole $B_j$ (see subsection~\ref{subsection:definitions} for details). Further, define \emph{period operator}: put $\Per\omega:=\{\Per_j\omega\}_{j=1}^\infty$ for $\omega\in\lo$. 

We study the following question:  \emph{for which domains $\Omega$ the set of sequences of the kind $\Per\omega$ with $\omega\in\lo$ {coincides} with $\ell^2$?} We say that such an $\Omega$ has \emph{complete interpolation property} (more precisely, \emph{complete interpolation by periods property}).

The question just stated is similar to classic problems on interpolation of a function from some analytic space by its values in prescribed points, see e.g.~\cite{Seip},~\cite{Nikolskiy}. Moreover, our problem is equivalent to the analogous problem on interpolation of  functions from unweighted Bergman space in $\Omega$ by their complex periods (see subsection~\ref{subsec:Bergman}). 

Note that the complete interpolation property of domain $\Omega$ can be reformulated in the language of homology theory. Namely, the system of curves in $\Omega$ linked to holes $B_j$ should be a Riesz basis (see~\cite{Bari}, ~\cite{Nikolskiy}) in the space of $L^2$-homologies in $\Omega$ (see remark at the end of subsection~\ref{subsection:Riesz}).

Note also that the question on interpolation by periods is reduced by duality to estimate on functions harmonic in $\Omega$ with locally-constant Dirichlet boundary data (proposition~\ref{criterii_adm}). This estimate can be of an independent interest.

The main result of this paper is description of domains $\Omega$ possessing complete interpolation property in terms of metric characteristics of mutual layout of sets $B_j, \, j\in\mathbb N$ (theorem~\ref{th:full_interp_criterion}).

The most difficult is to derive the first condition of theorem~\ref{th:full_interp_criterion} (uniform local finiteness property) from complete interpolation property; this is done in theorem~\ref{th:uniform_local_finiteness}. Nevertheless, the first three conditions of theorem~\ref{th:full_interp_criterion} look naturally from the viewpoint of classical theory of interpolation in analytic spaces, where, roughly speaking, Blaschke condition and separatedness conditions correspond to them. The fourth condition of theorem~\ref{th:full_interp_criterion} (\emph{capacity connectedness}) seems to be the most exotic. This condition is connectedness of a graph consisting of condensers  formed by connected components of $\R^2\setminus \Omega$ and having not very small capacity.

\paragraph*{Motivation to the problem.} \addcontentsline{toc}{subsection}{Motivation}The interest to the planar problem on interpolation by periods is raised by the following high-dimensional problem on energy minimizing electric current. In $\mathbb R^3$ consider a closed domain $K$ bounded by a surface homeomorphic to sphere with three handles (see fig.~\ref{fig:just_krendel}). Put three closed oriented curves $\gamma_1, \gamma_2, \gamma_3$ supported on boundary $\partial K$  and bounding oriented surfaces $S_1, S_2, S_3$ (\emph{sections}) respectively in such a manner that the set $K\setminus\left(S_1\cup S_2\cup S_3\right)$ is  simply connected. Pick $a_1, a_2, a_3 \in \R$. Let us find electric current $\vec I$ supported on $K$ and such that:

\begin{wrapfigure}[19]{l}{0.46\textwidth}
	\vspace{-0.5mm}
	{\includegraphics[trim=0.6cm 1.3cm 1cm 1.2cm, clip=true, scale=1.2]
		{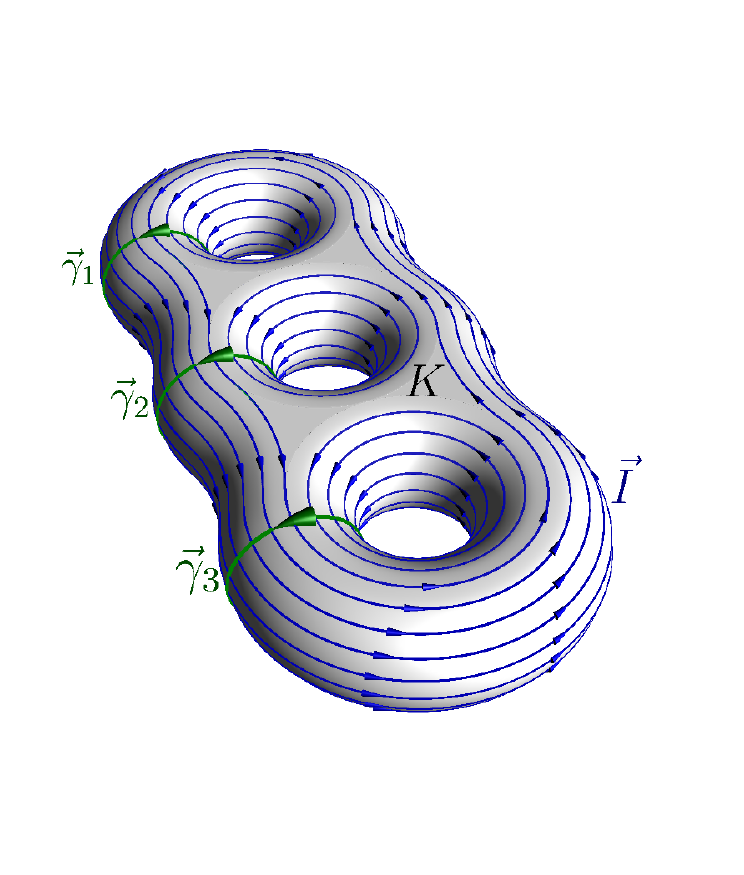}}
	\caption{Equilibrium current on multiply-connected set}
	\label{fig:just_krendel}
\end{wrapfigure}

\begin{enumerate} 
\item $\divv\vec I=0$;
\item $\int_{S_j} \langle\vec I, \vec{n}\rangle \,d\mathcal H^2=a_j$ for $j=1,2,3$; here $\vec n$ is unit normal field orienting surfaces $S_j$ and $\mathcal H^2$ is area measure on these surfaces;
\item current $\vec I$ gives minimum to the energy
$$
W(\vec I):=\int\limits_{\R^3}\int\limits_{\R^3}\dfrac{\langle{\vec I(x), \vec I(y)}\rangle}{|x-y|}\,dx\,dy.
$$
\end{enumerate}

\noindent Such an energy minimizing current  $\vec I$ is called \emph{equilibrium current}. This notion was introduced in~\cite{Yamaguchi}.  If $W(\vec I)<+\infty$ then  $\vec I$ is a distributional vector field (or de Rham current, see~\cite{DeRham}) whose components in some coordinate system can be represented as a sum of distributional partial derivatives of $L^2$-functions. One can properly define flows of such currents through surfaces $S_j$.

If the boundary $\partial K$ is smooth enough then the equilibrium current $\vec I$ can be found by solution of the Laplace equation for multivalued function in $\R^3\setminus K$. Nevertheless, we wish to establish theory of equilibrium currents in the same generality in which the classical capacity theory is developed. For this, one has to state the problem on equilibrium currents for arbitrary $K\subset \R^3$.  Non-regularity of the boundary $\partial K$ can be of several kinds. First, this can be a "differential"{} non-smoothness, for example, if there are "cusps"{} on~$\partial K$; a number of attempts is devoted to the behaviour of solutions of the Laplace operator near such singularities. Another type of non-smoothness of $K$ is its topological complexity. This occurs if $K$ has infinitely many holes, that is if homology spaces of $\R^3 \setminus K$ are infinite-dimensional. We immediately pass to the question: \emph{what should be the appropriate norming conditions on the flow of $\vec I$ through sections $S_j$? to what extent does such a problem depend on the choice of system of sections?} Author does not know anything about this.

One can also ask the following: \emph{which compact sets $K\subset \R^3$ can support non-zero current with finite energy?} It easy to see that then $K$ can also support a non-zero scalar charge with finite energy. It can be proved that \emph{there exist sets $K\subset \R^3$ supporting scalar charges of finite energy but not supporting solenoidal currents with finite energy}. One can take an appropriate Cantor-type set for such $K$. \begin{wrapfigure}[25]{r}{0.45\textwidth}
	\vspace{-3mm}
	{\includegraphics[trim=0.0cm 0.3cm 0cm 0.2cm, clip=true, scale=0.075]
		{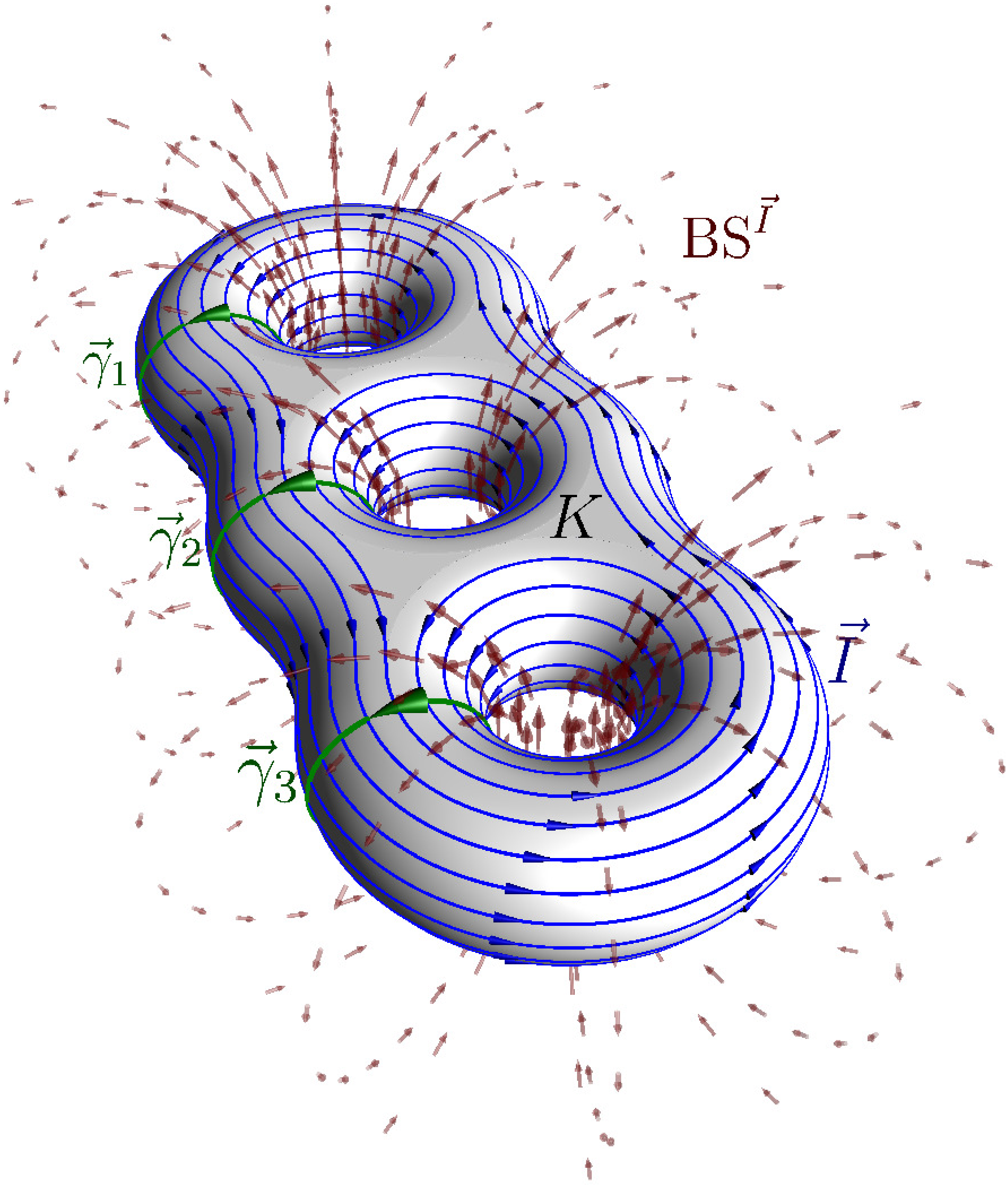}}
	\caption{Biot-Savart field}
	\label{fig:krendel_Field}
\end{wrapfigure}Thus, solenoidality condition $\divv\vec I=0$ makes our problem to be different from the classical question on the sets of positive capacity.

Let us also note that, unlike the classical capacity theory, if $K$ supports a current with finite energy then $K$ may not support a current with finite energy and with finite mass.

Now let us point out the relation between the problem on equilibrium current and the problem on interpolation by periods. An electric current $\vec I$ generates \emph{magnetic Biot-Savart field} (see fig.~\ref{fig:krendel_Field}) given by equality
$$
\BS^{\vec I}=\frac{1}{4\pi}\rot\left(\vec I\star 1/|x|\right).
$$
Here $\star$ denotes convolution and operator $\rot$ is understood in the sense of distributions. The role of such potentials in approximation theory was studied in~\cite{HavinPresa} and~\cite{MalinnikovaHavin}.

If $\vec I$ is an equilibrium current then field~$\BS^{\vec I}$ has to vanish 
in the interior of~$K$. Put $\vec f=\BS^{\vec I}$. One can show that the 
above-stated problem on equilibrium current is equivalent to the following 
problem on its magnetic field: \emph{to find a field $\vec f$ in  $\R^3\setminus K$ 
such that }
\begin{enumerate}
	\item \emph{$\rot \vec f=0$ in $\R^3\setminus K$;}
	\item \emph{circulations of  $\vec f$ over curves $\gamma_j$ 
	have prescribed values $4\pi a_j, \, j=1,2,3$; }
	\item \emph{under the above two conditions field~$\vec f$ gives minimum to the energy norm $\int_{\R^3\setminus K} |\vec f(x)|^2\,dx$.}
\end{enumerate} 
Thus, the problem on equilibrium electric current comes to the problem on interpolation of a vector field (or of a differential form of degree $1$) in the exterior domain $\R^3\setminus K$ by its periods.

Author does not know anything about such a high-dimensional problem. In this paper we just study its two-dimensional analogue. 
Let us remark that, according to the high-dimensional motivation, we would like to avoid complex variable technique in the two-dimensional problem. But these methods will not be actually used, most of the proofs will be obtained by real-valued estimates.

\paragraph*{Organization of the paper.} \addcontentsline{toc}{subsection}{Organization of the paper} In section~\ref{section:statement_simple} we give the main definitions on the problem of interpolation by periods and state its simplest properties. The above-defined complete interpolation property is equivalent to the satisfying of following two conditions: first, operator $\Per\colon \lo\to \ell^2$ should be bounded (in this case we say that $\Omega$ possesses \emph{Bessel property}); second, operator $\Per\colon \lo\to \ell^2$ should be surjective (then we say that  $\Omega$ has \emph{interpolation property}). \emph{Bessel constant } $C_B(\Omega)$ is defined as the norm of operator $\Per\colon\lo\to\ell^2$ whereas \emph{interpolation constant}~$C_I(\Omega)$ is the minimal norm of right-inverse operator of~$\Per$. We also use \emph{weak Bessel constant} $\tilde C_B(\Omega)$ which is the supremum of norms of functionals $\Per_j\colon \lo\to \R$ over $j=1,2,\dots$. The property of complete interpolation  by periods means that $C_I(\Omega) <+\infty$ and $C_B(\Omega) <+\infty$. If $\tilde C_B(\Omega) < +\infty$ then we say that $\Omega$ has \emph{weak Bessel property}.

The problem on interpolation by periods  turns to be conformally invariant. Also, Bessel constant, weak Bessel constant and interpolation constant change in a controlled manner under application of a quasiconformal diffeomorphism to $\Omega$ (proposition~\ref{predl:conformal_quasiconformal_invariance}).

In section~\ref{section:prerequisites} we prove the basic results on interpolation by periods. In subsection~\ref{subsection:reproducing_kernels} we find \emph{reproducing kernels} of functionals  $\Per_j$, i.e. such forms $\kappa_j \in \lo$ that $\Per_j \omega = \langle{\omega, \kappa_j}\rangle_{\lo}$ for any $\omega\in \lo$. To find these $\kappa_j$ one should solve Dirichlet problem for Laplace operator in domain $\Omega$ with locally-constant boundary data, see theorem~\ref{th:reprodicing_existence}. It turns (proposition~\ref{predl:negative_prod}) that $\langle{\kappa_j, \kappa_{j'}}\rangle < 0$ for any different $j, j'\in\mathbb N$. This immediately implies (theorem~\ref{th:weak_strong_bessel}) that if $\tilde C_B(\Omega) < +\infty$ then $C_B(\Omega) < +\infty$; in other words, weak Bessel property of domain $\Omega$ implies its Bessel property. Bessel property of $\Omega$ is thus equivalent to the estimate $\sup\limits_{j\in\mathbb N}\Cap2\left(B_j, \R^2\setminus\left(\Omega\cup B_j\right)\right) <+\infty$; here $\Cap2(\cdot, \cdot)$ is capacity of a condenser with two plates, see subsection~\ref{subsection:capacity}. It seems that there does not exist simpler metric conditions on holes $B_j$ which describe Bessel property without additional conditions on interpolation or on uniform local finiteness property which we define below.

Complete interpolation  property means that system $\{\kappa_j\}_{j=1}^\infty$ is a \emph{Riesz basis} in its closed linear span. In subsection~\ref{subsec:amdissible_functions} we interpret this condition taking in account the explicit form of reproducing kernels. For given $a_1, a_2, \dots \in \R$ consider function $u\colon \R^2 \to \R$ for which $\Delta u=0$ in $\Omega$,  $u=a_j$ on $B_j$ for $j=1,2,\dots$ and $u=0$ in $\mathbb D^{(c)}=\R^2\setminus \mathbb D$. It turns (proposition~\ref{criterii_adm}) that the complete interpolation property is equivalent to the estimate 
$$
C_I^{-2}(\Omega)\cdot\sum\limits_{j=1}^\infty a_j^2 \le \int_{\R^2} |\nabla u|^2 \dl \le C_B^{2}(\Omega)\cdot\sum\limits_{j=1}^\infty a_j^2
$$
for any $a_1, a_2, \dots\in \R$; here $\lambda_2$ is two-dimensional Lebesgue measure in $\R^2$. Thus, estimates on closed forms are reduced to estimates on harmonic functions. Since harmonic functions minimize Dirichlet integral under given boundary Dirichlet data, it is convenient to work with class of Sobolev functions locally constant in $\R^2 \setminus \Omega$ and equal zero in $\mathbb D^{(c)}$. We call such functions \emph{admissible} (for domain $\Omega$), the space of admissible functions is everywhere denoted by $\mao$.

The main goal of section~\ref{section:partial_criterion} is to derive a partial 
criterion of complete interpolation property of domain $\Omega$. Namely,  necessary 
and sufficient conditions will be obtained in the case when one of the 
following assumptions is satisfied: either 
\begin{equation}
\label{eq:hyper_lower_intro}
	\tag{$\diamond$}
	\inf\limits_{j\in\mathbb N}\diam_H(B_j) > 0
\end{equation}
($\diam_H$ is the hyperbolic diameter in $\mathbb D$), or 
\begin{equation}
\label{eq:strong_sep_intro}
\tag{$\diamond\diamond$} 
\exists\, \eps>0\colon
\dist(B_j, B_{j'}) \ge \eps \cdot \max\{\diam B_j, \diam B_{j'}\}  
\mbox{\emph{ for any different }} j, j'\in\mathbb N.
\end{equation}
The last condition is called \emph{strong separatedness} of holes $B_j$. Theorem~\ref{th:predv_criterii} states that if all the holes $B_j, \, j=1,2,\dots, $ are disks and one of conditions~(\ref{eq:hyper_lower_intro}) or~(\ref{eq:strong_sep_intro}) is fulfilled, then complete interpolation property of domain $\Omega$ is equivalent to the both conditions~(\ref{eq:hyper_lower_intro}) and~(\ref{eq:strong_sep_intro}). At the same time, no one of  conditions~(\ref{eq:hyper_lower_intro}) and~(\ref{eq:strong_sep_intro}) is not necessary for complete interpolation property of $\Omega$ (example~\ref{example:inverse} in subsection~\ref{section:examples}). Here we use \emph{patch method} in order to construct admissible functions. Also, an analogue of Carleson measures arises in proofs (subsection \ref{subsection:interp_suff}).


Let $\mathcal B(0,1/2)\subset\mathbb D$ be a disk centered in the origin and of radius $1/2$. Examples from the subsection~\ref{section:examples} lead to the following question: \emph{is it possible to construct domains~$\Omega$ with $C_I(\Omega)$ and $C_B(\Omega)$ bounded from the above such that disk $\mathcal B(0,1/2)$ intersects a very big number of holes $B_j$?} The answer is negative. Let us say that domain $\Omega$ has \emph{uniform local finiteness property} if there exists such $N=N(\Omega) <+\infty$ that any disk in hyperbolic metric in $\mathbb D$ of hyperbolic diameter $1$ intersects no more than $N$ of holes~$B_j$.
Theorem~\ref{th:uniform_local_finiteness} from the section~\ref{section:uniform_local_finiteness} is the crucial point of our study. It states that complete interpolation property of domain $\Omega$ implies its uniform local finiteness property. In the proof we, taking in account only estimate $\sup\limits_{j\in\mathbb N}\Cap2\left(B_j, \R^2\setminus\left(\Omega\cup B_j\right)\right) <+\infty$, construct a function $u\in\mao$ for which values of $u$ on holes $B_j$ are not too small, but the integral $\int\limits_{\R^2} |\nabla u|^2 \dl$ is not very large. Such a function is obtained as a distance to $\mathbb D^{(c)}$ in the inner metric generated by conformal metric $\mathds 1_{\Omega} |dz|$. Lower estimates for $u|_{B_j}$ are derived from capacity estimates by subdivision of geodesics in the above-mentioned metric into appropriate arcs (lemmas~\ref{rho_estim},~\ref{lemma:chain_construct} and \ref{capacity_estim_lemma}).

Under uniform local finiteness condition it is possible to give  simple metric criteria  for Bessel property and for interpolation property of domain $\Omega$ (separately).  So, under uniform local finiteness condition, Bessel property turns to be equivalent to the condition of \emph{weak} separatedness of holes $B_j$: $\dist(B_j, B_{j'}) \ge \eps \cdot \min\{\diam B_j, \diam B_{j'}\}$ for any $j, j'=1,2,\dots, \, j\neq j',$ and some $\eps >0$ not depending on $j$ and $j'$. This is simple to prove by use of equivalence of Bessel and weak Bessel properties (theorem~\ref{th:bessel_sufficient}). We also give a scheme of a metric proof of this criterion which does not use abstract estimates of reproducing kernels but relies on existence of annular regions in $\Omega$  wide enough. A structure of partial order on the set of holes naturally arises in this proof: $B_j \succ B_k$, roughly speaking, if $\diam B_k \ll \diam B_j$ and $\dist(B_j, B_k) \ll \diam B_j$. From our viewpoint, the existence of such a structure is the most precise metric expression of the  essence of boundedness of the operator $\Per$.

Interpolation property of domain $\Omega$ turns to be equivalent to the above-mentioned \emph{capacity connectedness} (see definitions of graphs $G(\Omega, S)$ and $g(\Omega, s)$ in the beginning of subsection~\ref{subsec:interp_criterion}), if the family of holes in $\Omega$ is uniformly locally finite. This is proved in theorems~\ref{th:interp_sufficient} and \ref{th:interp_tree} in subsection~\ref{subsec:interp_criterion}. In the proofs we implement the connectedness of $G(\Omega, S)$ by constructing in $\R^2$ planar "roads"{} corresponding to the edges of $G(\Omega, S)$ and connecting holes $B_j$.  Interpolation property turns to be non-monotone by the set of holes in the following sense: this property can fail if we erase some holes in $\Omega$ (see remark after theorem \ref{th:interp_tree}).

Until subsection~\ref{subsection:nonsmooth}, we assume that every hole~$B_j$ is a closure of a domain with boundary smooth enough. In the subsection~\ref{subsection:nonsmooth} we generalize criterion of complete interpolation (theorem~\ref{th:full_interp_criterion}) to the case of arbitrary holes $B_j$, assuming just that any~$B_j$ is a continuum not separating the plane. For this aim we approximate such $B_j$'s by sets with smooth boundaries.

Finally, in subsection~\ref{subsec:open_final} we formulate several open questions generalizing the planar  problem on interpolation by periods. In particular, one of such generalizations delivers a non-trivial quasiconformal invariant of countably-connected Riemann surfaces --  \emph{the existence of an integer Riesz basis in  $L^2$-\emph{(}co\emph{)}homology space}.

Some technical proofs skipped in the main exposition are brought in the Appendix (section~\ref{section:appendix}).

All our estimates are constructive in the following sense: if we prove finiteness of some quantity then our arguments allow to obtain an explicit estimate on such a quantity.

\paragraph*{Some notation.} Complex plane $\mathbb C$ is always identified with $\R^2$. Symbol $\mathbb D$ denotes open unit disk in $\mathbb C$:  $\mathbb D = \{z\in \mathbb C \colon |z| < 1\}$. If $E\subset\mathbb C$, then $E^{(c)}$ denotes  $\mathbb C \setminus E$. Thus $\mathbb D^{(c)}=\{z\in\mathbb C \colon |z| \ge 1\}$. 

Symbol $\lambda_d$ denotes Lebesgue measure in  $\R^d$, usually $d=2$, sometimes $d=1$. "Almost everywhere"{} or "a.e."{} means "almost everywhere with respect to two-dimensional Lebesgue measure $\lambda_2$".  Measure  $\mathcal H^1$ is one-dimensional Hausdorff measure, that is length measure on various curves.

Symbol $\dist$ denotes \emph{Euclidean} distance on the plane $\mathbb C = \R^2$ whereas $\dist_H$ stands for distance in the \emph{hyperbolic} metric $\dfrac{2|dz|\hspace{5pt}}{(1-|z|)^2}$ in the unit disk $\mathbb D$. Symbols $\diam E$ and $\diam_H(E)$ denote Euclidean and hyperbolic diameters of a set $E \subset \mathbb C$ or $E\subset\mathbb D$ respectively.
If $z \in \mathbb C$ or $z\in\mathbb D$ and $r\ge 0$ then  $\mathcal B(z, r)$ and $\mathcal B_H(z, r)$ are open balls in Euclidean (respectively, hyperbolic) metric centered in $z$ and of radius $r$; $\bar{\mathcal B}(z, r)$ and $\bar{\mathcal B}_H(z, r)$ are the corresponding closed balls.

Symbol  $\card$ denotes the number of elements of a finite or a countable set; $\clos$ and~$\Innt$ denote closure and interior of a set in a topological space respectively.

If $\Omega\subset \R^2$ is an open set and $p\in [1, +\infty]$, then $W^{1,p}(\Omega)$ is the usual Sobolev space of scalar functions $u\colon \Omega\to \R$, integrable with exponent $p$ in $\Omega$ and having in~$\Omega$ generalized partial derivatives lying in $L^p(\Omega)$. The space $W^{1,p}_{\loc}(\Omega)$ is the class of functions belonging to $W^{1,p}(\tilde \Omega)$ for any strictly inner open subset $\tilde{\Omega} \subset\Omega$.

If $X, Y$ are some quantities then estimate $X\asymp Y$ means that $C_1 X \le Y \le C_2 X$ with some \emph{absolute} constants $C_1, C_2 >0$. In all the other cases of comparability it will be said what does the constants of comparability depend on.

\newpage

\section{Statement of interpolation problem. Simplest properties.}

\label{section:statement_simple}

\subsection{Periods of a $1$-form square-integrable  in a regular domain}
\subsectionmark{Periods of forms in regular domain}

\label{subsection:definitions}

Let us fix the class of domains under consideration.

Let $\{B_j\}_{j=1}^\infty$ be a countable family of compact subsets of open unit disk $\mathbb D \subset \R^2\simeq\mathbb C$. We will assume that:

\begin{enumerate}
\item sets $B_j$ are pairwise disjoint;
\item every $B_j$  is connected and does not separate the plane;
\item every $B_j$ is a closure of a domain with smooth boundary (this condition is technical, we will give it up in subsection~\ref{subsection:nonsmooth});
\item sets  $B_{j}$ accumulate only to the boundary $\partial\mathbb D$ of unit disk $\mathbb D$; in other words, if $j=1,2,\dots$ is fixed then $\inf\limits_{j'\neq j} \dist(B_{j'}, B_{j})$ is positive (symbol $\dist(E_1, E_2)$ denotes Euclidean distance between sets $E_1, E_2 \subset \R^2$).
\end{enumerate}
Consider the set $\Omega = \mathbb D \setminus \left(\cup_{j=1}^\infty B_j\right)$ -- this is a countably-connected domain; sets $B_j$, bounded connected components of $\Omega^{(c)}$, are called \emph{holes} in domain $\Omega$. The set $\Omega$ of such kind will be called a \emph{regular domain}. In case when all the $B_j$ are closed disks we say that~$\Omega$ is a \emph{regular domain with round holes}. Sometimes we also will consider domains with finite number of holes. If the holes in such a domain satisfy the above-listed conditions then we also call such a domain regular (investigation of such domains differs from studying of countably-connected domains by an obvious change of notation).

Recall that a (real) differential form of degree $1$ ($1$-form) on some set  $\Omega\subset \R^2$ is an object of the kind $\omega=\omega_x dx+\omega_y dy$ where  $\omega_x, \omega_y\colon \Omega\to \R$ are functions called \emph{components} of the form $\omega$. Put $|\omega|:= \sqrt{\omega_x^2+\omega_y^2}$. \emph{Hodge star} operator $*$ maps $1$-form~$\omega$ to a $1$-form $*\omega=\omega_x dy - \omega_y dx$; we have $**\omega=-\omega$. 

In what follows we consider forms with Lebesgue measurable components. A smooth $1$-form $\omega$ defined on an open set $\Omega$ is called \emph{closed} if $d\omega=0$ where $d$ is the usual differential of the form. If components of $\omega$ are non-smooth but locally integrable in $\Omega$ then let us say that form $\omega$ is closed if $d\omega=0$ in the sense of distributions. It means that 
\begin{equation}
\label{eq:gen_closed}
\int\limits_\Omega \omega \wedge d\eta = 0
\end{equation}
for any "test"{} function $\eta \in C_0^\infty(\Omega)$. Form $\omega$ with components in $L^1_{\loc}(\Omega)$ is called \emph{exact} if there exists a function $u \in W^{1,1}_{\loc}(\Omega)$ such that $du=\omega$ in the sense of distributions. If the components of $\omega$ belong to $L^2(\Omega)$ (this is our main case) then, in the case of non-smooth boundary $\partial\Omega$  the primitive function $u$ can be, in general, non-integrable in the whole set $\Omega$. On smooth forms, the classic notions of closedness and exactness coincide with their distributional analogues.  Further, $1$-form $\omega$ defined on an open set~$\Omega$ is called \emph{co-closed} (or \emph{co-exact}) if form $*\omega$ is closed (respectively, exact). We do not make use of codifferential operator $*d*$.

Let  $U\subset \mathbb R^2$ be an open subset,  $E\subset \mathbb R^2$ and   $\varphi\colon U\to E$ be a $C^1$-smooth mapping. Denote by $J_\varphi(z)$ the Jacobian of mapping $\varphi$ in $z\in U$, and let $D\varphi(z)$ be its differential in $z$; $|D\varphi(z)|$ is the norm of that differential as norm of a matrix $2\times 2$ with respect to Euclidean norm in $\R^2$. 

If a form $\omega$ is defined on $E$ then we may define a form $\varphi^\sharp\omega$ on $U$ by the standard pull-back of $\omega$ obtained by change of coordinates. If $1$-form  $\tilde\omega=\tilde\omega_x dx+\tilde\omega_y dy$ is understood as a row-vector with components ${\omega}_x$ and ${\omega}_y$ then $(\varphi^\sharp\omega)(z)=\omega(\varphi(z))\cdot D\varphi(z)$.  If $E$ is open then pull-back and differential commute. In particular, pull-back of a closed (or exact) form if closed (respectively, exact).

Let 
$$
L^{2,1}(\Omega)=
 \{\omega  \mbox{ --  $1$-form in }  \Omega,\, \|\omega\|^2_{L^{2,1}(\Omega)}:=\int\limits_\Omega |\omega|^2 \dl < +\infty\} 
$$ 
be the set of all differential forms of degree $1$ defined in $\Omega$ which are square-integrable over area measure (symbol $\lambda_2$ denotes two-dimensional Lebesgue measure in the plane). The space $L^{2,1}(\Omega)$ is endowed with the scalar product
\begin{multline}
\label{eq:scalar_define}
\langle\omega_1, \omega_2\rangle_{L^{2,1}(\Omega)} := \int\limits_\Omega \omega_1\wedge *\omega_2 = \int\limits_\Omega (\omega_{1x}\omega_{2x}+\omega_{1y}\omega_{2y})\dl, ~~ \omega_1, \omega_2 \in L^{2,1}(\Omega),\,\\ \omega_1 = \omega_{1x}dx+\omega_{1y}dy, \, \omega_2 = \omega_{2x}dx+\omega_{2y}dy,
\end{multline}
and is thus a Hilbert space.

Let 
$$
\lo = \{\omega \in L^{2,1}(\Omega) \colon d\omega = 0\}
$$
be the subspace consisting of all closed forms from $L^{2,1}(\Omega)$; here equation $d\omega=0$ is understood in the sense of distributions, see~(\ref{eq:gen_closed}). Subspace $\lo$ is closed in $L^{2,1}(\Omega)$. Scalar product in $\lo$ is inherited from $L^{2,1}(\Omega)$.

Let  $j=1,2,\dots$ be an index. Since sets $B_{j'}$ with $j'\neq j$ do not accumulate to $B_j$ then one can find a smooth closed oriented curve $\gamma_j\subset \Omega$ winding around $B_j$ once in the positive direction and not winding around other holes $B_{j'}$. 
For forms $\omega\in\lo$ smooth in $\Omega$ define \emph{period functional} $\Per_j$: put
$$
\Per_j(\omega) := \int_{\gamma_j} \omega.
$$ 
This functional will not change if we replace curve $\gamma_j$ by an another curve compactly homological (see~\cite{DeRham}) to $\gamma_j$  in $\Omega$.

\begin{predl}[see also \cite{Fuglede60}] 
	The functional $\Per_j$ has the unique continuous extension from closed and 
	smooth in $\Omega$   $L^2$-forms to the whole $\lo$.
\end{predl}

Density of smooth \emph{closed} forms in $\lo$ follows, for example, from Hodge decomposition (see subsection~\ref{Hodge_section}). The continuity of functional $\Per_j$ with respect to $L^2$-norm can be easily derived from the fact that curve $\gamma_j$ is contained in $\Omega$ with some neighbourhood. Integrals of a closed form over homological loops coincide while an appropriate average of the family of curves compactly homological to $\gamma_j$ in $\Omega$ will be absolutely continuous with respect to $\lambda_2$ and will have bounded density with respect to $\lambda_2$. 

Note that if $\omega\in \lo$ then integral $\int_{\gamma} \omega$ is well-defined and coincides with $\Per_j\omega$ for \emph{$2$-quasievery} smooth loop $\gamma\subset\Omega$ homological in $\Omega$ to $\gamma_j$ (see ~\cite{Fuglede60}).

Now define \emph{period operator} $\Per\colon \lo \to \ell^2$: put
$$
\Per \omega := \{\Per_j(\omega)\}_{j=1}^\infty, ~~ \omega \in \lo,
$$
(in what follows $\ell^2$ is the standard space of \emph{real} square-summable number sequences). In other words, operator $\Per$ maps a closed form  to the sequence of its periods around holes $B_j$. If, considering several regular domains, we will need to emphasize the dependence of operator $\Per$ on domain $\Omega$ then we will denote it by $\Per^{(\Omega)}$, and by $\Per_j^{(\Omega)}$ -- its components that are corresponding period functionals. Mainly we will work with one domain and such a notation will not be needed.

As in the case of smooth forms, by use of extension of a primitive over paths one can show that if $\omega\in\lo$ and $\Per\omega=0$ then there exists such a function $u \in W^{1,2}_{\loc}(\Omega)$ that $\omega=du$ in the sense of distributions.


\begin{define} 
\label{def:interp_problem}
Let $\Omega$ be a regular domain.
	\begin{enumerate}
		\item  Let us say that domain $\Omega$ possesses  \emph{Bessel property} if operator $\Per\colon \lo\to \ell^2$ is bounded. The norm of this operator will be called \emph{Bessel constant} of domain~$\Omega$ and denoted by $C_B(\Omega)$.
		
		\item If the quantity $\tilde C_B(\Omega) := \sup\limits_{j\in\mathbb N} \|\Per_j\|_{\left(\lo\right)^*}$ is finite then we  say that domain $\Omega$ possesses \emph{weak Bessel} property. The quantity $\tilde C_B(\Omega)$ will be called \emph{weak Bessel constant} of domain $\Omega$.
		
		\item Let us say that domain $\Omega$ has \emph{interpolation property} if for any sequence $a \in \ell^2$ there exists such a form $\omega\in \lo$ that $\Per \omega = a$ and $\|\omega\|_\lo \le C\cdot \|a\|_{\ell^2}$ with some $C < +\infty$ not depending on $a$. The least such $C$ will be called \emph{interpolation constant} of $\Omega$ and denoted by $C_I(\Omega)$. In other words, interpolation property means that operator $\Per\colon \lo\to \ell^2$ has a bounded right-inverse whereas $C_I(\Omega)$ is the least norm of right-inverse of $\Per$.
		
		\item If domain $\Omega$ has both Bessel and interpolation properties then we say that $\Omega$ has \emph{complete interpolation property}.
	\end{enumerate}
\end{define}

Bessel property, of course, implies weak Bessel property with $\tilde{C}_B(\Omega) \le C_B(\Omega)$. We will see below (theorem~\ref{th:weak_strong_bessel}) that weak Bessel property implies Bessel property with ${C}_B(\Omega) \le \sqrt 2 \tilde C_B(\Omega)$.

Every functional $\Per_j$ is bounded in $\lo$, thus $\Per$ is a closed operator. So, interpolation property is equivalent to surjectivity of operator $\Per \colon \lo \to \ell^2$.

The \emph{problem on interpolation by periods} is to describe regular domains $\Omega$ having Bessel, interpolation  and complete interpolation properties and also to estimate constants $C_B(\Omega)$ and $C_I(\Omega)$ via metric characteristics of mutual layout of sets $B_j$.

If reader wants to get immediate constructions of domains with complete interpolation property, then he (or she) may see examples \ref{example:annulus} and \ref{example:conformal_symmerty} from subsection \ref{section:examples}. Estimates from these examples are straightforward and use only the definition given above and also conformal invariance of our problem (proposition \ref{predl:conformal_quasiconformal_invariance}).

Note that we do not state the question on interpolation by normed periods 
$$
{\Per_j(\omega)}/{\|\Per_j\|_{\left(\lo\right)^*}},
$$
as it is usually done in classic problems on analytic interpolation. In these problems normalization of the sequence of values of analytic functions allows to account the speed of tendency of interpolation nodes to the boundary, whereas reproducing kernels of value functionals can often be found in an explicit form (or, at least, one can calculate asymptotic of norms of these functionals). We will see below that norms of functionals~$\Per_j$ (or, that is the same, norms of their reproducing kernels) cannot be found in a complete explicit manner, also, these norms may turn to be either arbitrarily small or arbitrarily large -- without any relation to the distance from holes $B_j$ to the unit circle $\partial \mathbb D$. One of our goals is to find conditions under which such norms are bounded from the above as well as from the below.

Now let us note that our problem can be stated in a more general manner. First, instead of a disk $\mathbb D$ we may take any other domain $V\subset \mathbb C$ or even Riemann surface (while sets $B_j$ would be holes in $V$). If $V$ is simply connected then, by Riemann uniformization theorem, we may to reduce such a problem to the case of a disk with holes (see also proposition~\ref{predl:conformal_quasiconformal_invariance} on conformal invariance of our problem); investigation of a uniformizing mapping may, in general, turn to be difficult (we note that techniques used for this aim is similar to the ones which we use below). We do not study the problem in such a generality. In subsection~\ref{section:cylinder} we, nevertheless, argue by a conformal mapping of an annulus in $\mathbb D$ onto an appropriate cylinder. This allows us to simplify calculations.

Second, we may calculate periods along another curves, not only winding around holes $B_j$ corresponding to them. Namely, let $\Omega\subset\mathbb C$  be an arbitrary domain (one may also consider Riemann surfaces) and $\gamma_j, \, j=1,2,\dots,$ be some smooth oriented closed curves in $\Omega$. The space $\lo$ in this case is defined as above. Every curve $\gamma_j$ lies in~ $\Omega$ together with some neighbourhood, and thus the functional $\Per_j\colon \lo \to \R$, $\Per_j(\omega) = \int_{\gamma_j}\omega,$ $\omega\in \lo,$ is well defined and continuous. The questions on Bessel and interpolation properties in this case are stated as above. We will say that \emph{a interpolation problem is stated in $\Omega$ for periods along curves $\gamma_j$}; constants $C_B(\Omega)=C_B(\Omega;\gamma_1, \gamma_2, \dots)$ and $C_I(\Omega)=C_I(\Omega;\gamma_1, \gamma_2, \dots)$ are defined as above. This problem does not change if we replace every curve $\gamma_j$ by another one equal to $\gamma_j$ in compact homology space $H_1^c(\Omega)$ (see~\cite{DeRham}).

The above described way to choose curves $\gamma_j$ corresponding to holes $B_j$ seems to be the most natural from the analytical viewpoint. Curves $\gamma_j$ are linearly independent and complete in the space of compact homologies (the last property means that if $\int_{\gamma_j}\omega =0$ for all $j=1,2,\dots$, then form $\omega\in\lo$ is exact). At the same time, if we start from the \emph{intrinsic} geometry of domain $\Omega$ (as a Riemann surface) then we can not make a canonical choice of a complete linearly independent system of compact homologies in~$\Omega$. In particular, if we choose curves $\gamma_j$ corresponding to holes $B_j$ in $\Omega$ and make an inversion $\varphi$ with the center in one of the holes $B_j$, then domain will turn inside out and curves $\varphi(\gamma_j)$ will not correspond to the holes in $\varphi(\Omega)$. We will meet such a situation considering examples in subsection~\ref{section:examples} (see example~\ref{example:inverse} and fig.~\ref{fig:Inverse} and also example~\ref{example:inversion_composition}). In the most cases we will still deal with the problem on interpolation by periods over curves corresponding to holes $B_j$. We assume this silently if the choice of curves $\gamma_j$ is not specified.

\subsection{(Quasi)conformal (quasi)invariance of the problem}
\label{subsection:qc}

Let us recall the definition of a quasiconformal mapping; we restrict ourselves by the smooth case. Let $\Omega, \tilde \Omega \subset \R^2$ be domains, $\varphi\colon \Omega \to \tilde{\Omega}$ be a  $C^\infty$-smooth diffeomorphism and $K \ge 1$ be a number. Mapping $\varphi$  is  called \emph{$K$-quasiconformal} if for any $z\in\Omega$ the following estimate is held:
\begin{equation}
\label{neq:qc_def}
|D\varphi(z)|^2 \le K\cdot J_\varphi(z).
\end{equation}
The least constant $K$ which possesses the estimate~(\ref{neq:qc_def}) is called \emph{the distortion coefficient} of $\varphi$ and is denoted by $K(\varphi)$. The definition of quasiconformality, in particular, implies that $J_\varphi(z)\ge 0$ for all $z\in\Omega$, that is, quasiconformal mapping preserves orientation. Also,  $1$-quasiconformal mapping is conformal. See, e.g.,~\cite{Ahlfors} for more information on quasiconformal mappings.

\begin{predl} 
	\label{predl:conformal_quasiconformal_invariance}
	Suppose that an interpolation problem is stated in domain $\Omega$ for periods along curves $\gamma_j,\,j=1,2,\dots$ \emph{(}see at the end of subsection~\ref{subsection:definitions}\emph{)} and that $\tilde\Omega$ is a  domain in $\mathbb R^2$. Suppose that $\varphi \colon \Omega\to\tilde{\Omega}$ is a orientation preserving diffeomorphism. State the interpolation problem in $\tilde{\Omega}$ for periods along curves $\varphi(\gamma_j), \, j=1,2,\dots$.
	\begin{enumerate}
		\item	If $\varphi$ is conformal then $C_B(\tilde\Omega) =C_B(\Omega)$ and $C_I(\tilde\Omega) =C_I(\Omega)$.
		
		\item If $\varphi$ is quasiconformal and $K=K(\varphi)$ is its distortion coefficient then 
	\begin{gather*}
	\sqrt{K^{-1}}\cdot C_B(\Omega) \le C_B(\tilde \Omega) \le  \sqrt{K}\cdot C_B(\Omega),\\
	\sqrt{K^{-1}}\cdot C_I(\Omega) \le C_I(\tilde \Omega) \le  \sqrt{K}\cdot C_I(\Omega),\\
	\sqrt{K^{-1}}\cdot \tilde C_B(\Omega) \le \tilde C_B(\tilde \Omega) \le  \sqrt{K}\cdot \tilde C_B(\Omega).
	\end{gather*}
	\end{enumerate}
\end{predl}

In other words,  problem on interpolation by periods \emph{is conformally  invariant} and \emph{quasiconformally quasiinvariant}. The proof is given in the Appendix.

\subsection{Forms minimizing $\|\omega\|_{\lo}$ under prescribed periods}

\label{subsection:minimizeers}

It is natural to reformulate the questions on Bessel and interpolation property of a regular domain $\Omega$ as follows: what is the minimal $L^2$-norm of a closed $1$-form with given periods? can this norm be estimated from the below (respectively, from the above) through sum of squares of  periods? We will need some simple properties of such norm minimizers.

\begin{predl}
	\label{predl:minimizer_exists}
	Let $\Omega$ be a regular domain and $a=\{a_j\}_{j=1}^\infty,\, a_j \in \R, \, j=1,2,\dots,$ be a sequence of scalars. Suppose that the set $Y=\{\omega\in\lo\colon\Per \omega=a\}$ is not empty. Then there exists a unique element $\omega_0(a)\in Y$ with minimal $\lo$-norm. 
\end{predl}

\noindent {\bf Proof.} The existence of a form $\omega_0(a)$ follows from the closedness of affine subspace $Y \subset\lo$. The uniqueness follows from the fact that a sphere in a Hilbert space does not contain a segment. $\blacksquare$

\medskip

\noindent {\bf Remark.} Bessel property of a  domain $\Omega$ is equivalent to the estimate $\|\omega_0(a)\|_{\lo}\ge C_B^{-1}(\Omega) \cdot \|a\|_{\ell^2}$ for any $a \in \Per(\lo)$. Interpolation property of a domain $\Omega$ is held if and only if, for any  $a \in \ell^2$, the form $\omega_0(a)$ is defined and $\|\omega_0(a)\|_{\lo}\le C_I(\Omega) \cdot \|a\|_{\ell^2}$.

\begin{predl}
\label{predl:minimizer_coexact}
	For any real sequence $a\in \Per(\lo)$ the form $\omega_0(a)$ defined in proposition~\ref{predl:minimizer_exists} is co-exact. Moreover, if $\omega \in \lo$ and $\Per\omega = 0$ then $\langle \omega_0(a), \omega\rangle_{\lo} =0$.
\end{predl}

The proof is given in the Appendix.

Recall the definition of the class of harmonic forms: 

\begin{define}
	A form $\omega\in L^{2,1}(\Omega)$ is called \emph{harmonic} if 
$$
\left\{
\begin{aligned}
&d\omega=0\\
&d*\omega=0
\end{aligned}
\right.
$$	
on $\Omega$ \emph{(}differential relations are understood in the sense of distributions\emph{)}. The space of all $1$-forms from $L^{2,1}(\Omega)$ harmonic in $\Omega$ is denoted by $\Harm^{2,1}(\Omega)$.
\end{define}

Forms $\omega_0(a)$ associated to real sequences $a$ in proposition~\ref{predl:minimizer_exists} are closed and co-exact and thus harmonic. The components of harmonic forms are harmonic functions, hence a form harmonic in $\Omega$ is smooth in $\Omega$.

\subsection{Interpolation by periods in  $\lo$ and in Bergman space}

\label{subsec:Bergman}

Let  $\Omega$ be a regular domain and 
$$
\mathscr A(\Omega)=\{f\colon \Omega \to \mathbb C\mid f \mbox{ is analytic in  } \Omega, \, \|f\|_{\mathscr A(\Omega)}^2=\int\limits_\Omega |f|^2\dl < +\infty\}
$$ 
be (unweighted) Bergman space in $\Omega$. To any function $f\in \mathscr A(\Omega)$ associate the sequence $\{\int_{\gamma_j}f(\zeta) \, d\zeta\}_{j=1}^\infty$ of its periods  along curves $\gamma_j$ chosen as above. The problem on interpolation of functions from Bergman space in $\Omega$ by their periods is stated analogously to the problem on forms (but period operator is considered as an operator from~$\mathscr A(\Omega)$ to the space $\ell^2_{\mathbb C}$ of \emph{complex} sequences). Let $C_{B, \mathscr A}(\Omega)$ and $C_{I, \mathscr A}(\Omega)$ be Bessel and interpolation constants for period operator in the space $\mathscr A(\Omega)$ defined in the same manner as in the problem on $1$-forms.

\begin{predl}
\label{predl:Bergman}
Let $\Omega$ be a regular domain.
\begin{enumerate}
	\item Domain $\Omega$ has Bessel property for $\lo$ if and only if it has Bessel property for $\mathscr A(\Omega)$. Also, $C_{B}(\Omega) = C_{B, \mathscr A}(\Omega)$.
	\item Domain $\Omega$ has interpolation property for  $\lo$ if and only if it has interpolation property for $\mathscr A(\Omega)$. Also, $C_{I}(\Omega) = C_{I, \mathscr A}(\Omega)$.
\end{enumerate}
\end{predl}

The proof is given in the Appendix.

Our problem is thus equivalent to the problem  in Bergman space. The multipliers technique used for investigation of interpolation problems in spaces of analytic functions seems not to be applicable in the case of periods as well as Blaschke products technique. Thus, in what follows, we will not make use of analytic functions.

\section{Preliminary results}
\label{section:prerequisites}

In subsection~\ref{Hodge_section} we consider Hodge decomposition of the space $L^{2,1}(\Omega)$. This will allow us to specify the subspace of  $L^{2,1}_c(\Omega)$ which contains reproducing kernels of functionals~$\Per_j$, $j=1,2,\dots$. We will find these kernels by solution of Dirichlet problem for Laplace operator in subsection~\ref{subsection:reproducing_kernels}. In subsection~\ref{subsection:Riesz} we apply Riesz conditions to the system of period reproducing kernels: interpolation problem will be reduced to estimates on functions from $\sob$ constant on sets $B_j, \, j=1,2,\dots$. In subsection~\ref{subsection:capacity} we derive from these estimates some capacity conditions necessary for the complete interpolation. In subsection~\ref{subsection:hyperbolic_metric}  we give metric estimates linking Euclidean and hyperbolic metrics to condenser capacities. Also, in subsection~\ref{subsection:hyperbolic_metric} we introduce conditions of weak (or strong) separatedness of holes; these conditions turn to be necessary (or sufficient) for Bessel property of $\Omega$.

\subsection{Hodge decomposition}
\label{Hodge_section}
For compact manifolds without boundary a Hodge decomposition with three components is usual. We will need four components of this decomposition. In this subsection all the orthogonal sums and orthogonal complements are understood in the sense of scalar product in $L^{2,1}(\Omega)$.

Let $\Omega$ be a regular domain (or even a Riemann manifold). The set 
$\Harm^{2,1}(\Omega)$ of {harmonic} forms from $L^{2,1}(\Omega)$ was introduced 
at the end of subsection~\ref{subsection:minimizeers}. It is well known that 
\begin{multline}
\hspace{-0.3cm}
\label{Hodge1}
L^{2,1}(\Omega) = \clos_{L^{2,1}(\Omega)}\left\{du\colon u\in C^\infty_0(\Omega)\right\} \oplus \clos_{L^{2,1}(\Omega)}\left\{*du\colon u\in C^\infty_0(\Omega)\right\}\oplus \Harm^{2,1}(\Omega) =\\= {\mathcal F}_1(\Omega)\oplus{\mathcal F}_2(\Omega)\oplus\Harm^{2,1}(\Omega).
\end{multline}
Indeed, Stokes' theorem immediately gives the orthogonality of the first two items of this decomposition. Further, if a form $\omega$ is orthogonal to any form from first two components then $\omega$ is closed and co-closed  in the sense of distributions and thus is harmonic (and smooth) in $\Omega$. 

\begin{predl}
\label{closed_orthogonal}
If $\omega\in L^{2,1}(\Omega)$ and $d\omega=0$, then $\omega\bot {\mathcal F}_2(\Omega)$ in $L^{2,1}(\Omega)$. In other words, 
$$
\lo = {\mathcal F}_1(\Omega) \oplus \Harm^{2,1}(\Omega).
$$
\end{predl}

\noindent {\bf Proof.} Equation $\int\limits_{\Omega} \langle\omega, *du\rangle \dl=0$ for all $u \in C_0^\infty(\Omega)$ is exactly the distributional form of equality $d\omega=0$.
$\blacksquare$ 

\medskip

Now let us refine decomposition~(\ref{Hodge1}). Put
$$
{\mathcal F}_3(\Omega) = \left\{\omega \in \Harm^{2,1}(\Omega)\colon \omega  \mbox{ is exact in } \Omega\right\}
$$
and, further,  ${\mathcal F}_4(\Omega) = \Harm^{2,1}(\Omega)\ominus {\mathcal F}_3(\Omega)$.
Then 
$$
\Harm^{2,1}(\Omega) = {\mathcal F}_3(\Omega)\oplus {\mathcal F}_4(\Omega),
$$
and finally  
\begin{equation}
\label{Hodge2}
L^{2,1}(\Omega) = {\mathcal F}_1(\Omega)\oplus {\mathcal F}_2\oplus {\mathcal F}_3(\Omega)\oplus {\mathcal F}_4(\Omega).
\end{equation}

\begin{predl}
\label{predl:Hodge_F4_presice}
Denote by ${\mathcal F}_5(\Omega)$ the set of all \emph{exact} forms from $L^{2,1}(\Omega)$. Then
$$
{\mathcal F}_4(\Omega) = \Harm^{2,1}(\Omega) \cap \left(L^{2,1}(\Omega)\ominus {\mathcal F}_5(\Omega)\right).
$$
In other words, forms from ${\mathcal F}_4(\Omega)$ are orthogonal not only to exact harmonic forms but also to all the square-integrable exact forms.
\end{predl}

The proof is given in the Appendix.

In the case when $\Omega$ is a compact manifold with boundary smooth enough, forms from $\omega\in {\mathcal F}_4(\Omega)$ can be described through their boundary data (namely, the normal component of $\omega$ on $\partial \Omega$ must vanish). Moreover, in this case the space ${\mathcal F}_4(\Omega)$ is finite-dimensional (see, e.g.,~\cite{GunterSchwarz}). If, to the opposite, the boundary $\partial \Omega$ of a domain $\Omega\subset\R^2$ is non-smooth then vanishing of normal components of forms makes no sense; moreover, the space ${\mathcal F}_4(\Omega)$ may turn to be infinite-dimensional. This will occur for infinite-connected regular domains.

The next proposition states that the fourth component in~(\ref{Hodge2}) is responsible for periods of closed forms.

\begin{predl} 
	\label{predl:4th_piece_proj}
	Let $\omega\in\lo$ be a form and $\omega_4$ be the projection of  $\omega$ to~$\mathcal F_4(\Omega)$ in~$L^{2,1}(\Omega)$. Then $\Per\omega = \Per\omega_4$. In particular, operator $\Per$ vanishes on $\mathcal F_1(\Omega)\oplus\mathcal F_3(\Omega)$.
\end{predl}

The proof is given in the Appendix.

\medskip

\noindent {\bf Remark.} If $\omega\in {\mathcal F}_4(\Omega), \, \omega\neq 0$, then $\Per_j(\omega)\neq 0$ for some $j=1,2,\dots$. Indeed, in the other case form $\omega$ would be exact whereas all the exact harmonic forms fall in ${\mathcal F}_3(\Omega)$ but not in~${\mathcal F}_4(\Omega)$.

\subsection{Period reproducing kernels}
\label{subsection:reproducing_kernels}

Let  $\Omega$ be a regular domain, $j=1,2,\dots$. We are going to find forms $\kappa_j\in\lo$ such that 
\begin{equation}
\label{eq:reproducing}
\langle\kappa_j, \omega\rangle_\lo = \Per_j(\omega)
\end{equation}
for any form $\omega\in\lo$. Such a form $\kappa_j$ is called \emph{reproducing kernel} of a functional~$\Per_j$ (in the space $\lo$). In the case of a domain (or even a Riemann surface) with smooth boundary these kernels were known before (see, e.g.,~\cite{Accola}). Proposition~\ref{predl:4th_piece_proj} implies that $\kappa_j\in {\mathcal F}_4(\Omega)$ (since the functional $\Per_j$ vanishes on ${\mathcal F}_1(\Omega)$ and on ${\mathcal F}_3(\Omega)$).

In what follows, $\sob$ is the Sobolev space consisting of functions $u\colon \mathbb D\to \R$ having in $\mathbb D$ square-integrable partial distributional derivatives of order $1$ and such that $u=0$ on $\partial\mathbb D$ in the sense of boundary values operator. The space $\sob$ is endowed with \emph{Dirichlet norm}:
$$
\|u\|_\sob :=  \left(\int\limits_{\mathbb D} |\nabla u|^2\dl\right)^{1/2}, ~~ u\in\sob.
$$
Any function from $\sob$ is assumed to be extended by zero into $\mathbb D^{(c)}$ to be a Sobolev function on the whole plane.

\begin{theorem}
	\label{th:reprodicing_existence}
	Let $\Omega$ be a regular domain and $j=1,2,\dots$ be an index. Then:
	\begin{enumerate}
	\item  There exists a function ${\mathfrak v}_j={\mathfrak v}_j(\Omega)\in W^{1,2}(\mathbb D)$ such that:
	   \begin{equation}
	   \label{reproducing_eq}
	   \left\{
	   \begin{aligned}
	   & \Delta \mathfrak v_j = 0 \mbox{ in } \Omega
	   ;\\
	   & \mathfrak v_j = 1\mbox{ almost everywhere on } B_j;\\
	   & \mathfrak v_j=0 \mbox{ almost everywhere on } B_{j'}, \, j'\neq j;\\ 
	   & \mathfrak v_j = 0 \mbox{ on } \partial \mathbb D \mbox{ in the sense of boundary values operator.} 
	   \end{aligned}\right.
       \end{equation}
       The function $\mathfrak v_j$ minimizes Dirichlet integral $\int_{\mathbb D} |\nabla v|^2\,\dl$ over all functions $v\in \sob$ with the same values on $\Omega^{(c)}$.
    \item  If $v\in W^{1,2}(\mathbb C), \, v=0$ almost everywhere in  $\Omega^{(c)},\, \Delta v=0$  in $\Omega$, then $v=0$.
    \item Function $\mathfrak v_j$ satisfying \emph{(\ref{reproducing_eq})} is unique. 
	\item  Reproducing kernel of the functional $\Per_j$ on $\lo$ is the form $\kappa_j=-(*d\mathfrak v_j)$.
	\end{enumerate}
\end{theorem}

The proof is given in the Appendix.

Note that if holes $B_j$ are, for example, slits then we  cannot state the Dirichlet problem by equality of functions to $0$ or $1$ almost everywhere on holes. In this case the boundary problem should be stated in the class of precised functions. Then the equality to $0$ or $1$ will be up to zero capacity, that is \emph{quasieverywhere}. We then have to require positivity of capacities of all the holes $B_{j'}$. This will be the case if any $B_{j'}$ is connected and consists of more than one point.

Our problem can also be stated in terms of vector fields ($1$-forms are naturally associated to vector fields). Periods of forms then correspond to circulations of fields over given curves. In the language of vector fields reproducing kernel $-(*d\mathfrak v_j)$ corresponds to the field $(\nabla \mathfrak v_j)^\bot$, that is the gradient of function $\mathfrak v_j$ turned over $\pi/2$ counter-clockwise.

\medskip

\noindent {\bf Remark.} By the remark after proposition~\ref{predl:4th_piece_proj} the system $\{\kappa_j\}_{j=1}^\infty$ will be complete in the space~${\mathcal F}_4(\Omega)$. 

\medskip

In the proof of the theorem~\ref{th:reprodicing_existence} we use the following inequality (we will need it in some other situations):

\begin{theorem}
	\label{th:Hardy_neq}
	There exists a constant $\mathfrak c>0$ such that 
	$$
	\int\limits_{\mathbb D} \frac{u^2(z)\dl(z)}{(1-|z|)^2} \le \mathfrak c\cdot \int\limits_{\mathbb D} |\nabla u|^2 \dl
	$$
	for any function $u \in \sob$. In other words, the embedding $\sob\hookrightarrow L^2(\mu_0)$ where  $\mu_0=\dfrac{\lambda_2}{(1-|z|)^2}$ is continuous.
\end{theorem}

Indeed, it is enough to estimate $\displaystyle\int\limits_{|z|\ge 1/2} \frac{u^2(z)\dl(z)}{(1-|z|)^2}$. But this is easy to do by an application of well-known Hardy inequality (see, e.g.,~\cite{HardyLittlewoodPolya})
$$
\int\limits_0^\infty \left(\frac1x\int\limits_0^x f(y)\,dy\right)^2\,dx \le 4 \int\limits_0^\infty f^2(y)\,dy, ~~ f\in L^2[0,+\infty),
$$
to the function $f(r) = \dfrac{\partial}{\partial r} u\left((1-r) e^{i\theta}\right)$ on the segment $[0,1/2]$ and integrating the obtained estimate over $\theta \in [0,2\pi)$. 

Note, by the way, that the full mass of measure $\mu_0$ is infinite, in particular, estimate from the theorem~\ref{th:Hardy_neq} is not held for $u\equiv 1$ (even if we replace the right-hand side of this inequality by the whole Sobolev norm). In other words, zero boundary data of function $u$ on $\partial\mathbb D$ is essential for validity of this estimate.


The following surprising observation will easily lead us to the result of theorem \ref{th:weak_strong_bessel} (see below):

\begin{predl}
	\label{predl:negative_prod}
	Let $\Omega$ be a regular domain and $j, j'=1,2,\dots, \, j\neq j'$. Then $\langle\kappa_j, \kappa_{j'}\rangle_{\lo} <0$ where $\kappa_j, \,\kappa_{j'}$ are period reproducing kernels found in theorem~\ref{th:reprodicing_existence}.
\end{predl}

The proof is given in the Appendix.

\subsection{Riesz property}

\label{subsection:Riesz}

So, the period reproducing kernels $\kappa_j \in \mathcal F_4(\Omega), \, j=1,2,\dots,$ are found. We are going to apply Riesz bases theory to these kernels. Let us start with the following 

\begin{predl}
	System $\{\kappa_j\}_{j=1}^\infty$ has a  \emph{biorthogonally adjoint system} in $\mathcal F_4(\Omega)$, that is such a sequence of forms $\mathfrak w_j \in \mathcal F_4(\Omega),\, j=1,2,\dots,$ that  $\langle \mathfrak w_j, \kappa_j\rangle_\lo =1$ and $\langle \mathfrak w_j, \kappa_{j'}\rangle_\lo =0$ for $j\neq j'$. 
\end{predl}

\noindent {\bf Proof.} 
Pick $z_j=x_j + i y_j \in \Innt B_j$. Then the form $$
\omega_j(x+iy) := \frac{(x-x_j)dy-(y-y_j)dx}{2\pi\left((x-x_j)^2+(y-y_j)^2\right)}, ~~ x+iy\in\Omega,
$$
belongs to $\lo$, $\Per_j\omega_j=1, \, \Per_j \omega_{j'}=0$ for $j'\neq j$. For $\mathfrak w_j$ we now may take the projection of $\omega_j$ to $\mathcal F_4(\Omega)$.
$\blacksquare$

\medskip

In notation of~\cite{Bari},  Bessel property of $\Omega$ means Bessel property of system $\{\mathfrak w_j\}_{j=1}^\infty \subset \mathcal F_4(\Omega)$, whereas interpolation property of $\Omega$ means Hilbert property of this system. Applying of results from~\cite{Bari}, it is not hard to prove the following

\begin{predl} \label{pred:Riesz_Criterii} Let $\Omega$ be a regular domain and $\mathfrak v_j,\, j=1,2,\dots,$ be the functions obtained by solution of a Dirichlet problem \emph{(\ref{reproducing_eq})}.
\begin{enumerate}
	\item Bessel property of $\Omega$ is equivalent to the following: for any sequence $\{a_j\}_{j=1}^\infty\in \ell^2$ the series $\sum_{j=1}^\infty a_j \mathfrak v_j$ converges strongly in $\sob$ and $$\left\|\sum_{j=1}^\infty a_j \mathfrak v_j\right\|_\sob^2 \le C_B^2(\Omega) \cdot \sum_{j=1}^\infty a_j^2.$$
	\item Interpolation property of $\Omega$ is equivalent to the estimate 
	\begin{equation}
	\label{estim:kernels_neq_interp}
	\left\|\sum_{j=1}^\infty a_j \mathfrak v_j\right\|_\sob^2 \ge C_I^{-2}(\Omega) \cdot \sum_{j=1}^\infty a_j^2
	\end{equation}
	for any finitely supported sequence $\{a_j\}_{j=1}^\infty$. 
	In case of Bessel property the series $\sum_{j=1}^\infty a_j \mathfrak v_j$ converges strongly in $\sob$ for any $\{a_j\}_{j=1}^\infty\in\ell^2$, and  estimate \emph{(\ref{estim:kernels_neq_interp})} stays true for all such  $\{a_j\}_{j=1}^\infty$.
\end{enumerate}
\end{predl}

Thus, complete interpolation property of $\Omega$ is equivalent to the following: for any $a=\{a_j\}_{j=1}^\infty\in \ell^2$ series $u=\sum_{j=1}^\infty a_j \mathfrak v_j$ converges strongly in $\sob$ and its sum~$u$ admits the estimate 
\begin{equation}
\label{eq:Riesz_kappa}
C_I^{-1}(\Omega)\cdot\|a\|_{\ell^2} \le \|u\|_\sob \le C_B(\Omega)\cdot \|a\|_{\ell^2}.
\end{equation}

\medskip

\noindent {\bf Remark.} The system $\{\kappa_j\}_{j=1}^\infty$ belongs to $\mathcal F_4(\Omega)$ and, according to the remark at the end of subsection~\ref{Hodge_section}, is complete in $\mathcal F_4(\Omega)$. Thus the estimate~(\ref{eq:Riesz_kappa}) means that in case of complete interpolation the system $\{\kappa_j\}_{j=1}^\infty$ must be a \emph{Riesz basis} (see~\cite{Bari},~\cite{Nikolskiy}) in~$\mathcal F_4(\Omega)$.  

By proposition~\ref{predl:Hodge_F4_presice}, term $\mathcal F_4(\Omega)$ in Hodge decomposition~(\ref{Hodge2}) can be identified with the space $H_{L^2}^1(\Omega)$ of Hilbert cohomologies (or $L^2$-cohomologies) in $\Omega$.  Let us define this space in the following manner:
$$
H^1_{L^2}(\Omega) := \{\omega\in L^{2,1}(\Omega)\colon d\omega=0\}/\{\omega\in L^{2,1}(\Omega)\colon \omega=du \mbox{ for some } u\in W^{1,2}_{\loc}(\Omega)\}.
$$
Let us also define the space of \emph{Hilbert homologies} $H_{1,L^2}(\Omega)$ as a space dual to $H^1_{L^2}(\Omega)$. 
Although $H_{1,L^2}(\Omega)$ is canonically isomorphic to $H^1_{L^2}(\Omega)$, from viewpoint of homology theory these two spaces should be understood in a different ways.

For $j=1,2,\dots$ the curve $\gamma_j$ linked to the hole $B_j$ in domain $\Omega$ delivers a continuous functional on $H_{L^2}^1(\Omega)$ and thus can be regarded as an element of $H_{1,L^2}(\Omega)$. Thus,  \emph{$\Omega$ has complete interpolation property if and only if curves $\gamma_j, \, j=1,2,\dots,$  form a Riesz basis in the space $H_{1,L^2}(\Omega)$}.

\subsection{Admissible functions}

\label{subsec:amdissible_functions}

Finite linear combinations of kernels $\kappa_j$ are locally constant on $\Omega^{(c)}$. Let us introduce the following function class. 

\begin{define} A function $u\in \sob$ is called \emph{admissible for domain $\Omega$} if, for any $j=1,2,\dots$, function $u$ is constant on $B_j$ almost everywhere.
The class of all functions admissible for $\Omega$ is denoted by $\mao$.
\end{define}

The set $\mao$ endowed with the norm $\|u\|_{\sob}$ is a Hilbert space, functionals $u\mapsto u|_{B_j},\,j=1,2,\dots,$ are continuous on $\mao$.

Let  $\{a_j\}_{j=1}^\infty$ be some sequence of real scalars. Among all $u\in\mao$ with prescribed values $u|_{B_j}\equiv a_j$ the minimizer of norm $\|u\|_\sob$ is a function harmonic in~$\Omega$ and equal $a_j$ on $ B_j$. For finitely supported sequences $\{a_j\}_{j=1}^\infty$ this function is exactly~$\sum_{j=1}^\infty a_j \mathfrak v_j$. The following proposition will be our main tool to study interpolation:

\begin{predl}
	\label{criterii_adm}
	Let $\Omega$ be a regular domain.
\begin{enumerate}
	\item Domain $\Omega$ has Bessel property if and only if for any sequence $a=\{a_j\}_{j=1}^\infty\in \ell^2$ there exists such a function $u\in\mao$ that $u|_{B_j} = a_j$ and $\|u\|_\sob^2\le C_B^2(\Omega) \cdot \|a\|_{\ell^2}^2$. It is enough to check this condition only for finitely supported sequences $\{a_j\}_{j=1}^\infty$.
	
	\item Domain $\Omega$ possesses interpolation property if and only if for any function $u\in\mao$ the following estimate is true:
	\begin{equation}
	\label{est:interp_criterion}
	\|u\|_{\sob}^2\ge C_I^{-2}(\Omega)\cdot\sum\limits_{j=1}^\infty (u|_{B_j})^2.
	\end{equation}
\end{enumerate}
\end{predl}

{\noindent \bf Remark.}  An admissible function $u$ can be understood as a $0$-form whose integrals over any two $0$-simplices in $\Omega^{(c)}$ homological in $\Omega^{(c)}$ coincide. This links estimates from proposition~\ref{criterii_adm} to Alexander-Pontryagin duality known from homology theory.

\medskip

\noindent {\bf Proof of proposition~\ref{criterii_adm}.} 
The only non-trivial assertion\footnote{ The reader may to be amazed by the length of the following proof of this fact which, at a glance, seems to be obvious. Author, nevertheless, could not find a simpler and more abstract argument.} here is that estimate~(\ref{estim:kernels_neq_interp}) for finite linear combinations of reproducing kernels implies inequality~(\ref{est:interp_criterion}) for an arbitrary function $u\in \mao$. (All the other statements are easily obtained by use of the fact that harmonic functions minimize Dirichlet integral if boundary values are given, and also by uniqueness of solution of Dirichlet problem with zero boundary data, see theorem~\ref{th:reprodicing_existence}.)

We may assume that function $u$ for which we want to prove inequality~(\ref{est:interp_criterion}) is harmonic in $\Omega$. Indeed, let $u\in\mao$. The usual variational argument shows that, first, there exists a function $u_h\in\sob$ for which $u_h=u$ almost everywhere in $\mathbb D\setminus\Omega$ and $\Delta u_h=0$ in  $\Omega$, and, second, that $\|u_h\|_{\sob}\le \|u\|_{\sob}$.
But then estimate~(\ref{est:interp_criterion}) for function $u$ will follow from this estimate for ${u_h}$.

So, assume that $\Delta u=0$ in $\Omega$. Function $u$ is continuous on any compact subset in~$\mathbb D$ and its partial derivatives of all orders are continuous in $\Omega$ up to all $\partial B_j, \, j=1,2,\dots$. This follows from smoothness of boundaries $\partial B_j$ (see~\cite{LU}). Define $u$ to be zero outside of $\mathbb D$. Let us assume that $\nabla u=0$ \emph{everywhere} in $\mathbb C\setminus\Omega$. This will allow to avoid problems with measurability in the construction below.

On $\mathbb C$, consider  inner metric generated by metric $|\nabla u||dz|$ (cf. with proof of theorem~\ref{th:uniform_local_finiteness} below). 
Namely,  for any Lipschitz curve $\Gamma \colon [0,T] \to \mathbb C \, (T\in\R)$ the integral $$L(\Gamma):= \int\limits_\Gamma |\nabla u| \, d\mathcal{H}^1=\int\limits_0^T  \left|\nabla u\left(\Gamma(t)\right)\right| \cdot |\Gamma'(t)|\,dt\in[0, +\infty]$$ is well-defined (the expression under the last integral is measurable since since function $|\nabla u(\Gamma)|$ vanishes outside of open set $\Gamma^{-1}(\Omega)\subset[0,T]$ and is continuous on this set).
For $z_1, z_2 \in \mathbb C$, put  $\rho(z_1, z_2) := \inf L(\Gamma)$ where  $\inf$ is taken over all Lipschitz curves in~$\mathbb C$ joining $z_1$ and $z_2$.
Then $\rho$ is a degenerated metric in~$\mathbb C$ in which holes $B_j$ collapse to points (because any two points $z_1, z_2 \in \Omega^{(c)}$ lying in the same hole~$B_j$ can be joined by a Lipschitz curve $\Gamma\subset B_j$, we have $L(\Gamma)=0$ for such a curve); the exterior $\mathbb D^{(c)}$ of unit disk collapses to point as well. Triangle inequality for metric $\rho$ is nevertheless held.

Pick some $\nu\in\mathbb N$. Put 
\begin{equation*}
u_\nu(z) := \dist_\rho\left(z,\, \mathbb D^{(c)} \cup \bigcup\limits_{j=\nu}^\infty B_j\right), ~~ z\in \mathbb C.
\end{equation*}
In other words,  $u_\nu(z)=\inf L(\Gamma)$, where $\inf$ is taken over all Lipschitz curves joining $z$ with a point in $\mathbb D^{(c)}$ or with a point in one of the sets $B_j, \, j=\nu, \nu+1, \dots$. If  $z_1, z_2\in\Omega$ then we may join these points by a segment $[z_1, z_2]$  of a line and conclude that 
\begin{equation}
\label{estim:Lipschitz_metric}
|u_\nu(z_1)-u_\nu(z_2)|\le\rho(z_1, z_2) \le L([z_1, z_2]) \le |z_1-z_2|\cdot \sup\limits_{[z_1, z_2]}|\nabla u(z)|.
\end{equation}
This inequality and  local boundedness of $\nabla u$ in $\mathbb D$ imply that function~$u_\nu$ is locally Lipschitz in $\mathbb D$. If $z_1$ belongs to $\Omega$ or to one of sets $\Innt B_j, \, j=1,2,\dots$,  then gradient~$\nabla u$ is continuous in $z_1$; this and~(\ref{estim:Lipschitz_metric}) imply that 
$$
\varlimsup\limits_{z_2 \to z_1} \frac{|u_\nu(z_1) - u_\nu(z_2)|}{|z_1-z_2|} \le |\nabla u(z_1)|.
$$
Thus
\begin{equation}
\label{estim:Lipschitz_interp2}
|\nabla u_\nu(z)| \le |\nabla u(z)|
\end{equation} 
for almost all $z \in \mathbb D$. 
In particular,  function $u_\nu$ is constant on any $B_j$, $j=1,2,\dots$. Also,~(\ref{estim:Lipschitz_interp2}) implies that $u_\nu \in W^{1,2}(\mathbb D)$ (since $|\nabla u|\in L^2(\mathbb D)$).

Next, pick a number $r_0\in (0,1)$. For $\theta\in[0,2\pi]$ join the point $r_0e^{i\theta}$ with $e^{i\theta}$ by a segment of line. We then have 
$$
u_\nu(r_0e^{i\theta}) \le \int\limits_{[r_0e^{i\theta}, e^{i\theta}]} |\nabla u|\,d\mathcal H^1=
\int\limits_{r_0}^1 |\nabla u(re^{i\theta})|\,dr.
$$
Integrating over $\theta$ and taking in account that $\dl(re^{i\theta}) = r\,dr\,d\theta\ge r_0\,dr\,d\theta$ we get
$$
\int\limits_{\partial\mathcal B(0,r_0)} u_\nu\,d\mathcal H^1 
\le \int\limits_{|z|\in[r_0,1]} |\nabla u(z)|\dl(z).$$
The quantity in the right-hand side tends to zero with $r_0\to1$. From this it follows that the trace of function $u_\nu$ on $\partial\mathbb D$ (in the sense of boundary values operator) equals to zero. 

So, $u_\nu \in \mao$ while the sequence $\{u_\nu|_{B_j}\}_{j=1}^\infty$ is finitely supported (because the function $u_\nu$ by its definition vanishes on $B_\nu \cup B_{\nu+1}\cup\dots$). Let us prove that 
\begin{equation}
\label{est:adm_interpolation}
C_I^2(\Omega)\cdot\|u\|^2_\sob \ge C_I^2(\Omega)\cdot\|u_\nu\|^2_\sob \ge \sum\limits_{j=1}^\infty \left(u_\nu|_{B_j}\right)^2
\end{equation}
Indeed, the first inequality follows from the estimate~(\ref{estim:Lipschitz_interp2}) which we already proved; to prove the second inequality consider function $\tilde{u}_\nu\in \mao$ which is equal to $u_\nu$ on any set $B_j$ and minimizing the norm $\|\tilde{u}_\nu\|_\sob$. Then $\Delta\tilde u_\nu=0$ in $\Omega$ by the usual variational argument, $\tilde u_\nu = u_\nu$ on $B_j, \, j=1,2,\dots,$ and $\tilde u_\nu\in \mao$. By theorem~\ref{th:reprodicing_existence} the function satisfying these conditions is   unique. But the finite linear combination $\tilde u_\nu'=\sum_{j=1}^\infty (u_\nu|_{B_j})\cdot \mathfrak v_j$ also satisfies to all the above-mentioned conditions, and thus $\tilde u_\nu=\tilde u_\nu'$. Application of~(\ref{estim:kernels_neq_interp}) to $\tilde u_\nu'$ gives us 
$$
\|u_\nu\|_\sob^2 \ge \|\tilde u_\nu\|_\sob^2 \ge C_I^{-2}(\Omega)\cdot \sum\limits_{j=1}^{\infty} (u_\nu|_{B_j})^2.
$$
Inequality~(\ref{est:adm_interpolation}) is established.
Now we are going to show that values $u_\nu|_{B_j}$ are not very small with respect to $u|_{B_j}$.

Fix a number $n\in \mathbb N$ and pick some $\nu>n$; the choice of $\nu$ will be specified in the below. 
Suppose that $T\in\R$ and $\Gamma\colon[0, T]\to {\mathbb C}$ is a Lipschitz curve joining some point $z_0\in B_n$ with a point on $\partial\mathbb D$ or with a point on one of the sets $B_j, \, j=\nu, \nu+1, \dots$. We have to estimate $L(\Gamma)$. We may assume that $\Gamma([0,T))\subset\mathbb D$. Put $\delta = 1-|\Gamma(T)|$. If our curve ends on $\partial \mathbb D$ then  $\delta=0$. If not, then we may take $\nu$ big enough such that $\delta$ will be as small as we want, because in this case $\Gamma$ ends at one of the holes $B_{\nu},B_{\nu+1}, \dots$ while holes accumulate only to $\partial\mathbb D$. In particular, we may pick $\nu$ big enough in order to have $\delta < \dist(B_n, \partial\mathbb D)$. For any $t<T$ the arc $\Gamma([0,t])$ lies strictly inside $\mathbb D$, then $u$ is piecewise-smooth on this arc. Then we may apply Newton-Leibniz formula to get
$$
u(\Gamma(t))-u(\Gamma(0)) = \int\limits_{\Gamma([0,t])}\left(\nabla u\right)_\tau \,d\mathcal H^1,
$$
where $\left(\nabla u\right)_\tau$ is tangent component of gradient along $\Gamma$. Thus for any $t<T$ we get estimate 
\begin{equation}
\label{est:NewtonLeibnitz}
L(\Gamma) \ge \left|(u|_{B_n}) - u(\Gamma(t))\right|.
\end{equation}

Put $a=\inf\{|u(z)|\colon z\in \Gamma([0,T])\}$. Denote by $\zeta$ the point on $\partial\mathbb D$ closest to $\Gamma(T)$. If $r\in(\delta, \dist(B_n, \partial\mathbb D))$ then the arc $\beta_r=\partial\mathcal B(\zeta, r) \cap \mathbb D$ lying in $\mathbb D$ separates $\Gamma(T)$ from~$B_n$ and thus has to intersect $\Gamma$. For almost every $r\in(\delta, \dist(B_n, \partial\mathbb D))$ the function $u$ is absolutely continuous on $\beta_r$, takes zero values at the endpoints of this arc, also $|u(z)|\ge a$ in some point $z\in\beta_r$. Thus
$$
\int\limits_{\beta_r}|\nabla u|^2\,d\mathcal H^1 \ge \frac a{\pi r}
$$
for almost every $r$. Integration over $r$ gives 
$$
\int\limits_{\mathbb D} |\nabla u|^2 \dl \ge \frac{a}{\pi}\cdot\left(\log\dist(B_n, \partial\mathbb D)-\log\delta\right).
$$
If $n$ is fixed, we can make $\delta$ as small as we need by picking $\nu>n$ large enough, and from the last estimate it follows that in this case the number $a$ also has to become small (not depending on the choice of $\Gamma$). Estimate~(\ref{est:NewtonLeibnitz}) now implies that $L(\Gamma) \ge \left |(u|_{B_j})\right| - a$. Since $u_\nu|_{B_n} = \inf L(\Gamma)$ over $\Gamma$ joining $B_n$ with $\partial\mathbb D$ or with one of the holes $B_\nu, B_{\nu+1}, \dots$, we conclude that $\varliminf\limits_{\nu\to\infty} \left|(u_\nu|_{B_n})\right| \ge \left|(u|_{B_n})\right|$.

Now fix some $m\in\mathbb N$. Inequality~(\ref{est:adm_interpolation})  implies, in particular, that 
$$
\|u\|^2_\sob \ge C_I^{-2}(\Omega)\cdot\sum\limits_{n=1}^m \left(u_\nu|_{B_n}\right)^2.
$$
Passing to $\varliminf$ with $\nu\to\infty$ we conclude that
$$
\|u\|^2_\sob \ge C_I^{-2}(\Omega)\cdot\sum\limits_{n=1}^m \left(u|_{B_n}\right)^2.
$$
Now it remains to pass to the limit over $m\to\infty$ and obtain the desired inequality
$$
\|u\|^2_\sob \ge C_I^{-2}(\Omega)\cdot\sum\limits_{n=1}^\infty \left(u|_{B_n}\right)^2
$$
for any  function $u\in\mao$ such that $\Delta u=0$ in $\Omega$, and, since that, for arbitrary $u\in\mao$. Proof is finished.
$\blacksquare$

\begin{sled}
\label{sled:monotonicity_domain}
Let $\Omega_1, \Omega_2$ be regular domains, $\{B_j^{(1)}\}_{j=1}^\infty, \{B_j^{(2)}\}_{j=1}^\infty$ be bounded connected components of their complements to $\mathbb R^2$. Suppose that $B_j^{(1)}\subset B_j^{(2)}$ for any $j=1,2,\dots$. Then $C_B(\Omega_2) \ge C_B(\Omega_1)$ and $C_I(\Omega_2) \le C_I(\Omega_1)$.
\end{sled}

\noindent {\bf Proof.} We may use the proposition just proved by comparison of classes of functions admissible for $\Omega_1$ and $\Omega_2$. We can, nevertheless, get a straightforward proof relying just on the definitions of Bessel and interpolation constants. $\blacksquare$ 

\medskip

\subsection{Capacity conditions}

\label{subsection:capacity}

Now let us relate our problem to $L^2$-capacity. A \emph{condenser} is a pair of $(E_1, E_2)$ of two sets on the plane; sets $E_1$ and $E_2$ are called \emph{plates} of this condenser. The following definition of capacity is a rude one but is enough for our purposes.

\begin{define}
\label{def:capacity}
	Let $E_1, E_2\subset \mathbb C$ be arbitrary sets. \emph{Capacity of condenser} formed by plates $E_1$ and $E_2$ is the quantity
	\begin{equation}
	\label{eq:capacity_def1}
		\Cap2(E_1, E_2) := \inf\left\{\int\limits_{\mathbb C}|\nabla u|^2\dl\colon u\in W^{1,2}_{\loc}(\mathbb C),\, u\ge 1 \mbox{ a.e.\hspace{0.61mm}on } E_1,\, u\le 0 \mbox{ a.e.\hspace{0.61mm}on } E_2\right\}.
	\end{equation}
\end{define}

The simplest properties of capacities are given in the following 

\begin{predl} 
\label{predl:capacity_simple}
Let $E_1, E_1', E_2, E_2', E_3 \subset\mathbb C$ be arbitrary sets.
	\begin{enumerate} 
		
		\item 	Infimum in \emph{(\ref{eq:capacity_def1})} will not change if we, in addition, substitute function from this $\inf$ to restriction $u(z) \in[0,1]$  for almost every $z\in\mathbb C$.
		
		\item Capacity is symmetric, that is, $\Cap2(E_1, E_2)=\Cap2(E_2, E_1)$. If $E_1\subset E_1'$, $E_2 \subset E_2'$ then $\Cap2(E_1, E_2)\le \Cap2(E_1', E_2')$. 
		
		\item Capacity is semiadditive: $\Cap2(E_1, E_2\cup E_3) \le \Cap2(E_1, E_2)+\Cap2(E_1, E_3)$.
	\end{enumerate}
\end{predl}

The proof is given in the Appendix.

We will also need conformal invariance of capacity. 

\begin{predl}
\label{predl:capacity_invariance}
	Suppose that $E_1, E_2\subset\mathbb D,\, z_0\in \mathbb D, \,\theta\in\R$, and $\varphi(z) = e^{i\theta}\cdot\dfrac{z-z_0}{1-z\bar z_0}$,  $z\in \mathbb D,$ is a conformal automorphism of $\mathbb D$.  Then $\Cap2(E_1, E_2) = \Cap2(\varphi(E_1), \varphi(E_2))$. Moreover,  $\Cap2(E, \mathbb D^{(c)}) = \Cap2(\varphi(E), \mathbb D^{(c)})$ for any $E\subset\mathbb D$.
\end{predl}

The proof is given in the Appendix.

For any $j=1,2,\dots$ function $\mathfrak v_j$ found in theorem~\ref{th:reprodicing_existence} satisfies equations $\mathfrak v_j|_{B_j}=1$, $\mathfrak v_j|_{B_{j'}}=0, \, j'\neq j,$  and minimizes Dirichlet integral $\int_{\mathbb C}|\nabla v|^2\dl$ under these conditions (that was the construction of this function in the proof of theorem~\ref{th:reprodicing_existence}). From this we obtain relation
\begin{equation}
\label{eq:per_cap_norm}
\|\Per_j\|_{\left(\lo\right)^*} = \|\kappa_j\|_\lo = \|\mathfrak v_j\|_\sob= \Cap2(B_j, \Omega^{(c)} \setminus B_j)^{1/2}.
\end{equation}
(We denote by $\Omega^{(c)}$ the set $\mathbb C \setminus \Omega=\mathbb D^{(c)}\cup \bigcup\limits_{j=1}^\infty B_j$ where $\mathbb D^{(c)}=\{z\in\mathbb C\colon |z| \ge 1\}$.) This immediately leads us to necessary conditions for Bessel and interpolation properties:

\begin{predl} Let $\Omega$ be regular domain.
	\label{necessary_cap}
\begin{enumerate} 
	\item The following equality is true: 
	 $$
	 \tilde C_B^2(\Omega)=\sup\limits_{j\in\mathbb N}\Cap2(B_j, \Omega^{(c)} \setminus B_j).
	 $$
	\item If  $\Omega$ possesses interpolation property then 
	$$
	\Cap2(B_j, \Omega^{(c)} \setminus B_j) \ge C_I^{-2}(\Omega)
	$$
	 for any $j=1,2,\dots$.
	\end{enumerate}
\end{predl}


One can easily give examples showing that $\Omega$ may not have \emph{interpolation} property even if quantities $\Cap2(B_j, \Omega^{(c)} \setminus B_j)$, $j=1,2,\dots$, are bounded from the above and separated from zero uniformly. 
Namely, let $\delta>0$ be small enough, and put $\Omega_\delta=\mathbb D\setminus \left(\bar{\mathcal B}(-2\delta,\delta)\cup\bar{\mathcal B}(2\delta,\delta)\right)$. According to theorem~\ref{th:interp_sufficient} proved in the below, $C_I(\Omega_\delta)\xrightarrow{\delta\to0}+\infty$, while capacities $\Cap2\left(\bar{\mathcal B}(\pm 2\delta,\delta), \mathbb D^{(c)} \cup \bar{\mathcal B}(\mp 2\delta,\delta)\right)$ are bounded from the above and from the below (see proposition~\ref{predl:capacity_diam_low_estim} and corollary~\ref{sled:weak_separ_cap} below). Acting with conformal shifts as in example~\ref{example:inverse} in subsection~\ref{section:examples} (see below), one can construct from a sequence of such domains with~$\delta$ tending to zero a domain with infinite interpolation constant for which capacities of the form $\Cap2(B_j, \Omega^{(c)} \setminus B_j)$ still satisfy upper and lower estimates.

Nevertheless, Bessel property and weak Bessel property turn to be equivalent.

\subsection{Weak Bessel property implies Bessel property}

\begin{theorem}
\label{th:weak_strong_bessel}
	If regular domain $\Omega$ has weak Bessel property, then it has Bessel property, moreover, $C_B(\Omega) \le \sqrt 2 \tilde C_B(\Omega)$.
\end{theorem}

\noindent {\bf Proof 1.}
Let us check the conditions of the first part of proposition~\ref{criterii_adm}. Let $\{a_j\}_{j=1}^\infty$ be a real-valued \emph{finitely supported} sequence. We are going to show that there exists a function $u\in \mao$ such that $u|_{B_j}=a_j, \, j=1,2,\dots,$ and $\|u\|_{\sob} \le \sqrt 2 \tilde C_B(\Omega) \cdot \|a\|_{\ell^2}$. First, let us assume that $a_j \ge 0$ for all $j=1,2,\dots$. Recall that functions $\mathfrak v_j$ were constructed in theorem~\ref{th:reprodicing_existence}. We have $\mathfrak v_j|_{B_{k}} \equiv \mathds 1_{\{j=k\}} \, (j,k\in\mathbb N)$. Put
$$
u = \max\limits_{j\in\mathbb N} a_j \mathfrak v_j.
$$
Then $u|_{B_j} \equiv a_j$ for all $j=1,2,\dots$ (since all $a_j$ are non-negative).
Further, sequence $a$ is finitely supported; it is well-known that maximum of two Sobolev functions is also a Sobolev function admitting gradient estimate. Hence $u\in\sob$ and 
$$
|\nabla u| \le \max\limits_{j\in\mathbb N} a_j |\nabla \mathfrak v_j|
$$
almost everywhere. From this we obtain
$$
\int\limits_{\mathbb D} |\nabla u|^2 \dl \le  \sum\limits_{j\in\mathbb N}a_j^2\int\limits_{\mathbb D}|\nabla \mathfrak v_j|^2\dl.
$$
Due to weak Bessel property we have $\|\mathfrak v_j\|^2_{\sob} \le \tilde C_B^2(\Omega)$ for all  $j=1,2,\dots$ (see equation~(\ref{eq:per_cap_norm})). 
Therefore, $\|u\|_{\sob} \le \tilde C_B(\Omega)\cdot \|a\|_{\ell^2}$. So, the function $u$ satisfies all the requirements we posed. 

Now suppose that sequence $a$ is finitely supported but does not, in general, have non-negative terms. Represent this sequence as $a=a^{(1)}-a^{(2)}$ where sequences $a^{(1)}=\{a^{(1)}_j\}_{j=1}^\infty$, \, $a^{(2)}=\{a^{(2)}_j\}_{j=1}^\infty$ are finitely supported, have non-negative components while $\|a\|_{\ell^2}^2=\|a^{(1)}\|_{\ell^2}^2+\|a^{(2)}\|_{\ell^2}^2 $. According to the already proved, there exist functions $u^{(1)}, u^{(2)} \in \mao$ for which $u^{(\alpha)}|_{B_j} \equiv a^{(\alpha)}_j$  $(j=1,2,\dots)$, $\|u^{(\alpha)}\|_{\sob} \le \tilde C_B(\Omega) \cdot \|a^{(\alpha)}\|_{\ell^2}$ \, ($\alpha=1,2$). For $u=u^{(1)}-u^{(2)}$ we have $u\in \mao$,  $u|_{B_j}=a_j$ for all $j=1,2,\dots$. Moreover, 
$$
\|u\|_{\sob} \le \|u^{(1)}\|_{\sob}+\|u^{(2)}\|_{\sob} \le \tilde C_B(\Omega)\left(\|a^{(1)}\|_{\ell^2}+\|a^{(2)}\|_{\ell^2}\right) \le \sqrt 2\tilde C_B(\Omega) \cdot \|a\|_{\ell^2}.
$$
Thus we have checked the conditions of the first part of proposition~\ref{criterii_adm}. Theorem is proved.
$\blacksquare$

\medskip

\noindent {\bf Proof 2.} Let us show that $\|\sum_{j=1}^{\infty} a_j \kappa_j\|_{\lo} \le \sqrt 2 \,\tilde C_B(\Omega)\cdot \|\{a_j\}_{j=1}^\infty\|_{\ell^2}$ for any finitely supported sequence $\{a_j\}_{j=1}^\infty$ (reproducing kernels $\kappa_j$ were found in theorem~\ref{th:reprodicing_existence}). By proposition~\ref{pred:Riesz_Criterii}, this is enough to prove our statement. First we get an estimate for non-negative sequences $\{a_j\}_{j=1}^\infty$. To this end, we notice that, by proposition~\ref{predl:negative_prod}, we have  $\langle\kappa_j, \kappa_{k}\rangle_{\lo} < 0$ for $j\neq k, \, j,k=1,2,\dots$. Hence
$$
\left\|\sum_{j=1}^{\infty} a_j \kappa_j\right\|^2_{\lo} = 
\sum\limits_{j=1}^\infty a_j^2 \|\kappa_j\|_{\lo}^2 + 2\cdot \sum\limits_{j>k} a_j a_k \langle\kappa_j, \kappa_{k}\rangle_{\lo} \le \tilde C_B^2(\Omega) \cdot \|\{a_j\}_{j=1}^\infty\|_{\ell^2}^2,
$$
because all $a_j$ have the same sign and $\|\kappa_j\|_{\lo} \le \tilde C_B(\Omega)$ for all $j=1,2,\dots$. Now the desired estimate for arbitrary, not necessarily non-negative sequences $\{a_j\}_{j=1}^\infty$ follows from the obtained estimate by subdivision of such a sequence into positive and negative parts. $\blacksquare$

\medskip

\noindent {\bf Remark.} Suppose that $H$ is an abstract Hilbert space and a countable system $\{\xi_j\}_{j=1}^\infty$ in $H$ is understood as a system of reproducing kernels in interpolation problem for operator $T\colon H\to \ell^2, \, Tx=\{\langle{x, \xi_j}\rangle_H\}_{j=1}^\infty, \, x\in H$. Is it true that weak Bessel property in such a problem implies Bessel property? No, of course. Is that true that interpolation and weak Bessel properties imply Bessel property? The answer is also negative. More precisely, we have the following 

\begin{predl}
	\label{predl:reverse_counterex}
	There exists a sequence of vectors $\xi_1, \xi_2, \dots$ in Hilbert space $\ell^2$ such that:
	\begin{enumerate}
		\item $\sup\limits_{j\in\mathbb N} \|\xi_j\|_{\ell^2} < +\infty$;
		\item $\|\sum_{j=1}^\infty a_j \xi_j\|_{\ell^2}\ge C_1\cdot \|\{a_j\}_{j=1}^\infty\|_{\ell^2}$ for any finitely supported sequence $\{a_j\}_{j=1}^\infty$ with some absolute constant $C_1>0$;
		\item a reverse estimate $\|\sum_{j=1}^\infty a_j \xi_j\|_{\ell^2}\le C_2\cdot \|\{a_j\}_{j=1}^\infty\|_{\ell^2}$ for finitely supported sequences $\{a_j\}_{j=1}^\infty$ is not true for none $C_2<+\infty$.
	\end{enumerate}
\end{predl}

The proof is given in the Appendix.

\subsection{Hyperbolic metric and metric capacity estimates}

\label{subsection:hyperbolic_metric}

Our problem is conformally invariant. Thus we may expect that  the necessary and sufficient conditions will be conformally invariant. Let us recall some notions of hyperbolic geometry (see also~\cite{Katok}).

The hyperbolic metric in the disk is the metric $\dfrac{2|dz|\hspace{4pt}}{1-|z|^2}$. Geodesics in this metric are arcs of circles (or segments of lines) orthogonal to $\partial\mathbb D$. M{\"o}bius automorphisms  
$$
z\mapsto e^{i\theta}\cdot \frac{z-z_0}{1-z \bar z_0}, ~~ \theta \in \R,\, z_0\in \mathbb D,
$$
of the disk are isometries in the hyperbolic metric. The following formula for hyperbolic distance $\dist_H$ is well-known:
\begin{equation}
\label{hyperb_dist}
\dist_H(z_1, z_2) = \arctanh\left|\frac{z_1 - z_2}{1-z_1 \bar z_2}\right|, ~~ z_1, z_2 \in\mathbb D.
\end{equation}
Denote by $\mathcal B_H(z, r)$ the open ball in hyperbolic metric of radius $r\ge 0$ and centered in a point $z \in\mathbb D$. Form formula~(\ref{hyperb_dist}) it follows, in particular, that $\mathcal B_H(0, 1) = \mathcal B(0, \tanh 1)$.

We will make use of \emph{hyperbolic diameters} of sets $B_j$ which are their diameters in the hyperbolic metric. Denote these quantities by $\diam_H(B_j)$ (symbol $\diam(B_j)$ stands for the usual Euclidean diameter of $B_j$).

Let us point out the relation between hyperbolic diameters of sets to their Euclidean metric properties.
Let $E$ be a compact set in $\mathbb D,\, b=\dfrac{\diam E}{\dist(E, \partial \mathbb D)}$. It is easy to prove that  
$$
2b \ge \diam_H (E) \ge \arctanh\left(\dfrac b{4+b}\right) = \dfrac12\log\left(\dfrac b2+1\right).
$$
This immediately provides the following

\begin{predl}
\label{sled:diamh_metric}
	Let $\Omega$ be a regular domain. 
	\begin{enumerate}
		\item $\inf\limits_{j\in\mathbb N}\diam_H(B_j)>0$ if and only if $\inf\limits_{j\in\mathbb N} \dfrac{\diam E}{\dist(E, \partial \mathbb D)}>0$; any of these two infima can be estimated from the below in terms of the other $\inf$. 
		
		\item $\sup\limits_{j\in\mathbb N}\diam_H(B_j)<+\infty$ if and only if $\sup\limits_{j\in\mathbb N}\dfrac{\diam E}{\dist(E, \partial \mathbb D)}<+\infty$; any of these two suprema can be estimated from the above in terms of the other  $\sup$.
	\end{enumerate}
\end{predl}

We make essential use of the following separatedness conditions of holes $B_j$.

\begin{define}
\label{def:separatedness}
	Pick $\eps>0$. Let us say that \emph{holes $B_j$ are $\,\eps$-strongly separated} if for any two distinct indices $j$ and $j'$ the following estimate is held:
	$$
	\dist(B_j, B_{j'}) \ge \eps \cdot \max\{\diam B_j, \diam B_{j'}\}.
	$$
	If for any two distinct indices $j$ and $j'$ the estimate 
	$$
	\dist(B_j, B_{j'}) \ge \eps \cdot \min\{\diam B_j, \diam B_{j'}\},
	$$
	is held then we say that \emph{holes $B_j$  are $\,\eps$-weakly separated}.
	
	We will say that holes are \emph{strongly (or weakly) separated} if they are $\eps$-strongly \emph{(}respectively, $\eps$-weakly\emph{)} separated for some $\eps>0$.
\end{define}

In the case of round holes strong separatedness condition is equivalent to that disks $(1+\eps) B_j$ (that are disks with the same centers and $(1+\eps)$ times enlarged radii) are pairwise disjoint for some $\eps>0$. If numbers $\diam_H(B_j)$ for all $j=1,2,\dots$ do not exceed some given value  then separatedness conditions can be reformulated in a M{\"o}bius invariant form, that is via hyperbolic metric. However, we are more convenient to work with Euclidean distances.

We will permanently need the following estimate:

\begin{predl}
	\label{predl:capacity_diam_low_estim}
	Suppose that any of two sets $E_1, E_2\subset\mathbb C$ is either a closure of  a domain with $C^\infty$-smooth boundary or  $\mathbb D^{(c)},\, E_1 \cap E_2 =\varnothing$. Then 
	$$
	\frac{\min\{\diam E_1, \diam E_2\}}{\dist(E_1, E_2)} \le \frac{\exp\left(24\pi \cdot \Cap2(E_1, E_2)\right)-1}{\pi}.
	$$
\end{predl}

 The proof is given in the Appendix.

\begin{predl}
	\label{predl:weak_separated_bessel}
	If a regular domain $\Omega$ has weak Bessel property then:
		\begin{enumerate}
			\item quantity  $\sup\limits_{j\in\mathbb N} \diam_H(B_j)$~is finite and can be estimated from the above only through~$\tilde C_B(\Omega)$; 			
			\item holes $B_j$ are $\eps$-weakly separated for some $\eps>0$, depending only on $\tilde C_B(\Omega)$.
		\end{enumerate}
\end{predl}

The proof is given in the Appendix.

\begin{predl} 
	\label{predl:strong_sep_Bessel}
	Let $\Omega$ be a  regular domain possessing weak Bessel property.
	If $\diam_H(B_j)\ge r$ for all $j=1,2,\dots$ with some constant $r>0$ then sets $B_j$ are $\eps\mbox{-strongly}$  separated with some $\eps>0$ depending only on $r$ and $\tilde C_B(\Omega)$.
\end{predl}

The proof is given in the Appendix.

Further,  let us estimate the capacity from the above.

\begin{predl}
	\label{predl:capacity_upper}
	Let $E\subset\mathbb C$ be a bounded set, $\Lambda=\diam E>0$, $\eps>0$ and  $U$ be $(\Lambda\eps)$-neighbourhood of  set $E$, \, $U^{(c)}=\mathbb C\setminus U$. Then $\Cap2(E, U^{(c)}) \le \dfrac{\pi (1+2\eps)^2}{\eps^2}$.
\end{predl}

The proof is given in the Appendix.

\begin{sled}
\label{sled:weak_separ_cap}
If sets $E_1, E_2 \subset \mathbb C$ are bounded and   $$\dist(E_1, E_2) \ge \eps \cdot \min\{\diam E_1, \diam E_2\}$$ for some $\eps>0$ then $\Cap2(E_1, E_2)\le \dfrac{\pi (1+2\eps)^2}{\eps^2}$.
\end{sled}

\medskip

\noindent {\bf Remark.} According to proposition \ref{predl:capacity_diam_low_estim} and corollary \ref{sled:weak_separ_cap}, condenser formed by \emph{connected} plates $E_1, E_2$ looks, in what concern upper or lower boundedness of its capacity, like a condenser formed by disks of diameters $\diam E_1$ and $\diam E_2$ placed on the distance $\dist(E_1, E_2)$ from each other.

\medskip

Finally, we have the following 

\begin{predl}
\label{predl:strong_sep_Bessel_suff}
	Let $\eps$ be positive and $\Omega$ be a  regular domain. If holes $B_j$ are $\eps$-strongly separated  and $\sup\limits_{j\in\mathbb N} \diam_H(B_j) < +\infty$ then $\Omega$ has Bessel property and the Bessel constant $C_B(\Omega)$ is bounded from the above by a quantity depending only on $\eps$ and $\sup\limits_{j\in\mathbb N} \diam_H(B_j)$.
\end{predl}

The proof is given in the Appendix.

\section{Partial complete interpolation criterion}

\label{section:partial_criterion}

In this section we will derive criterion of complete interpolation in a regular domain $\Omega$ in the case when one of the following conditions is fulfilled: either $\inf\limits_{j\in\mathbb N}\diam_H(B_j)>0$ or sets~$B_j$ are strongly separated. If one of these conditions is satisfied then complete interpolation is equivalent to both of them. However, we will construct an example showing that no one of these conditions is necessary for complete interpolation.

The main result of this section, theorem \ref{th:predv_criterii}, is overlapped by theorem~\ref{th:full_interp_criterion}. Nevertheless, patch method and analogue of Carleson measures appearing in this section can be of an independent interest.

\subsection{A sufficient interpolation condition}

\label{subsection:interp_suff}

The goal of this subsection is to prove that if hyperbolic diameters of holes are bounded from the below then domain $\Omega$ has interpolation property. This is enough if holes are disks. If holes have arbitrary form then we will also impose separatedness condition.

Let $\Omega$ be a regular domain. We bring into consideration a (positive) measure in $\mathbb D$: put
\begin{equation}
\label{mu_omega_def}
\mu_\Omega:=\sum\limits_{j=1}^\infty \frac{\mathds{1}_{B_j}\cdot \lambda_2}{\lambda_2(B_j)}.
\end{equation}
Now, we will say that a positive measure $\mu$ in unit disk $\mathbb D$ \emph{has} (MC) \emph{property} if for any function $u\in \overset{\circ}{W}{}^{1,2}(\mathbb D)$ we have
\begin{equation}
\label{carleson_def}
\int\limits_{\mathbb D} |\nabla u|^2\dl \ge {c} \cdot \int\limits_{\mathbb D} u^2\,d\mu
\end{equation}
with some $c>0$ not depending on $u$. Abbreviation (MC) stands for "Maz'ya-Carleson". The choice of such a term is related to the fact that sharp results on embeddings of Sobolev spaces into $L^p(\mu)$ with $\mu$ be some measure in $\R^n$, can be found in the well-known monograph by V.G. Maz'ya~\cite{Mazya}. To the other hand, in classic problems on analytic interpolation, Carleson measures arise in the similar context (see~\cite{Nikolskiy}). 

Recall that, by theorem~\ref{th:Hardy_neq}, measure $\mu_0=\dfrac{\lambda_2}{(1-|z|)^2}$, possesses (MC) property. 

\begin{lemma}
	\label{carleson_sufficient}
	Let $\Omega$ be a  regular domain. If the measure  $\mu_\Omega$ given by 
	\emph{(\ref{mu_omega_def})} satisfies \emph{(MC)}~condition 
	\emph{(}estimate \emph{(\ref{carleson_def})}\emph{)}, then  domain $\Omega$ 
	has interpolation property whereas $C_I(\Omega)$ can be estimated from the 
	above through the constant from \emph{(MC)} condition for the measure 
	$\mu_\Omega$.
\end{lemma}

\noindent {\bf Proof.} Notice that $\mu_\Omega(B_j)=1$ because sets $B_j$ are disjoint. Thus, by the second condition of proposition~\ref{criterii_adm}, interpolation property of  $\Omega$ can be reformulated in the following way:
$$
\int\limits_{\mathbb D} |\nabla u|^2\dl \ge C_I^{-2}(\Omega) \cdot \int\limits_{\mathbb D} u^2\,d\mu_\Omega
$$
for any admissible function $u\in\mao$. But if the measure $\mu_\Omega$ has (MC) property then such an estimate is true even for all $u\in \sob$. $\blacksquare$

\medskip

Now let us prove a sufficient interpolation condition for domains with round holes. 

\begin{theorem}
	\label{th:disks_interp_sufficient}
	Suppose that holes $B_j$ are disks and numbers $\diam_H(B_j)$ are bounded from the below. Then domain $\Omega$ has interpolation property and $C_I(\Omega)$ does not exceed some constant depending only on $\inf\limits_{j\in\mathbb N} \diam_H(B_j)$.
\end{theorem}

\noindent {\bf Proof.} By the proposition~\ref{sled:diamh_metric}, $c \cdot \diam B_j \ge \dist(B_j, \partial\mathbb D)$ with some $c$ depending only on $\inf\limits_{j\in\mathbb N}\diam_H(B_j)$. Hence  
$1-|z|\le (c+1)\cdot \diam B_j$ for each $z\in B_j$. We have 
\begin{equation}
\label{neq:disk_mu0_compare}
\frac{\pi}{4 \lambda_2(B_j)} = \frac{1}{\diam^2 B_j} \le \frac{(c+1)^2}{\left(1-|z|\right)^2}, ~~~ z\in B_j,
\end{equation}
and, from this and by the definition of measures $\mu_0$ and $\mu_\Omega$, it follows that 
$$
\mu_\Omega \le \frac{4 (c+1)^2}{\pi}\cdot \mu_0.
$$
Now (MC) property of measure $\mu_\Omega$ is provided by this property of measure $\mu_0$, and  interpolation property of $\Omega$ is immediate from lemma~\ref{carleson_sufficient}. In addition, $C_I(\Omega)$ can be estimated from the above by a constant depending only on $\inf\limits_{j\in\mathbb N}\diam_H(B_j)$. $\blacksquare$

\medskip

In the subsection~\ref{section:examples} we will give an example~\ref{example:inverse} that shows that interpolation property (even together with Bessel property) does not imply (MC) property of measure~$\mu_\Omega$ (even for domains with round holes). In other words, estimate~(\ref{carleson_def}) for functions~$u$ which are constant on any set $B_j$ does not, in general, imply such an estimate for arbitrary $u\in \sob$.

If holes $B_j$ are not disks then the estimate~(\ref{neq:disk_mu0_compare}) in the proof of the last theorem may fall. Let us give an example showing that for domains with  holes of arbitrary form the  assertion of theorem~\ref{th:disks_interp_sufficient} is not true indeed. For this aim, pick a number $n\in \mathbb N$ large enough and cut from $\mathbb D$ line segments of the form $[1/3 \cdot e^{2\pi ik/n}, 1/2\cdot e^{2\pi ik/n}]$ with $k=0, 1, \dots , n-1$. Further, let holes $B_k$ with these $k$ be closed $\eps$-neighbourhoods of these segments where $\eps=\eps_n>0$ is small enough in order to have pairwise disjointedness of such neighbourhoods. Put $\Omega_n =\mathbb D \setminus \bigcup_{k=0}^{n-1}B_k$. Quantities $\diam_H B_k$ are separated from zero by a constant not depending on $n$ (and on $\eps$). To the other hand, it is easy to choose a function $u\in \sob$ which equals $1$ on $\mathcal B(0, 3/4)$. For such a function we have $\sum_{k=0}^{n-1} (u|_{B_k})^2 \ge n$. By the second assertion of proposition~\ref{criterii_adm}, this means that $C_I(\Omega) \ge \sqrt{n} \cdot c$ with an absolute constant $c>0$. Acting by conformal shifts like in example~\ref{example:inverse} in subsection~\ref{section:examples} (see below), it is possible to construct a domain $\Omega$ such that holes in $\Omega$ have hyperbolic diameters bounded from the below, but $C_I(\Omega)=+\infty$.

Now we will prove that domain $\Omega$ will still have interpolation property if we add holes separatedness condition to the estimate $\inf\limits_{j\in\mathbb N}\diam_H(B_j)>0$. We are mainly interested in the complete interpolation property of domain $\Omega$; the estimate $\sup\limits_{j\in\mathbb N} \diam_H(B_j)<+\infty$ and weak separatedness of holes are necessary for that. In the proof of proposition~\ref{predl:strong_sep_Bessel}, it is established that  estimates $\inf\limits_{j\in\mathbb N} \diam_H(B_j)>0$,  $\sup\limits_{j\in\mathbb N} \diam_H(B_j)< +\infty$ together with weak separatedness of holes~$B_j$ imply their strong separatedness. This condition is what we add to estimate $\inf\limits_{j\in\mathbb N} \diam_H(B_j)>0$.

\begin{theorem}
	\label{th:arbitrary_interp_sufficiency_separated}
	Let domain $\Omega$ be regular, holes $B_j$ be  $\eps$-strongly separated with some $\eps>0$ and  $\inf\limits_{j\in\mathbb N} \diam_H(B_j) >0$. Then $\Omega$ possesses interpolation property and the constant~$C_I(\Omega)$ can be estimated from the above through $\eps$ and $\inf\limits_{j\in\mathbb N} \diam_H(B_j)$.
\end{theorem}

\begin{wrapfigure}[20]{r}{0.55\textwidth}
	\centering
	{\includegraphics{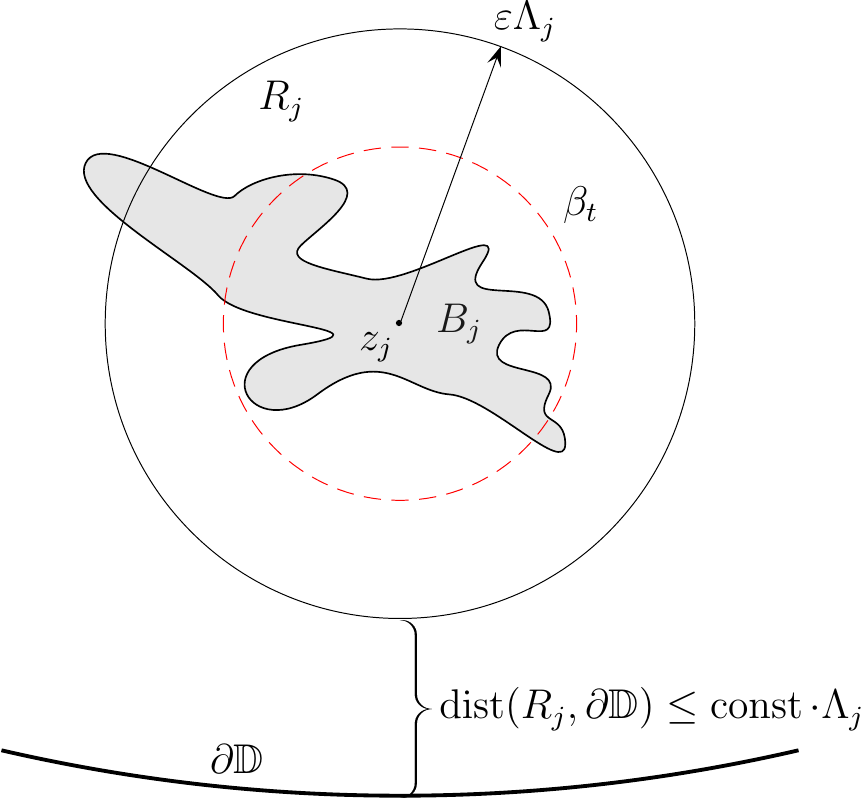}}
	\caption{Disk $R_j$ looks like the hole $B_j$}
	\label{fig:Rj_arbitrary}	
\end{wrapfigure}

\noindent {\bf Proof.} Pick $u\in\mao$. Denote by $a_j$ the constant value $u|_{B_j}$. According to the second assertion of proposition~\ref{criterii_adm}, we have to prove that 
$$
\|u\|_{\sob}^2 \ge c\cdot \sum\limits_{j=1}^\infty a_j^2,
$$
where a constant $c>0$ depends on the constants in the statement of theorem but not on $u$. We will not establish (MC) property of measure $\mu_\Omega$ though is seems to be possible. Instead we reduce the desired estimate to  theorem \ref{th:Hardy_neq} by replacement of sets~$B_j$ by appropriate disks $R_j$ placed near $B_j$ and pairwise disjoint (see fig. \ref{fig:Rj_arbitrary}). In addition, $\diam R_j$ will be comparable to $\diam B_j$. The possibility of choice of such $R_j$ is provided by separatedness condition of holes $B_j$. On sets $R_j$ function $u$ will not still be constant. If disk~$R_j$ is close to $B_j$ then $a^2_j$ can be estimated via $\dfrac{1}{\lambda_2(R_j)}\displaystyle\int\limits_{R_j} u^2\dl$ and  Dirichlet integral of  $u$ over an appropriate set. The estimate for $\displaystyle\sum\limits_{j=1}^\infty\dfrac{1}{\lambda_2(R_j)}\displaystyle\int\limits_{R_j} u^2\dl$ will provided, like in the previous theorem,  by the (MC) property of measure $\mu_0$.

Put $\Lambda_j = \diam B_j$ and let  $U_j$ be the  $(\eps\Lambda_j/2)$-neighbourhood of set $B_j$. By the hypothesis,  sets~$U_j$ are pairwise disjoint. We may assume that $\eps < 1$. Pick arbitrary $z_j \in B_j, \, j=1,2,\dots$.  Put $R_j := \mathcal B(z_j, \eps \Lambda_j/2)$ and, for  $t \in [0, \eps \Lambda_j/2]$, define a circumference $\beta_t=\partial \mathcal B(z_j, t)$ concentric to $R_j$. Notice that $\beta_t \cap B_j\neq \varnothing$ because $z_j \in B_j, \, \diam B_j > \eps \Lambda_j$ and $B_j$ is connected. Define the function $u$ to be zero outside of~$\mathbb D$, now $u$ is defined on $R_j$ even if $R_j \not\subset \mathbb D$.

Now, correct the function $u\in \sob$ on a set of zero measure such that $u$ will be absolutely continuous on circumferences $\beta_t$ with almost all radii $t\in [0, \eps \Lambda_j/2]$. Suppose that  $t$ is one of such radii. Function $u$ takes value $a_j$ on~$\beta_t$, hence, for any point $z\in \beta_t$ we have estimate
$$
|a_j| \le |u(z)| + \int\limits_{\beta_t}|\nabla u|\,d\mathcal H^1 \le 
|u(z)| + \sqrt{2\pi t}\left(\int\limits_{\beta_t}|\nabla u|^2\,d\mathcal H^1\right)^{1/2},
$$
from which
$$
a_j^2 \le 
2u^2(z) + 4\pi t\int\limits_{\beta_t}|\nabla u|^2\,d\mathcal H^1.
$$
Integration by $z\in \beta_t$ over  $\mathcal H^1$ gives
$$
2\pi t a_j^2 \le 2\int\limits_{\beta_t} u^2(z)\,d\mathcal H^1+8\pi^2 t^2 \int\limits_{\beta_t}|\nabla u|^2\,d\mathcal H^1.
$$
Now integrate the obtained inequality by $t\in [0,\eps \Lambda_j/2]$; taking in account that $t \le \eps \Lambda_j/2$, we get
$$
\frac{\pi \eps^2 \Lambda_j^2 a_j^2}{4} \le 2 \int\limits_{R_j}u^2 \dl+2\pi^2 \eps^2 \Lambda_j^2 \int\limits_{R_j}|\nabla u|^2\dl.
$$
Finally, we obtain the estimate
$$
a_j^2 \le \frac{2}{\lambda_2(R_j)} \int\limits_{R_j}u^2 \dl+8\pi \int\limits_{R_j}|\nabla u|^2\dl.
$$
Disks  $R_j$ are pairwise disjoint since they are situated in disjoint sets $U_j$. Summation of the last inequalities over $j$ gives us
\begin{equation}
\label{sum_Rj_estim}
\sum\limits_{j=1}^\infty a_j^2 \le 2\int\limits_{\mathbb C} u^2 \, d\mu +8\pi\cdot \|u\|^2_{\sob}=
2\int\limits_{\mathbb D} u^2 \, d\mu +8\pi\cdot \|u\|^2_{\sob},
\end{equation}
where $\mu=\displaystyle\sum\limits_{j=1}^\infty  \dfrac{\mathds{1}_{R_j}\cdot 
\lambda_2}{\lambda_2(R_j)}$ is a measure on $\mathbb C$; in the last equation 
we shrank the integration domain from  $\mathbb C$ to~$\mathbb D$ since $u=0$ 
in $\mathbb D^{(c)}$ almost everywhere. (Recall that in our construction disks 
$R_j$ may not lie in $\mathbb D$.)

The hypothesis of theorem implies the existence of a constant $c_1<+\infty$ (depending only on $\inf\limits_{j=1,2,\dots} \diam_H(B_j)$), for which
$$
\frac{\dist(B_j, \partial\mathbb D)}{\Lambda_j} \le c_1
$$
for all  $j=1,2,\dots$.
Then $\dist(z_j, \partial\mathbb D)\le \diam B_j+ \dist(B_j, \partial\mathbb D) \le (1+c_1) \Lambda_j$. For any point $z \in R_j\cap \mathbb D$ we have $1-|z| \le (1-|z_j|)+|z-z_j|\le (1+c_1+\eps/2)\Lambda_j$ and hence
$$
\frac{1}{(1-|z|)^2} \ge  \frac{1}{((1+c_1+\eps/2) \Lambda_j)^2} = \frac{\pi\eps^2}{4(1+c_1+\eps/2)^2}\cdot\frac1{\lambda_2(R_j)}.
$$
Thus $\mathds 1_{\mathbb D}\cdot \mu \le \dfrac{4(1+c_1+2\eps)^2}{\pi\eps^2} \cdot \mu_0$. Now, (MC) property of measure $\mu_0$ (theorem~\ref{th:Hardy_neq}) implies that  $\int\limits_{\mathbb D} u^2 \,d\mu \le C \cdot \|u\|_{\sob}^2$, where constant $C>0$ depends only on constants from the statement of the theorem. From this and~(\ref{sum_Rj_estim}) we finally conclude that 
$$
\sum\limits_{j=1}^\infty a_j^2 \le (8\pi+2C)\cdot \|u\|^2_\sob.
$$
Here  function $u\in \mao$ is arbitrary and $a_j=u|_{B_j}$. So, we checked the hypothesis of the second assertion of proposition~\ref{criterii_adm}. Theorem is proved.
$\blacksquare$

\medskip

\subsection{Lower estimate for $\diam_H(B_j)$}

\label{section:cylinder}

In the previous subsection we proved that condition $\inf\{\diam_H(B_j)\mid j=1,2,\dots\} > 0$ is sufficient for interpolation (in the case of holes of an arbitrary form we should strong separatedness condition). Now we are going to prove that is holes are strongly separated then interpolation property implies that $\inf\{\diam_H(B_j)\mid j=1,2,\dots\} > 0$. The proof is given only for domains with round holes. Our argument uses \emph{patch method}. 

\begin{theorem}
	\label{patch_theorem}
	Let  $\Omega$ be a  domain with round holes. Suppose that $B_j$ are  $\eps$-strongly separated for some $\eps>0$ and, moreover, that  $\sup\limits_{j\in\mathbb N} \diam_H(B_j) < +\infty$. 
	
	If domain $\Omega$ possesses interpolation property then $\inf\limits_{j\in\mathbb N} \diam_H(B_j)$ is positive and can be estimated from the below only through $\eps, \, \sup\limits_{j\in\mathbb N} \diam_H(B_j)$ and $C_I(\Omega)$.
\end{theorem}

\noindent {\bf Proof.} If $B\subset\mathbb C$ is an open or closed disk centered in some point $z\in\mathbb C$ and of radius $r\ge 0$ and also  $a\ge 0$ then we denote by $aB$ respectively open or closed disk centered in  $z$ and of radius $ar$.

By proposition~\ref{sled:diamh_metric} there exists $\eps_0>0$ depending only on $\sup\limits_{j\in\mathbb N} \diam_H(B_j)$ such that $\dist(B_j, \partial\mathbb D) \ge \eps_0 \cdot \diam B_j$ for all $j=1,2,\dots$. Decreasing, if necessary,  $\eps_0$ we may assume that $\eps_0< \eps/2$. 
By the choice of  $\eps_0$, disks $(1+\eps_0)B_j$ are pairwise disjoint and lie 
in $\mathbb D$ for all $j=1,2,\dots$.

We have to prove that hyperbolic diameter of any hole is bounded from the below. Let us prove this for one fixed hole assuming that it this  disk is $B_1$. By use of a conformal automorphism of disk $\mathbb D$ we reduce the question to the case when $B_1=\bar{\mathcal B}(0,r)$. (It is not hard to show that the strond separatedness property is preserved under M{\"o}buis transforms, if $\sup\limits_{j\in\mathbb N} \diam_H(B_j) < +\infty$.) Then we have to prove that radius $r$ is bounded from below. 

As a heuristic consideration, let us assume that there is just one hole in the domain. Then the solution of Dirichlet problem~(\ref{reproducing_eq}) is the function 
$$
u_1(z) = 
\begin{cases}
\dfrac{\log|z|}{\log r}, & |z|>r;\\
\vspace{-0.4cm}\\
~~~ 1, &|z|\le r
\end{cases} ~~~ (z\in\mathbb D).
$$
Calculation of Dirichlet integral of this function shows that norm $\|u_1\|_\sob$ is small if~$r$ is small, that is, function $u_1$ does not 
satisfy  the second condition of proposition~\ref{criterii_adm}. If there are another holes in  $\Omega$ then function $u_1$ will not still be admissible for $\Omega$. The idea of patch method of our proof is to correct the function $u_1$ on sets $(1+\eps_0)B_{j}$ for all $j=1,2,\dots$ such that it will become constant on $B_j$; such a function will be admissible for $\Omega$. Thus, function $u_1$ will be perturbed by certain corrections, and we will have to get upper estimates for Dirichlet integral of sum of these perturbations. 

It is convenient to use cylindrical coordinates in order to compute corrections. Denote by $\Cylr$ the cylinder of height $|\log r|$ built over the circumference of radius $1$. Let us introduce coordinates on $\Cylr$: angle $\theta$ along the circle and height $y\in (0,|\log r|)$. Coordinates $(\theta, y)$ define a family of local charts on $\Cylr$ which endows $\Cylr$ with a conformal structure.

Let us map the annulus $\mathbb D \setminus B_1$ onto $\Cylr$: put
\begin{equation}
\label{eq:conformal_cyl_def}
\varphi(z) := (\arg z, -\log|z|), ~~~ z\in \mathbb D \setminus B_1.
\end{equation}
Then $\varphi$ is a conformal mapping. In particular, it preserves harmonicity of functions. Define a function $\tilde u_1 \colon \Cylr \to \R$ as $\tilde u_1 = u_1 \circ \varphi^{-1}$,  this is push-forward of $u_1$ to the cylinder. Notice that  $u_1(\theta, y)=y/|\log r|$.

\begin{figure}
\hspace{4.875cm}
	{\includegraphics[scale=0.7]
		{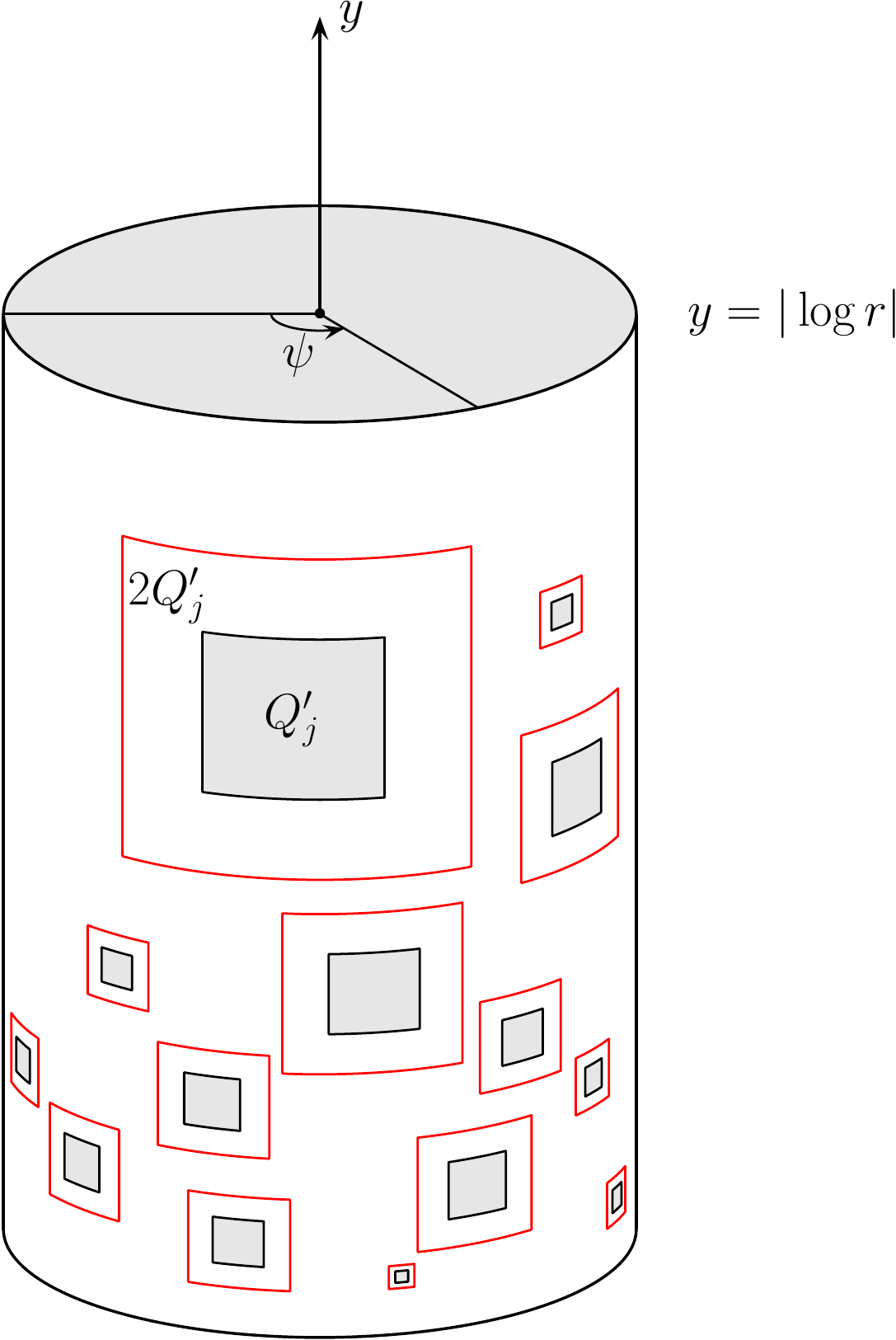}}
\caption{Cylinder with square holes}
	\label{fig:cylinder}
\end{figure}

Images $\varphi(B_j)$ are sets in $\Cylr$ with smooth boundaries. It is convenient to do estimates on cylinder assuming that these images are squares $Q'_j$ on cylinder such that their doublings $2Q'_j$ are pairwise disjoint (see fig.~\ref{fig:cylinder}). Images of disks will not, of course, have right angles, but we will achieve that by shrinking holes on cylinder. A good separatedness of squares $Q'_j$ will be attained by application to  $\Omega$ of a quasiconformal mapping with controlled distortion coefficient and then using proposition~\ref{predl:conformal_quasiconformal_invariance} on quasiconformal quasiinvariance. The class of quasiconformal mappings is less rigid than of one of conformal mappings and then allows to do such technical transitions. A quasiconformal diffeomorphism will be applied in every set $(1+\eps_0) B_j$.

Pick an \emph{absolute} constant $c_0\in [1,+\infty)$ large enough (with no dependence neither on $\eps_0$ nor on $C_I(\Omega)$); the choice of this constant will be clear after calculation of asymptotics of elementary functions in proof of lemma~\ref{lemma:planimetry_asymp_concentr}.

\begin{lemma}
Suppose that disks $(1+\eps_0)B_j$ are pairwise disjoint and lie in $\mathbb D$ for all $j=1,2,\dots$ and that $c_0\ge 1$ is some number. Then there exists a quasiconformal diffeomorphism $\psi\colon \mathbb D\to \mathbb D$ mapping disks $B_j$ to disks $\psi(B_j)$ such that $c_0^2\psi(B_j)\subset B_j$ for all $j=1,2,\dots$. Diffeomorphism $\psi$ can be chosen in such a manner that its distortion coefficient $K(\psi)$  will be bounded from the above by a value depending only on $\eps_0$ \emph{} and~$c_0$\emph{}.
\end{lemma}

In other words, holes in the  domain $\psi(\Omega)$ will be disks strongly separated with a separatedness constant large enough. 

\medskip

\noindent {\bf Proof.} 
For $z\in \mathbb D \setminus \bigcup\limits_{j=1}^\infty (1+\eps_0)B_j$ put $\psi(z)=z$. Let $B_j=\bar{\mathcal B}(z_j, r_j)$. 
In the disk $(1+\eps_0) B_j$ define $\psi$ radially:
$$
\psi(z_j+Re^{i\theta}):= z_j+f_j(R)\, e^{i\theta}, ~~ \theta \in [0, 2\pi], \, R\in [0, (1+\eps_0)r_j],
$$
where $f_j$ is a smooth increasing bijection of $[0, (1+\eps_0)r_j]$ onto $[0, 
(1+\eps_0)r_j]$ while $f_j(r_j)={r_j}c_0^{-2}$ (recall that $c_0\ge 1$ by the 
hypothesis); we may also assume that $f$ is linear on $[0,r_j]$. In order to 
have smoothness of  $\psi$ on $\partial\left((1+\eps_0)B_j\right)$ we should 
also require the equalities  $f_j^{(l)}\left((1+\eps_0)r_j\right)=0$ for all 
$l=1,2,\dots$.

Our construction is scale-invariant: namely, we may assume that ${f_j(t)=r_j f_0(t/r_j)}$ for some fixed function $f_0\colon [0,1+\eps_0]\to[0,1+\eps_0]$ depending only on $\eps_0$ and $c_0$ but not on $r_j$. Since scalings do not change distortion coefficient then  $K(\psi|_{(1+\eps_0)B_j})$ does not depend on $r_j$ but only on the choice of $f_0$.  By the construction, $\psi$ will be a diffeomorphism of $(1+\eps_0)B_j$ onto itself up to boundary $\partial (1+\eps_0)B_j$ where $\psi$ is identical. This implies that $\psi\colon \mathbb D\to\mathbb D$  is a quasiconformal diffeomorphism and its distortion coefficient $K(\psi)$ can be estimated from the above only through $\eps$ and $c_0$. For any $j=1,2,\dots$ set $\psi(B_j)$ is a disk whereas concentric disk $c_0^2\psi(B_j)$ coincides with $B_j$ by the construction (since $f_j(r_j)=r_j c_0^{-2}$). Lemma is proved.
$\blacksquare$

\medskip

Let $\psi$ be the  diffeomorphism constructed above. By proposition 
\ref{predl:conformal_quasiconformal_invariance},  $C_I(\psi(\Omega)) \le 
\sqrt{K(\psi)} \cdot C_I(\Omega)$ where distortion coefficient $K(\psi)$ 
depends only on $\eps_0$ (that is, on constants $\eps$ and 
$\sup\limits_{j\in\mathbb N} \diam_H(B_j)$ from the hypothesis of the theorem) 
and on absolute constant $c_0$. Moreover, hyperbolic diameters of holes in 
domain $\psi(\Omega)$ do not exceed  hyperbolic diameters of corresponding 
holes in domain $\Omega$. Hence we may prove our theorem assuming that disks 
$c_0^2 B_j$ are pairwise disjoint and lie in  $\mathbb D$ -- in the other case 
just pass to domain $\psi(\Omega)$.

Now continue the proof of the theorem in the above-mentioned assumption on good separatedness of holes. Recall that we want to get lower estimate for the radius $r$ of the hole $B_1=\bar{\mathcal B}(0,r)$ and conformal mapping $\varphi\colon \mathbb D \setminus B_1\to \Cylr$ is defined by~(\ref{eq:conformal_cyl_def}).

Let $\Lambda \in [0,\min(\pi,|\log r|/2)), \, \theta_0\in [0,2\pi], \, y_0\in (\Lambda, |\log r|-\Lambda)$. \emph{Square} on cylinder $\Cylr$ (centered in point $(\theta_0, y_0)$ and of side $2\Lambda$) is the set $Q=\{(\theta, y)\colon |\theta -\theta_0|\le \Lambda,\, |y-y_j|\le \Lambda\}$. If, moreover, $\Lambda < \min(\pi/2,|\log r|/4)$ and  $y_0\in (2\Lambda, |\log r|-2\Lambda)$ then \emph{doubling} of such a square $Q$ is the set $2Q:=\{(\theta, y)\colon |\theta -\theta_0|\le 2\Lambda,\, |y-y_j|\le 2\Lambda\}$.

\begin{lemma}
	\label{lemma:planimetry_asymp_concentr}
	Suppose that absolute constant $c_0\ge 1$ is large enough and that  round holes $B_j$ in a regular domain  $\Omega$ are such that sets $c_0^2B_j$, $j=1,2,\dots$, are pairwise disjoint and lie in~$\mathbb D$.  Then, for 	$j=2,3,\dots$, there exist squares $Q_j'$ on the cylinder $\Cylr$ such that 	
	$$
	\varphi(B_j)\subset Q_j'\subset 2Q_j' \subset \varphi(c_0^2 B_j).
	$$
\end{lemma}

\begin{figure} 
\hspace{1.7cm}
{\includegraphics[scale=1.3]
	{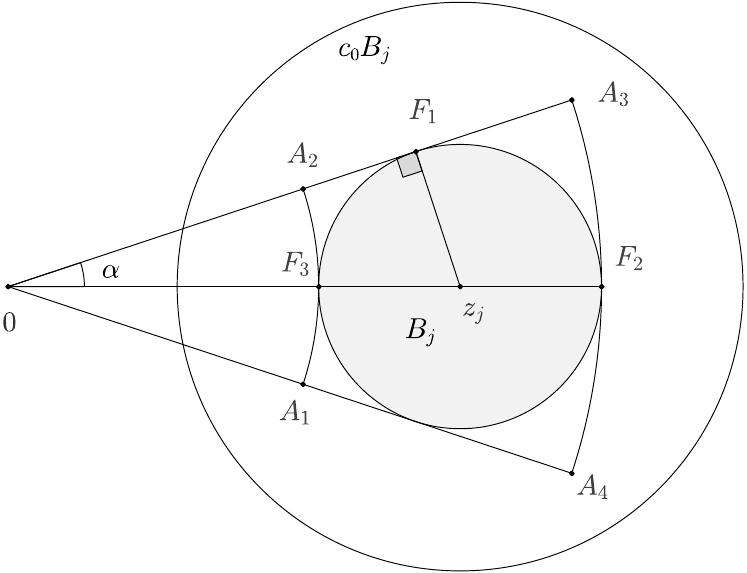}}
\caption{Neighbourhood of hole $B_j$}
	\label{fig:one_concentr}
\end{figure}

\noindent {\bf Proof.}
Consider disk $B_j=\bar{\mathcal B}(z_j,r_j),\, j\neq 1$. We have $B_j\subset c_0 B_j \subset c_0^2 B_j$ while disk~$c_0^2 B_j$ does not contain the origin  --  since $0\in B_1$ and disks $c_0^2 B_j$ are disjoint. Let us do a planimetric construction  given at fig.~\ref{fig:one_concentr} relying on the fact that $0\notin c_0 B_j$.

Since $0\notin B_j$, there exist two lines passing through $0$ and touching the circle~$\partial B_j$; also, there exist two circles centered in the origin and touching $\partial B_j$. Denote the intersection point of these tangent lines and circles by $A_1, A_2, A_3$ and $A_4$ as at fig.~\ref{fig:one_concentr}.  Let~$Q_j$ be curved quadrilateral circumscribed around circle $\partial B_j$  with vertices $A_1, A_2, A_3, A_4$ and bounded by arcs $\stackrel \frown{ A_1 A_2}$,  $\stackrel \frown {A_3 A_4}$ and segments $A_2 A_3$, $A_4 A_1$. Then $\varphi(Q_j)$ is a rectangle on~$\Cylr$.

Denote by $F_1$ the point of touching of the line $0A_2A_3$ and the circle $\partial B_j$ and by  $F_3$ and $F_2$ the point of touching of arcs $\stackrel \frown{ A_1 A_2}$ and  $\stackrel \frown {A_3 A_4}$ with $\partial B_j$  respectively.

Let $x=|z_j|$. Our draft is scale-invariant, hence all the ratios and all the angles on this draft depend only on $r_j/x$. At the same time, $r_j/x\le c_0^{-2}$ (since $0\notin c_0^2 B_j$). Let us prove that $Q_j \subset c_0 B_j$ if $c_0$ is big enough. 
We may assume that ray $0F_3 z_j F_2$ is the semi-axis $[0,+\infty)$, in the other case apply a rotation. 

Notice that angle $\alpha:=\angle A_2 0 F_3=\arcsin(r_j/x)$. Let $z\in Q_j$. Then $|\arg z| \le \alpha$, $|z|\in [x-r_j, x+r_j]$. Hence, the length of the arc of the circle centered in origin and joining $z$ with a point on real axis (more precisely, on the segment~$[F_2 F_3]$)  does not exceed $(x+r_j)\arcsin(r_j/x)$ and thus is no more than 
$$
\frac{\pi}{2}\cdot\frac{r_j}{x}\cdot \left(x+r_j\right) \le \pi r_j,
$$
since $x\ge r_j$. All the points on the segment $F_2 F_3$ lie in $B_j$. Thus $|z - z_j|\le (\pi+1) r_j$ for any $z\in Q_j$. Hence we proved that $Q_j\subset c_0 B_j$ if $c_0>\pi+1$.

\begin{figure} 
\hspace{0.6cm}
{\includegraphics[scale=1.3]
	{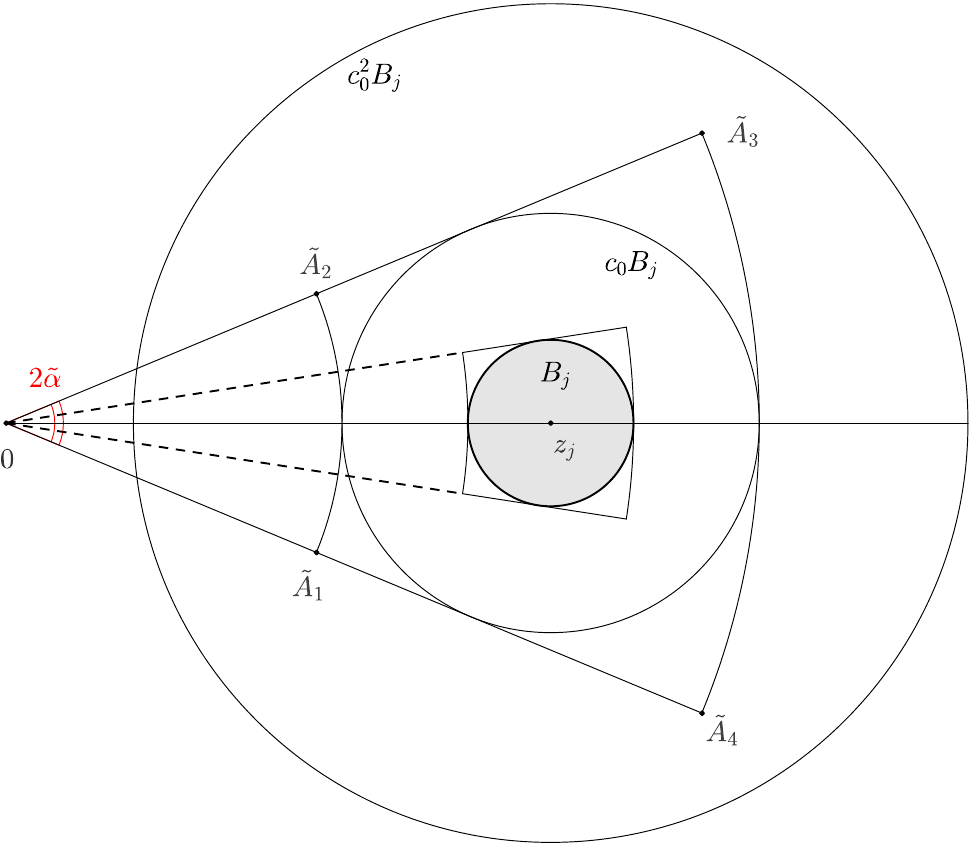}}
\caption{Two neighbourhoods  of hole $B_j$}
\label{fig:concentr_two}
\end{figure}

Now notice that disk $c_0^2 B_j$ does not contain the origin. Thus we can perform an analogous construction for disk $c_0 B_j$ by circumscribing around it curved quadrilateral~$\tilde Q_j$ such that $\varphi(\tilde Q_j)$ is a rectangle on $\Cylr$. The construction is given at fig.~\ref{fig:concentr_two}. Namely, this time we consider two lines tangent to $\partial(c_0 B_j)$ and passing through $0$  and also two circles centered in origin which are also tangent to $\partial(c_0 B_j)$. These four curves bound a curved quadrilateral~$\tilde Q_j$ formed by arcs $\stackrel \frown {\tilde A_1 \tilde A_2}$ and $\stackrel \frown {\tilde A_3 \tilde A_4}$ and also by segments  $\tilde A_2 \tilde A_3$ and $\tilde A_1 \tilde A_4$. By $\tilde \alpha$ we denote the half of the angle under which disk $c_0 B_j$ is seen from point $0$. Notice that $\tilde{\alpha} = \arcsin(c_0 r_j/x)$ -- analogously to the previous construction. As in the above, we prove that $\tilde Q_j \subset c_0^2 B_j$ if $c_0$ is big enough.


\begin{figure} 
	\hspace{2.45cm}
	{\includegraphics[scale=1]
		{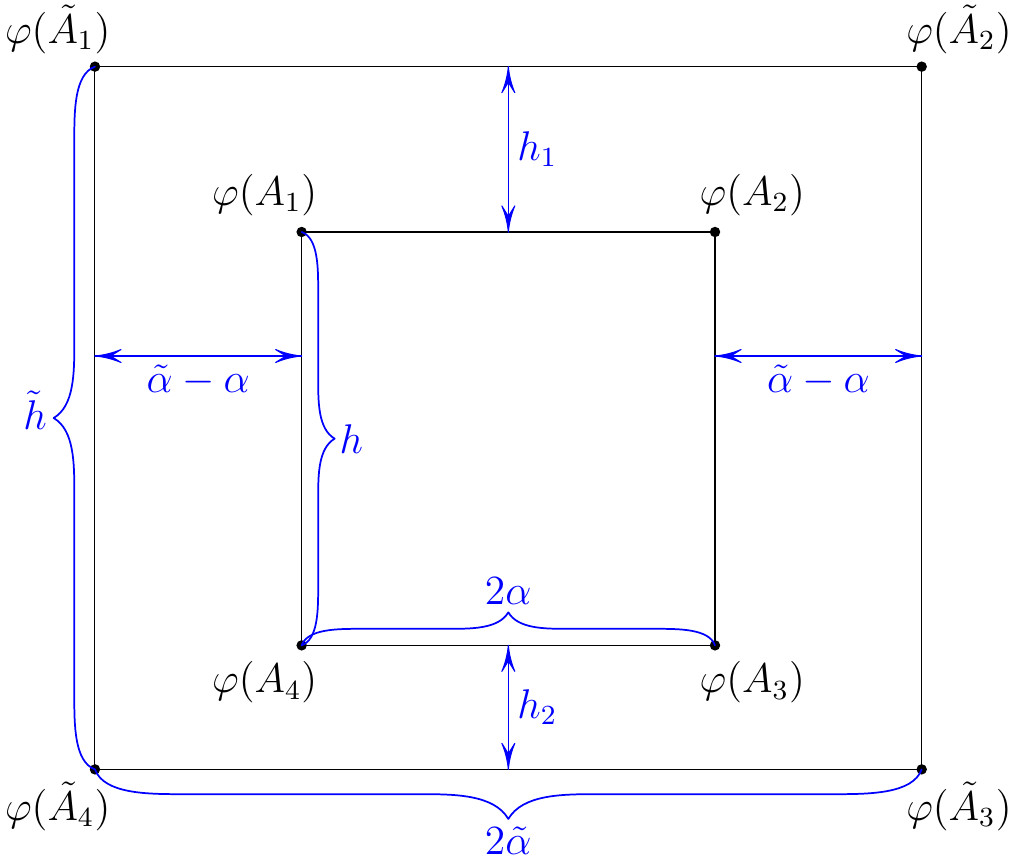}}
	\caption{Images of neighbourhoods of hole $B_j$}
	\label{fig:rectan_cyl}
\end{figure}

Under mapping $\varphi$ quadrilaterals  $Q_j$ and $\tilde Q_j$ pass to rectangles on cylinder (that is, rectangles in one of the charts $(\theta, y)$), geometric parameters of their mutual placement are quantities $\alpha, \tilde{\alpha}$ which we already introduced  and also $h, \tilde h, h_1, h_2$ pointed out at fig.~\ref{fig:rectan_cyl}. All these parameters depend only on $r_j/x$. Moreover, if 
$c_0>2$ then for $x>c_0^2 r_j$ the following asymptotic estimates are held (the constants of comparability are absolute and do not depend on $c_0$ for  $c_0>2$):

\begin{enumerate}
	\item $2\alpha=2\arcsin(r_j/x)\asymp  r_j/x$ and $h=\log\left(\dfrac{x+r_j}{x-r_j}\right)\asymp r_j/x$; 
	\item $\tilde \alpha - \alpha =\arcsin(c_0 r_j/x)-\arcsin(r_j/x)\asymp (c_0-1)r_j/x$; 
	\item $h_1=\log(x-r_j)-\log(x-c_0 r_j)\asymp (c_0-1)r_j/x$;
	\item $h_2=\log(x+c_0 r_j)-\log(x+r_j)\asymp (c_0-1)r_j/x$.
\end{enumerate}
If factor $c_0-1$ in the last three estimates is big enough (in account of comparability constants in all the four estimates) then distances from $\varphi(Q_j)$ to $\partial \varphi(\tilde Q_j)$ are large enough with respect to the lengths of sides of rectangle $\varphi(Q_j)$. In this case we may take any square with side $\max\{2\alpha, h\}$ containing $\varphi(Q_j)$ as $Q_j'$. Lemma is proved. $\blacksquare$

\medskip

Now we proceed the proof of theorem~\ref{patch_theorem}. We will correct the function $\tilde u_1=u_1\circ \varphi^{-1}$ on cylinder by changing its values on squares $2Q_j'$ constructed in the lemma just proved. Namely, let us find, for every $j=2,3,\dots$, a smooth function $v_j\colon 2Q_j'\to \R$ with the following properties:

\begin{enumerate}
	\item $v_j=\tilde u_1$ on $\partial(2Q_j')$;
	\item $v_j\equiv a_j=\const$ on $Q_j$ with constants $a_j\in\R$ as we will.
\end{enumerate}
Moreover, we need Dirichlet integrals $\int\limits_{2Q_j'}|\nabla v_j|^2\,d\theta\,dy$ to be able to estimate. For this goal it is natural to require that $\Delta v_j=0$ in $2Q_j'\setminus Q_j$ and that $a_j$ be such that Dirichlet integral of solution of the corresponding  Dirichlet problem will be minimal over all $a_j\in \R$. However, we will avoid such a boundary problem.

Notice -- and this is the keypoint of the current proof! -- that $\tilde u_1(\theta,y)=y/|\log r|$ is a linear function. Thus the problems on finding functions~$v_j$ for different $j$ are obtained from each other by linear change of variables in domain of arguments and in the domain of values (together with boundary data). 

In the plane $(\theta, y)$ pick a smooth function $v_0\colon Q_0=[-1,1]\times [-1,1] \to \R$ such that $v_0(\theta,y)=y$ for $(\theta,y)\in {\partial Q_0}$ and $v_0=0$ on $\frac12 Q_0:=[-1/2,1/2]\times [-1/2,1/2]$. The choice of $v_0$ depends on nothing. Let $Q_j'$ be square with side $\Lambda_j$ and centered in  $(\theta_j, y_j)\in\Cylr$. Put
$$
v_j(\theta,y):=\frac{\Lambda_j}{|\log r|}\, v_0\left(\frac{\theta-\theta_j}{\Lambda_j}, \frac{y-y_j}{\Lambda_j}\right)+\frac{y_j}{|\log r|}, ~~ (\theta, y) \in 2Q_j'.
$$
In is easy to see that then $v_j=\tilde u_1$ on $\partial (2Q_j')$, moreover, $v_j\equiv y_j/{|\log r|}$ on $Q_j'$.

Put 
$$
\tilde u(\theta, y) = 
\begin{cases}
v_j(\theta, y), &(\theta,y)\in 2 Q_j', \, j=2,3,\dots;\\
\tilde u_1(\theta, y), &(\theta, y) \in \Cylr \setminus \bigcup\limits_{j=2}^\infty 2Q_j'.
\end{cases}
$$
Consider pull-back $u=\tilde u\circ\varphi$ of function $\tilde u$ from $\Cylr$  into $\Omega\setminus B_1$ and define this function to be equal $1$ on $B_1$. Let us show that $u\in\mao$. On compact subsets in $\mathbb D\setminus B_1$ the function $u$ is piecewise-smooth and continuous, since the function $\tilde u$ is piecewise smooth and continuous on $\Cylr$; moreover, $u$ is continuous up on $\partial B_1$. Thus $u\in W^{1,2}_{\loc}(\mathbb D)$. Let us estimate Dirichlet integral of $u$. Since $u\equiv 1$ on $B_1$, then $\int\limits_{\mathbb D} |\nabla u|^2 \dl=\int\limits_{\Cylr} |\nabla\tilde u|^2\,d\theta\,dy$  due to conformal invariance of Dirichlet integral (gradient~$\nabla \tilde u$ is taken by $\theta$ and $y$). We have 
\begin{multline*}
\int\limits_{\Cylr} |\nabla\tilde u|^2\,d\theta\,dy \le 
\int\limits_{\Cylr} |\nabla\tilde u_1|^2\,d\theta\,dy+\sum\limits_{j=2}^\infty\int\limits_{2Q_j'} |\nabla v_j|^2\,d\theta\,dy=\\=\frac{2\pi}{|\log r|}+\sum\limits_{j=2}^\infty\frac{\Lambda_j^2}{|\log r|^2}\cdot \int\limits_{Q_0}|\nabla v_0|^2\,d\theta\,dy
\end{multline*}
(the last equality is true by invariance of Dirichlet integral under scalings). Now notice that the area of $\Cylr$ equals $2\pi|\log r|$. Squares $Q_j'$ are disjoint and thus their total area $\sum\limits_{j=2}^\infty \Lambda_j^2$ does not exceed $2\pi|\log r|$. Thus 
$$
\int\limits_{\Cylr} |\nabla\tilde u|^2\,d\theta\,dy \le \frac{C}{|\log r|},
$$
with absolute constant $C>0$ (since it depends only on the choice of $v_0$). Convergence of the series $\sum\limits_{j=2}^\infty \Lambda_j^2$ implies that $\Lambda_j\xrightarrow{j\to +\infty}0$. From the construction of $\tilde u$ we thus immediately have the continuity of $u$ up to $\partial \mathbb D$ and boundary values of $u$ on $\partial \mathbb D$ are zero. Hence, $u\in \sob$.

So, function $u=\tilde u\circ \varphi$ is admissible for $\Omega$ and  takes value $1$ on $B_1$ and value $y_j/|\log r|$ on~$B_j$, $j=2,3,\dots$ (since $\varphi(B_j)\subset Q_j'$). If domain $\Omega$ possesses  interpolation property then, by the second assertion of proposition~\ref{criterii_adm},
\begin{equation}
\label{eq:Blaschke_preliminary}
\frac{C}{|\log r|} \ge C_I^{-2}(\Omega)\cdot \left(1+\sum\limits_{j=2}^\infty \frac{y_j^2}{|\log r|^2}\right),
\end{equation}
from which $r$ is bounded from the below by a quantity depending only on $C_I(\Omega)$, $\sup\limits_{j\in\mathbb N} \diam_H(B_j)$ and $\eps$ (dependence on $\sup\limits_{j\in\mathbb N} \diam_H(B_j)$ and $\eps$ arose in the application of a quasiconformal mapping). Theorem is proved. $\blacksquare$

\medskip

\noindent {\bf Remark.} Inequality~(\ref{eq:Blaschke_preliminary}) implies that if $\Omega$ has interpolation property then the series $\sum_{j=2}^\infty y_j^2$ is convergent. If quantities $\diam_H(B_j)$ are bounded from the above then $y_j$ is comparable with $1-|z_j|$ where $z_j$ is the center of disk $B_j$. Hence, if the hypothesis of theorem~\ref{patch_theorem} is satisfied then series $\sum_{j=1}^\infty (1-|z_j|)^2$ converges, moreover, its sum may be estimated through $C_I(\Omega)$,  $\sup\limits_{j\in\mathbb N} \diam_H B_j$ and the constant $\eps$ of strong separatedness. This is an analogue of Blaschke condition, we will discuss this condition in subsection~\ref{subsection:Blaschke} under more general assumptions.

\subsection{Some examples}
\label{section:examples}

\paragraph{Round annulus.} 
\label{example:annulus}
Let $B_1=\mathcal B(0, r)$, $0<r<1$ and $\Omega=\mathbb D \setminus B_1$ be annulus. It is easy to find form $\omega\in\lo$ with $\Per_1\omega=1$ which minimizes $\int_{\Omega}^{}|\omega|^2\dl$ over all such forms: 

\begin{wrapfigure}[14]{r}{0.4\textwidth}
		\centering\vspace{-2mm}
	{\includegraphics[scale=1.2]
		{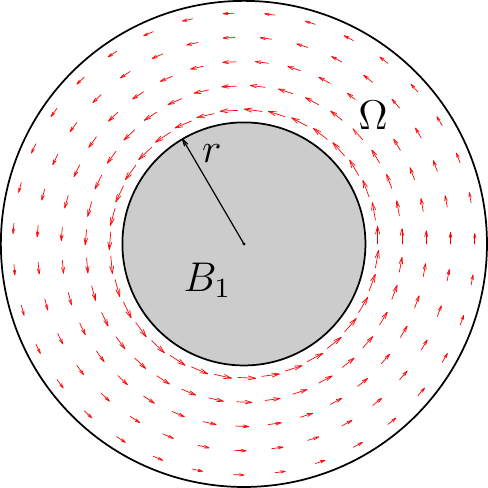}}
	\caption{Round annulus}
	\label{fig:annulus}
\end{wrapfigure}

$$
\omega(x+iy)=\frac{x\,dy-y\,dx}{2\pi\left(x^2+y^2\right)}, ~~ x+iy \in \Omega.
$$
\noindent Then $\int_{\Omega}|\omega|^2\dl=\dfrac{|\log r|}{2\pi}$. Hence 
$$
C_B^{-1}(\Omega)=C_I(\Omega)=\sqrt{\dfrac{|\log r|}{2\pi}}.
$$

\noindent Thus if $r\to 1$ then $C_B(\Omega) \to\infty$. In this case the annulus is too narrow, and then a form with unit period may have a small $L^2$-norm. If, to the opposite, $r\to 0$  then $C_I(\Omega) \to\infty$. In this case a form with unit period is forced to have large $L^2$-norm. If $r=0$ then  in the degenerated annulus $\mathbb D\setminus \{0\}$ there is no closed square-integrable form with non-zero period along a loop winding around the origin. In other words, $L^2$-(co)homologies of the punctured disk are trivial (this property can be understood as removability of the point). That is why we consider domains with holes consisting of more than one point.

\begin{wrapfigure}[20]{l}{0.5\textwidth}
	\centering
	{\includegraphics
		[trim=0.9cm 17.6cm 12.8cm 4.3cm, clip=true, scale=1]
		{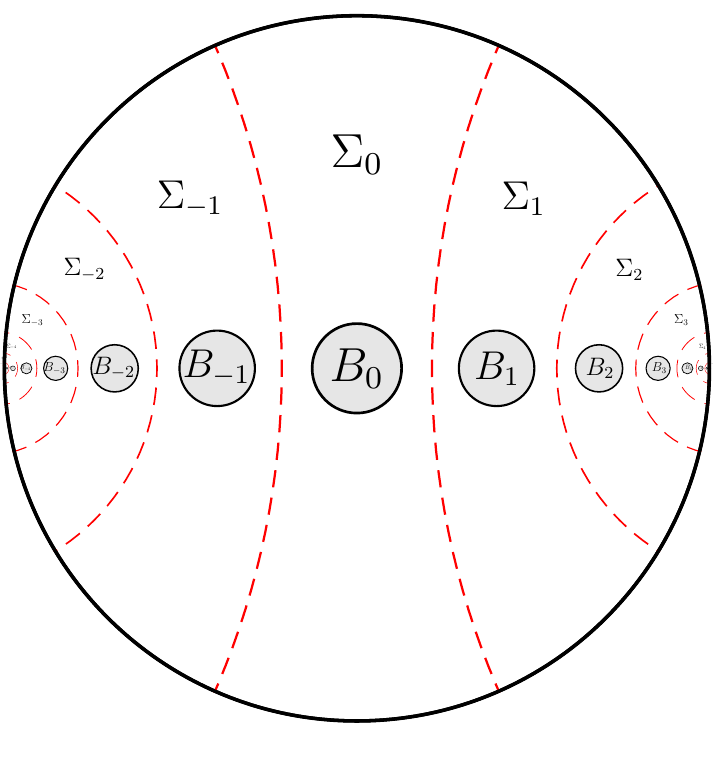}}
	\caption{A domain with conformal symmetry}
	\label{fig:ConformalInvariance}
\end{wrapfigure}

\paragraph{A domain with conformal symmetry.} 
\label{example:conformal_symmerty} 
Our next example is a domain $\Omega$ having a conformal symmetry. Let 
$a\in(-1,0)$ and $\varphi(z) = \frac{z-a}{1-az}$ be  conformal automorphism of 
disk  $\mathbb D$. There exists a domain $\Sigma_0\subset\mathbb D$ such that sets 
$\Sigma_j:=\varphi^j(\Sigma_0)$ 
are pairwise disjoint when $j=0, \pm 1, \pm 2, \dots$, whereas  $\mathbb D \subset \bigcup\limits_{j=-\infty}^\infty \clos \Sigma_j$. Such a set $\Sigma_0$ is called \emph{a fundamental domain} of mapping $\varphi$; one can choose set $\Sigma_0$ such that it will be bounded by four arcs of circles and be symmetric with respect to the origin and also such that $0\in\Sigma_0$ (see fig.~\ref{fig:ConformalInvariance}).

Take a closed disk $B_0$ centered in origin lying strictly inside $\Sigma_0$. Put $B_j = \varphi^j(B_0)$, $\, j=\pm 1, \pm 2, \dots$, and, further,  $\Omega := \mathbb D \setminus \bigcup\limits_{j=-\infty}^\infty \varphi^j(B_0)$. 

All the sets $B_j$ are disks obtained from each other by an isometry of the hyperbolic space. Thus, their hyperbolic diameters coincide. It is easy to prove that sets $B_j$ are $\eps\mbox{-strongly}$ separated with some $\eps>0$. Thus, by the proposition~\ref{predl:strong_sep_Bessel_suff} and theorem~\ref{th:disks_interp_sufficient}, domain $\Omega$ possesses complete interpolation property.

The estimates for complete interpolation property of domain $\Omega$ can be worked out directly. To prove Bessel property notice that if $\omega\in\lo$ then
$$
|\Per_j\omega|^2 \le C\cdot \int\limits_{\Sigma_j\setminus B_j}|\omega|^2\dl
$$	
for any $j\in\mathbb Z$ and constant $C<+\infty$ does not depend on $\omega\in \lo$ and on $j$ due to conformal invariance of the last estimate. Summation over all $j\in \mathbb Z$ gives the conclusion
$$
\sum_{j=-\infty}^\infty|\Per_j\omega|^2 \le C\cdot \int\limits_{\Omega}|\omega|^2\dl.
$$	
Thus $C_B(\Omega) < +\infty$.

To prove interpolation property of $\Omega$ without theorem~\ref{th:disks_interp_sufficient}, let us pick a smooth form $\omega_0\in L^{2,1}_c(\mathbb D \setminus B_0)$ with period $1$ around hole $B_0$. For $j=0, \pm 1, \pm 2, \dots$ define forms $\omega_j=(\varphi^{-j})^{\sharp}\omega_0$,  pull-backs of $\omega_0$ under action of degrees of $\varphi$. Take arbitrary  $a=\{a_j\}_{j=-\infty}^\infty\in\ell^2$. Consider the series
\begin{equation}
\label{eq:explicit_interp_symmetry}
\omega = \sum\limits_{j=-\infty}^\infty a_j \omega_j.
\end{equation}
Since $\Per_j\omega_k=\mathds 1_{\{j=k\}}$ then $\Per \omega=a$ if, of course, series~(\ref{eq:explicit_interp_symmetry}) converges in $\lo$. To prove this and also to estimate $\|\omega\|_{\lo}$ via $\|a\|_{\ell^2}$, it enough to prove that Gram matrix of system $\{\omega_j\}_{j=-\infty}^\infty$ defines a bounded linear operator in $\ell^2$. For this aim, notice that $\langle{\omega_j, \omega_k}\rangle_{\lo}$ depends only on $|j-k|$ due to conformal invariance of our construction. Put $I_j := \int_{\Sigma_j} |\omega_0|^2 \dl$. Since quantities $\diam \Sigma_j$ decrease exponentially when $|j|\to\infty$, the estimate of the form $I_j \le C_1 \cdot \exp(-C_2 |j|)$ is held ($C_1, C_2\in (0,+\infty)$ do not depend on $j$). We have the estimate $|\langle{\omega_0, \omega_k}\rangle_{\lo}|\le \sum_{j=-\infty}^{\infty}\sqrt{I_j I_{k-j}}$; dividing the last series into intervals from $-\infty$ up to  $[k/2]$ and from $[k/2]+1$ up to~$+\infty$ and applying Cauchy-Schwartz inequality, we derive an estimate of the form $|\langle{\omega_0, \omega_k}\rangle_{\lo}|\le C_3 \cdot \exp(-C_4 |k|)$ with  $C_3, C_4\in (0,+\infty)$ not depending on $k$. So,  $|\langle{\omega_j, \omega_k}\rangle_{\lo}|\le C_3 \cdot \exp(-C_4 |j-k|)$, now, by Young's inequality for convolutions, Gram matrix of system $\{\omega_j\}_{j=-\infty}^\infty$ gives a bounded linear operator in $\ell^2$.

Note that complete interpolation property can be proved for other domains having more plentiful  group of conformal automorphisms (see also dyadic example at the end of subsection~\ref{subsection:Blaschke}).

\paragraph{Inverse domain.} 
\label{example:inverse}
The following construction  (and also its modification from example~\ref{example:inversion_composition} below) delivers a counterexample to several naive conjectures on our problem. Also, in this example we, the only time in this paper, show formally  how to construct a domain with countable number of holes from a countable family of holes with finite number of holes by an application to them of a M{\"o}bius automorphisms making conformal copies of these domain to be far from each other in hyperbolic metric and thus almost independent in the sense of interpolation.

Let us show that there exists a domain $\Omega$ having complete interpolation property and such that holes $B_j$ in $\Omega$ are disks but such that numbers $\diam_H(B_j)$ take arbitrarily small values and sets $B_j$ are not strongly separated (proposition~\ref{predl:strong_sep_Bessel} and theorem~\ref{patch_theorem} imply that under complete interpolation condition any of the last two properties implies the other).

\begin{figure}[!b] 
\hspace{-1.4cm}
{\includegraphics
{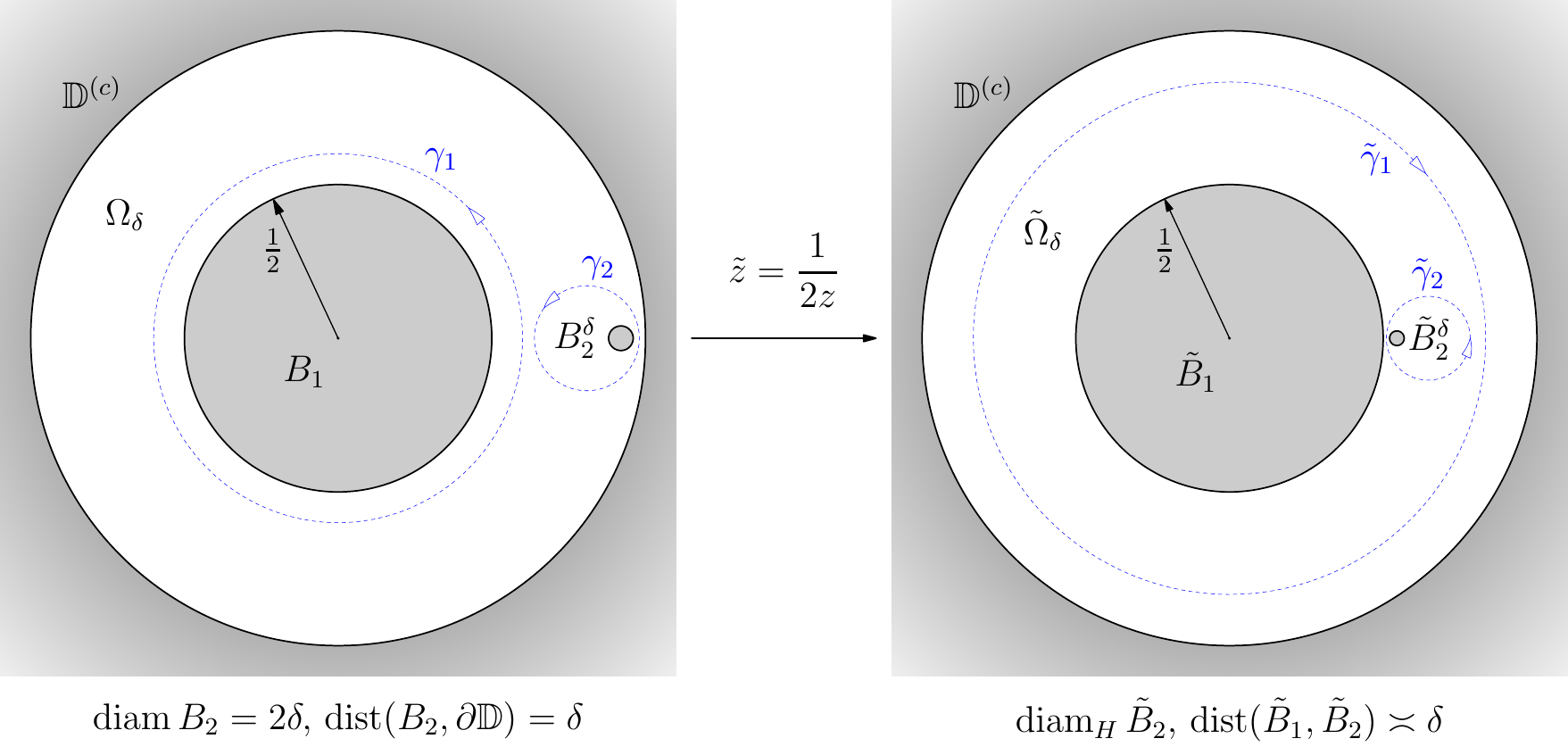}}
\caption{Inverse domain}
\label{fig:Inverse}
\end{figure} 

Pick $\delta>0$ small enough. Consider two round holes $B_1:=\mathcal B(0,1/2)$, $B_2^\delta:= \mathcal B(1-2\delta, \delta)$.
Holes $B_1, B_2^\delta$  are $\eps$-strongly separated with some $\eps>0$, not depending on $\delta$. Moreover, their hyperbolic diameters are bounded from the above and separated from zero uniformly by $\delta$. Put $\Omega_\delta =\mathbb D\setminus \left(B_1\cup B_2^\delta\right)$ (the left-hand domain at fig. \ref{fig:Inverse}). By  proposition~\ref{predl:strong_sep_Bessel_suff} and theorem~\ref{th:disks_interp_sufficient}, domain $\Omega_\delta$ has complete interpolation property, moreover, constants $C_B(\Omega_\delta)$ and $C_I(\Omega_\delta)$ are bounded from the above by some number not depending on $\delta$.

Now, let $\tilde \Omega_\delta$ be a domain obtained by application of inversion $\psi(z)=1/2z$ to~$\Omega_\delta$. Domain $\tilde \Omega_\delta$ is a regular domain in disk $\mathbb D$, denote by $\tilde B_1$ and $\tilde B_2^\delta$ the holes in $\tilde \Omega_\delta$ as at fig.~\ref{fig:Inverse}. Let loops $\gamma_1$ and $\gamma_2$ in domain $\Omega_\delta$ be winding around holes $B_1$ and $B_2^\delta$ respectively, put $\tilde \gamma_1=\psi(\gamma_1)$ and  $\tilde \gamma_2=\psi(\gamma_2)$. By proposition~\ref{predl:conformal_quasiconformal_invariance},  interpolation problem in $\Omega_\delta$ for periods along loops $\gamma_1$ and $\gamma_2$ is conformally equivalent to interpolation problem in domain $\tilde \Omega_\delta$ for periods along loops $\tilde \gamma_1$ and $\tilde \gamma_2$, that is, Bessel  and interpolation constants in these problems coincide. If $\Per_1^{\left(\tilde \Omega_\delta\right)}$ and $\Per_2^{\left(\tilde \Omega_\delta\right)}$ are functionals of periods of a form from $L^{2,1}_c(\tilde \Omega_\delta)$ around holes $\tilde B_1$ and $\tilde B_2$ respectively then
\begin{equation}
\label{eq:period_change_basis}
\left( \begin{array}{c}
\Per_1^{\left(\tilde \Omega_\delta\right)}\tilde \omega\\
\Per_2^{\left(\tilde \Omega_\delta\right)}\tilde \omega
\end{array}
\right) = 
\left(\begin{array}{cc}
-1 & -1\\
0 & 1
\end{array}\right)
\left(\begin{array}{c}
\int_{\tilde \gamma_1}\tilde \omega\\
\int_{\tilde \gamma_2}\tilde \omega
\end{array}\right), ~~ \tilde \omega \in L^{2,1}_c(\tilde \Omega_\delta).
\end{equation}
The $2\times2$ matrix in the right-hand side of the last relation is bounded and non-degenerate. Hence Bessel  and interpolation constants for the problem in domain $\tilde \Omega_\delta$ and periods around holes $\tilde B_1$ and $\tilde B_2$ also stay bounded when $\delta\to 0$. At the same time, for any fixed $\eps>0$ holes in $\tilde\Omega_\delta$ will not be $\eps$-strongly separated if number $\delta$ is small enough, also  $\diam_H \tilde B_2^\delta \xrightarrow{\delta\to 0}0$. 

\medskip

Now let us show how to construct a domain $\Omega$ possessing complete interpolation property and such that hyperbolic diameters of holes in this domain take arbitrarily small values. Also, holes in domain $\Omega$ will not be strongly separated. To do this, apply conformal shifts to domains $\tilde\Omega_\delta$ that we constructed above. Namely, let $\varphi_s(z) = \dfrac{z-s}{1-z\bar s}$ be a  M{\"o}bius transform ($s\in \mathbb D$). \emph{For any sequence $\{\delta_n\}_{n=1}^\infty\subset(0,1/100)$ there exists a real sequence $\{s_n\}_{n=1}^\infty\subset(-1,0)$ tending to $-1$ rapidly enough such that  sets $\varphi_{s_n}\left(\tilde B_1^{\delta_n}\cup \tilde B_2^{\delta_n}\right)$ are pairwise disjoint for $n=1,2,\dots$ while domain $\Omega = \mathbb D \setminus \bigcup\limits_{n=1}^\infty\varphi_{s_n}\left(\tilde B_1^{\delta_n}\cup \tilde B_2^{\delta_n}\right)$ has complete interpolation property.}

First we show how to provide Bessel property of domain $\Omega$ constructed in such a manner. For $\delta\in(0,1/100)$ it is easy to show that if $\omega$ is a closed form in domain $\tilde\Omega_\delta$ then
$$
\left\|\Per^{\left(\tilde \Omega_\delta\right)} \omega\right\|_{\ell^2}^2 \le C_1 \cdot  \int\limits_{\tilde\Omega_\delta \cap \mathcal B(0,3/4)} |\omega|^2\dl.
$$
with some constant $C_1 < +\infty$ not depending on $\delta$. (To prove this, 
it is enough to enlarge the hole  $B_1^\delta$ in domain $\Omega_\delta$ such that its radius will become equal to $2/3$ and apply proposition \ref{predl:strong_sep_Bessel_suff}.) By conformal invariance, an analogous estimate will be true also for a set of the form $\varphi_s(\tilde\Omega_\delta \cap \mathcal B(0,3/4))$ for periods around holes $\varphi_s(\tilde 
B_1^\delta)$ and  $\varphi_s(\tilde B_2^\delta)$ (for any $s\in\mathbb D$). Now it is sufficient to 
choose numbers $s_n$, $n=1,2,\dots$, consecutively in order to have pairwise disjointness of 
sets $\varphi_{s_n}\left(\tilde\Omega_{\delta_n} \cap \mathcal B(0,3/4)\right)$, this will 
immediately lead us to the sequence of periods in domain $\Omega$.

Now let us show how to provide interpolation property of domain $\Omega$. For any $\delta\in(0,1/100)$ there exist forms $w_{1,\delta}, \, w_{2,\delta} \in L^{2,1}_c(\tilde \Omega_\delta)$ for which $\Per^{\left(\tilde \Omega_\delta\right)} w_{1,\delta} = (1,0)$,  $\Per^{\left(\tilde \Omega_\delta\right)} w_{2,\delta} = (0,1)$, whereas $\|w_{1,\delta}\|_{L^{2,1}_c(\tilde \Omega_\delta)}\le C_2$, $\|w_{2,\delta}\|_{L^{2,1}_c(\tilde \Omega_\delta)} \le C_2$ (constant $C_2 < +\infty$ does not depend on $\delta$). That is because interpolation constants of domains~$\tilde\Omega_\delta$ are bounded uniformly for $\delta < 1/100$. Put $w_1^n:= (\varphi_{s_n}^{-1})^\sharp w_{1,\delta_n}$,  $w_2^n:= (\varphi_{s_n}^{-1})^\sharp w_{2,\delta_n}$. Then $w_1^n, w_2^n \in \lo$. If $a=(a_1^1, a_2^1, a_1^2, a_2^2, \dots)\in \ell^2$ is a real sequence and
$$
\omega = \sum\limits_{j=1}^\infty (a_1^j w_1^j+a_2^j w_2^j),
$$
then $\Per^{(\Omega)} \omega = a$. To estimate $\|\omega\|_{\lo}$ through $\|a\|_{\ell^2}$ (and to prove $L^2$-convergence of the series defining this form) one has to prove that Gram matrix of system $\{w_1^n, w_2^n\}_{n=1}^\infty$ defines a bounded linear operator from $\ell^2$ into $\ell^2$ (see~\cite{Bari}). Let us choose numbers $s_n$ consecutively for $n=1,2,\dots$. If  numbers $s_1, \dots, s_{n-1}$ are chosen and  $s_n\to -1$ then for $\alpha, \beta=1,2$ and $m=1,\dots, n-1$ it is easy to prove that 
$$
\int\limits_{\varphi_{s_m}(\tilde \Omega^{\delta_m})\cap \varphi_{s_n}(\tilde \Omega^{\delta_n})} |w_\alpha^m|\cdot |w_\beta^n|\dl \to 0.
$$
Hence numbers $s_n$ can be chosen such that scalar products $\langle w_\alpha^n, w_\beta^m\rangle_{\lo}$ ($m\neq n, \, \alpha, \beta=1,2$) will be arbitrarily small in any norm. In particular, under appropriate choice of sequence $\{s_n\}_{n=1}^\infty$ Gram matrix of system $\{w_1^n, w_2^n\}_{n=1}^\infty$ can be forced to be arbitrarily close in Hilbert-Schmidt norm to the matrix formed by blocks 
$$
\left(
\begin{array}{cc}
\langle{w_1^n, w_1^n}\rangle_{\lo}& \langle{w_1^n, w_2^n}\rangle_{\lo}\\
\langle{w_2^n, w_1^n}\rangle_{\lo}& \langle{w_2^n, w_2^n}\rangle_{\lo}
\end{array}
\right).
$$
But such a matrix defines a bounded operator due to estimate on $\|w_{1,\delta}\|_{L^{2,1}_c(\tilde \Omega_\delta)}$ and $\|w_{2,\delta}\|_{L^{2,1}_c(\tilde \Omega_\delta)}$, this implies boundedness of Gram matrix of system $\{w_1^n, w_2^n\}_{n=1}^\infty$. Thus, under an appropriate choice of numbers $s_n, \, n=1,2,\dots,$ domain $\Omega$ will have interpolation property.

It is easily seen that if $\delta_n \xrightarrow{n\to\infty}0$ then  $\diam_H \varphi_{s_n}(\tilde B_2^{\delta_n})\xrightarrow{n\to\infty}0$ and
$$
\frac{\dist(\varphi_{s_n}(\tilde B_1^{\delta_n}), \varphi_{s_n}(\tilde B_2^{\delta_n}))}{\diam \varphi_{s_n}(\tilde B_1^{\delta_n})} \xrightarrow{n\to\infty}0.
$$

Thus, \emph{there exist regular domains $\Omega$ with round holes possessing complete interpolation property and such that holes in $\Omega$ are not strongly separated, and also hyperbolic diameters of these holes take arbitrarily small values.}

Let us note one more property of the constructed domain $\Omega$. Recall that if $\Omega'$ is some regular domain and  $B_1, B_2, \dots$ are holes in $\Omega'$ then $\mu_{\Omega'}=\sum\limits_{j=1}^\infty \frac{\mathds{1}_{B_j}\cdot \lambda_2}{\lambda_2(B_j)}$ is the measure associated to this domain. Notice that the contribution of hole $\tilde B_2^\delta$ into $\mu_{\tilde \Omega_\delta}$ is weakly${}^*$ close to point unit mass if $\delta\to 0$. Thus it is easy to prove that norm of embedding operator of $\sob$ into $L^2(\mu_{\tilde \Omega_\delta})$ tends to $+\infty$ when $\delta\to 0$; at the same time, a large norm in  $L^2(\mu_{\tilde \Omega_\delta})$ can be given by functions constant on $\tilde B_2^\delta$ and having unit norms in $\sob$. Due to conformal invariance of Dirichlet integral, norms of embeddings of~$\sob$ into $L^2\left(\dfrac{\mathds 1_{B}\cdot \lambda_2}{\lambda_2(B)}\right)$, where $B=\varphi_{s_n}(\tilde B_2^{\delta_n})$, will also tend to infinity when $\delta_n\to 0$ (for any choice of numbers $s_n$). But this implies that \emph{embedding operator from $\sob$ into $L^2(\mu_\Omega)$ will not be bounded}. Thus, complete interpolation property of a  domain $\Omega$ does not, in general, imply (MC) property of measure $\mu_\Omega$ (though estimate $\|u\|_{L^2(\mu_\Omega)}\le C_I(\Omega) \cdot \|u\|_{\sob}$ for functions constant on any hole in $\Omega$ is equivalent to interpolation property of such a domain).

\paragraph{Composition of inversions.} 
\label{example:inversion_composition}
The above-constructed example of an inverse domain leads us to the natural question: \emph{what will be if we make a larger number of inversions?} Given by a number $M\in\mathbb N$ and numbers  $s_1, s_2, \dots, s_M\in(-1,0)$ close to $-1$, let us construct a domain $\Omega_M$ with a step-by-step process. On $m$-th step ($m=1,2,\dots, M$), we will construct a domain $\Omega_m$ whose holes are disks $B_{1,m}, B_{2,m}, \dots, B_{m,m}\subset \mathbb D$; also, we will construct closed oriented curves $\beta_{1,m},\beta_{1,m},\dots,\beta_{m,m}$ in $\Omega_m$. 

Put $B_{1,1}=\bar{\mathcal B}(0,1/2)$, $\Omega_1=\mathbb D \setminus B_{1,1}$ and 
let $\beta_{1,1}$ be a circle placed in~$\Omega_1$ and winding around $B_{1,1}$ 
counter-clockwise. As in the previous example, put $\psi(z)=\dfrac 1{2z}$,  
$\varphi_s(z) =\dfrac{z-s}{1-\bar{s}z}$,  $s,z\in\mathbb D$.  Suppose that 
$m=1,2,\dots, M-1$ and that disks $B_{1,m}, B_{2,m}, \dots, B_{m,m}$ are 
already constructed. For $s_m$ close enough to~$-1$ put 
$B'_{j,m+1}:=\varphi_{s_m}(B_{j,m})$,  $\beta'_{j,m+1}:=\varphi_{s_m}(\beta_{j,m})$ for $j=1,2,\dots,m$, further, $B'_{m+1,m+1}:=\bar{\mathcal B}(0,1/2)$ and $\Omega_{m+1}':=\mathbb D\setminus\bigcup_{j=1}^{m+1} B_{j,m+1}'$; curve $\beta_{m+1,m+1}'$ is defined as a loop in $\Omega_{m+1}'$ winding around $B_{m+1,m+1}'$ counter-clockwise and not winding around other holes in $\Omega_{m+1}'$.  Now put $B_{j,m+1}:=\psi(B'_{j,m+1})$, $\beta_{j,m+1}:=\psi(\beta_{j,m+1}')$ for $j=1,2,\dots, m,$ $\beta_{m+1,m+1}:=\psi(\beta_{m+1,m+1}')$, $B_{m+1,m+1}:=\bar{\mathcal B}(0,1/2)$, $\Omega_{m+1}:=\mathbb D\setminus\bigcup_{j=1}^{m+1} B_{j,m+1}$.

Performing this construction consecutively for $m=1,2,\dots, M-1$, we will obtain domains $\Omega_m, \Omega_m', \, m=1,2,\dots, M$. For such $m$ let us state problem on interpolation by periods in $\Omega_m$  along curves $\beta_1, \beta_2, \dots, \beta_m$, and also this problem in $\Omega_m'$  and along curves $\beta_1', \beta_2', \dots, \beta_m'$ (see at the end of subsection~\ref{subsection:definitions}). 

Take a number $C > \max\{C_B(\Omega_1; \beta_1), C_I(\Omega_1; \beta_1)\}$. Let us prove by induction by $m$  that, under appropriate choice of numbers $s_m$, the inequality $\max\{C_B(\Omega_m; \beta_1, \beta_2, \dots, \beta_m), C_I(\Omega_1; \beta_1, \beta_2, \dots, \beta_m)\}<C$ is held for all $m=1,2,\dots,M$. Indeed, if this inequality is fulfilled for the domain $\Omega_m$ then by use of conformal invariance of the problem and arguing like in example~\ref{example:inverse} it is easy to show that $\max\{C_B(\Omega_{m+1}'; \beta_1', \beta_2', \dots, \beta_m', \beta_{m+1}'), C_I(\Omega_{m+1}'; \beta_1', \beta_2', \dots, \beta_m', \beta_{m+1}')\}<C$ if  $s_m$ is close enough to $-1$. But then $$\max\{C_B(\Omega_{m+1}; \beta_1, \beta_2, \dots, \beta_m, \beta_{m+1}), C_I(\Omega_{m+1}; \beta_1, \beta_2, \dots, \beta_m, \beta_{m+1})\}<C,$$ again due to conformal invariance.

Let $j=1,2,\dots, M$ and  $\gamma_j$ be a loop winding around hole $B_{j,M}$ in $\Omega_M$ in the positive direction and not winding around the other holes in this domain. It easy to see that in the space $H_{1,c}(\Omega_M)$ of compact homologies we have $\beta_1=\gamma_1$, $\beta_j=-\sum_{k=1}^{j}\gamma_k$ for $j=2,3,\dots, M$. (For example, the second relation should be understood as follows: $\int_{\beta_j}\omega = -\sum_{k=1}^{j}\int_{\gamma_k}\omega$ for any $\omega\in L^{2,1}_c(\Omega_M)$.) Let $A$ be $M\times M$ matrix which maps period vector of a form $\omega\in L^{2,1}_c(\Omega_M)$ along curves $\gamma_1, \gamma_2, \dots, \gamma_M$ into vector of $\omega$'s periods along loops $\beta_1, \beta_2, \dots, \beta_M$ (cf. with~(\ref{eq:period_change_basis})). Then norm of matrix $A$ with respect to Euclidean norm in  $\R^M$ is greater than $\sqrt M$. Since $C_B(\Omega_{m+1}; \beta_1, \beta_2, \dots, \beta_M) < C$, this implies that $C_I(\Omega_{m+1}; \gamma_1, \gamma_2, \dots, \gamma_M) \ge \sqrt M/C$. 

So, if in $\Omega_M$ we state interpolation problem for periods around holes in this domain as we usually do, then interpolation constant for such a problem will be large if $M$ is large; one can prove the same about Bessel constant for this problem. This naturally leads us to the following question: \emph{do there exist such domains $\Omega$ that constants $C_B(\Omega)$ and $C_I(\Omega)$ are not very large, but disk $\mathcal B(0,2/3)$ intersects a very large  number of holes in $\Omega$?} The negative answer will be given in section~\ref{section:uniform_local_finiteness}.

The example just constructed also leads to the statement of the problem on integer Riesz basis in Hilbert homologies which will be discussed in subsection~\ref{subsec:open_final}.

\subsection{Partial criterion}

Let us summarize results of proposition~\ref{predl:strong_sep_Bessel_suff}, 
theorem~\ref{th:disks_interp_sufficient}, proposition 
\ref{predl:strong_sep_Bessel}, theorem~\ref{patch_theorem} and example  
\ref{example:inverse}  from subsection~\ref{section:examples} in the following 
theorem:
\begin{theorem}
\label{th:predv_criterii}
	Let $\Omega$ be regular domain with round holes $B_j, \, j=1,2,\dots$. Suppose also that $\sup\limits_{j\in\mathbb N}\diam_H(B_j) < +\infty$.
	\begin{enumerate}
		\item If holes $B_j$ are $\eps$-strongly separated with some $\eps>0$ and also $\inf\limits_{j\in\mathbb N}\diam_H (B_j) > 0$, then  $\Omega$ has complete interpolation property. Moreover, constants $C_B(\Omega)$ and $C_I(\Omega)$ can be estimated from the above only in terms of $\sup\limits_{j\in\mathbb N}\diam_H(B_j),\, \inf\limits_{j\in\mathbb N}\diam_H(B_j)$ and $\eps$.
		
		\item Suppose that $\inf\limits_{j\in\mathbb N}\diam_H(B_j)>0$ and $\Omega$ possesses complete interpolation property. Then holes $B_j$ are $\eps$-strongly separated with some $\eps>0$ depending only on $\inf\limits_{j\in\mathbb N}\diam_H(B_j)$ and $C_B(\Omega)$. 
		
		\medskip
		
		If, vice versa, holes $B_j$ are $\eps$-strongly separated for some $\eps>0$ and $\Omega$ has complete interpolation property, then $\inf\limits_{j\in\mathbb N}\diam_H(B_j)$ is strictly positive and can be estimated from the below only through $\sup\limits_{j\in\mathbb N}\diam_H(B_j), \, \eps$ and $C_I(\Omega)$.
		 
	\item There exist regular domains $\Omega$ of the form $\mathbb D\setminus \bigcup_{j=1}^\infty B_j$ where each $B_j$,  $j=1,2,\dots,$  is a disk  which have complete interpolation property and such that $\inf\limits_{j\in\mathbb N}\diam_H (B_j)=0$ and holes $B_j$ are not  $\eps$-strongly separated for any $\eps>0$.
	\end{enumerate}
\end{theorem}

The first two assertions of this theorem  give criterion of complete interpolation for domains with round holes in the case when one of the two following conditions is satisfied: either $\inf\limits_{j\in\mathbb N}\diam_H(B_j)>0$, or holes $B_j$ are strongly separated. But no one of these two condition is necessary for complete interpolation -- this is the third assertion of the theorem. 

First two assertions of theorem~\ref{th:predv_criterii} are, in fact, true for domains with domains with holes of arbitrary form, not necessary round.  Propositions~\ref{predl:strong_sep_Bessel_suff} and~\ref{predl:strong_sep_Bessel} were proved without assumption that all the sets $B_j$ are disks. Sufficiency of hypothesis of the first assertion of theorem~\ref{th:predv_criterii} for interpolation property was proved in theorem~\ref{th:arbitrary_interp_sufficiency_separated}.  Validity of theorem~\ref{patch_theorem} on lower estimate on hyperbolic diameters also stays true when~$B_j$ are not necessarily disks. This can be derived from theorem~\ref{th:interp_sufficient} (see below) taking in account proposition~\ref{predl:strong_sep_Bessel_suff} and theorem~\ref{th:uniform_local_finiteness} on uniform local finiteness which we pass to right now.

\medskip

\section{Uniform local finiteness of family of holes}\markright{Uniform local finiteness}
\label{section:uniform_local_finiteness}

Theorem~\ref{th:predv_criterii} easily implies that if domain $\Omega$ with round holes has complete interpolation property and if either holes $B_j$ are $\eps$-strongly separated or $\inf\limits_{j\in\mathbb N} \diam_H(B_j) >0$ then for any $R<1$ number of holes $B_j$ intersecting $\mathcal B(0,R)$ can be estimated from the above through $C_B(\Omega), C_I(\Omega)$, $R$ and either $\eps$ or $\inf\limits_{j\in\mathbb N} \diam_H(B_j)$ respectively; due to conformal invariance, the same estimate will be true if we replace disk $\mathcal B(0, R)$ by any its conformal copy obtained by a M{\"o}bius transform. Our goal is to prove the same estimate for arbitrary complete interpolation domains without assumptions neither  on strong separatedness nor on the inequality $\inf\limits_{j\in\mathbb N} \diam_H(B_j) >0$.

\begin{define}
	\label{def:ULF}
	Let us say that family of holes in domain $\Omega$ is \emph{uniformly locally finite} if there exists such $N\in\mathbb N$ that any disk in hyperbolic metric having in this metric radius $1$ intersects no more than $N$ of sets $B_j$. The least such $N$ is called the \emph{constant of  uniform local finiteness} of domain $\Omega$ and is denoted by $N(\Omega)$.
\end{define}

Let $d>0$. It is easy to see that if uniform local finiteness property is satisfied then any set of hyperbolic diameter no more than $d$ intersects no more than $\tilde N$ of sets $B_j$ where $\tilde N<+\infty$ can be estimated from the above through $d$ and $N(\Omega)$. Denote this number by $N(\Omega, d)$ having in mind that it depends only on $d$ and $N(\Omega)$.

\begin{theorem}
	\label{th:uniform_local_finiteness}
	If a regular domain $\Omega$ has interpolation and weak Bessel properties then it possesses uniform local finiteness property. Moreover, the constant $N(\Omega)$ can be estimated from the above by means of only $C_I(\Omega)$ and $\tilde C_B(\Omega)$.
\end{theorem}

If $\tilde C_B(\Omega) < +\infty$ then $C_B(\Omega) < +\infty$ (theorem~\ref{th:weak_strong_bessel}). In the argument below we, nevertheless, will make use only of estimate $\tilde C_B(\Omega) < +\infty$. 

\medskip  

\noindent {\bf Proof of theorem.} Due to conformal invariance of the problem it is sufficient to estimate the number of holes intersecting $\mathcal B(0,1/2)$. 

First we make use of interpolation property of domain $\Omega$. We are going to construct a function $u\in \mao$ in such a way that values $u|_{B_j}$ will admit a lower  estimate while norm $\|u\|_{\sob}$ will not be too large. The estimate from the second assertion of proposition~\ref{criterii_adm} will then lead us to the desired.

\medskip

\emph{Heuristic consideration.} Let us look for a function $u\in\mao$ such that 
\begin{equation}
\label{grad_estim}
|\nabla u| \le 1
\end{equation}
almost everywhere. This will, of course, give an estimate for $\|u\|_{\sob}$.  If~(\ref{grad_estim}) is fulfilled then 
\begin{equation} 
\label{heuristic_neq}
\begin{cases}
&\left|(u|_{B_j}) - (u|_{B_{j'}})\right| \le \dist(B_j, B_{j'}),\\
&\left|(u|_{B_j})\right| \le \dist\left(B_{j}, \partial\mathbb D\right).
\end{cases}
\end{equation}
Consider this system of inequalities as a system with countable numbers of unknown quantities $u|_{B_j}$ (and forget about values of $u$ and $\nabla u$ on $\Omega$ for the moment).
What can be the largest value of $u|_{B_{j_0}}$ under such restriction if 
$u|_{\partial \mathbb D}=0$? One can give a sharp answer. Notice that if $j_0$ 
is a fixed index and $(j_1, j_2, \dots, j_\nu)$ is a finite sequence of indices 
then 
$$
u|_{B_{j_0}} \le \dist\left(B_{j_0}, B_{j_1}\right)+\dist\left(B_{j_1}, B_{j_2}\right)+\dots+\dist\left(B_{j_{\nu-1}}, B_{j_\nu}\right)+\dist\left(B_{j_\nu}, \partial\mathbb D\right).
$$
Put 
$$
u|_{B_{j_0}} = \inf\left\{ \dist\left(B_{j_0}, B_{j_1}\right)+\dist\left(B_{j_1}, B_{j_2}\right)+\dots+\dist\left(B_{j_{\nu-1}}, B_{j_\nu}\right)+\dist\left(B_{j_\nu}, \partial\mathbb D\right)\right\},
$$
where $\inf$ is taken over all finite chains of indices $(j_1, j_2, \dots, 
j_\nu)$ of arbitrary length. It is easy to see that numbers $u|_{B_j}$ defined 
in such a way satisfy~(\ref{heuristic_neq}). The second condition is obvious 
(just take a chain of one index $j$). To check the first condition, let us 
notice that if $j_0, j_0'$ are two indices and  $(j_1, j_2, \dots, j_\nu)$ is 
one of chains in the $\inf$ defining $u|_{B_{j_0}}$ then chain  $(j_0, j_1, \dots, j_\nu)$ can be substituted in the $\inf$ from the definition of $u|_{B_{j_0'}}$.  Hence $u|_{B_{j_0'}}\le u|_{B_{j_0'}}+\dist\left(B_{j_0}, B_{j_0'}\right)$; the symmetric inequality is true as well which implies the first condition in~(\ref{heuristic_neq}). Nevertheless, we have not obtained any lower estimates on numbers $u|_{B_j}$, moreover, function $u$ has not been defined on $\Omega$. The exposed consideration leads us to introducing a special metrics in $\mathbb D$. We have already used this idea in the proof of proposition \ref{criterii_adm}.

\medskip

Now pass to formal argument. Everywhere in this proof we consider polygonal chains in $\mathbb C$ with finite number of vertices; that is, piecewise-linear continuous mappings $\Gamma\colon [0,T]\to \mathbb C, \, T<+\infty$. We also assume that such a polygonal chain $\Gamma$ is endowed with an orientation and is parametrized naturally, that is $|\Gamma'(t)|= 1$ on $[0,T]$ almost everywhere with respect to one-dimensional Lebesgue measure $\lambda_1$. If mapping $\Gamma$ is injective then we call such a chain \emph{simple}. If $\Gamma(0)=z_1, \Gamma(T)=z_2$ then we say that chain $\Gamma$ \emph{joins} points $z_1$ and $z_2$.

In $\mathbb C$, consider a degenerated metric $\mathds 1_{\Omega} |dz|$. For a polygonal chain  $\Gamma$ in  $\mathbb C$ put 
$$
L(\Gamma) := \lambda_1(\Gamma^{-1}(\Omega)).
$$
(Set $\Gamma^{-1}(\Omega)$ is open in $\R$ and is thus Lebesgue measurable.) If 
$\Gamma\colon [0,T]\to \mathbb C$ is a  {simple} polygonal chain  in 
$\mathbb C$ then $L(\Gamma)=\mathcal H^1(\Gamma[0,T]\cap \Omega)$ where symbol 
$\mathcal H^1$ stands for one-dimensional Hausdorff measure, that is length. If 
$\Gamma[0,T]\subset \Omega$ and a polygonal chain $\Gamma$ is simple, then 
$L(\Gamma) = \mathcal H^1(\Gamma[0,T])$;  if, to the opposite,  
$\Gamma[0,T]\subset\mathbb D^{(c)}$ or $\Gamma[0,T]\subset B_j$ for some 
$j=1,2,\dots$ then $L(\Gamma)=0$. 

Define a degenerated metric 
$$
\rho(z_1,z_2) := \inf\left\{L(\Gamma) \colon \Gamma \mbox{ is a  polygonal chain joining } z_1 \mbox{ with } z_2\right\}, ~~~ z_1, z_2\in \mathbb C.
$$
As such a $\Gamma$, we may take line segment $[z_1, z_2]$ joining $z_1$ and 
$z_2$ and endowed with natural parametrization, therefore, 
\begin{equation}
\label{lip_metric}
\rho(z_1, z_2) \le L([z_1, z_2]) \le |z_1-z_2|.
\end{equation} 

Recall that any set $B_j$ is a closure of a domain with boundary smooth enough.  Thus any two points $z_1, z_2\in B_j$ can be joined by a polygonal chain $\Gamma$ lying in $B_j$; for such a chain we have $L(\Gamma)=0$. Thus, $\rho(z_1, z_2)=0$. So,  all the holes $B_j$ collapse into points in the degenerated metric $\rho$.

\medskip

Note that any line segment lying in $\Omega$ will be geodesic (that is a locally shortest curve) in metric $\rho$. Also, all curves lying in some hole or in  $\mathbb D^{(c)}$ will be geodesics as well.

\medskip

For $z\in \mathbb C$ put
\begin{multline*}
u(z) := \dist_{\rho} (z, \mathbb D^{(c)}) =\\= \inf\left\{L(\Gamma) \colon \Gamma \subset \mathbb C \mbox{ is a  polygonal chain joining } z \mbox{ with a point in } \mathbb D^{(c)}\right\}.
\end{multline*}
From~(\ref{lip_metric}) we have that 
\begin{equation}
\label{eq:inner_metric_lip}
u(z_2)\le u(z_1)+|z_1-z_2|
\end{equation}
for any $z_1, z_2\in\mathbb C$. Hence function $u$ is Lipschitz and $|\nabla u| \le 1$ almost everywhere in $\mathbb C$. By definition, $u=0$ on $\partial\mathbb D$. Thus,  $u\in \sob$.  Further,   $\rho(z_1, z_2)=0$ if $z_1, z_2 \in B_j$ for some $j$, and then $u$ is constant on each hole $B_j$. So, $u\in \mao$ and $\|u\|_{\sob}\le \sqrt\pi$. 

Now we are going to prove, relying only on weak Bessel property, that values $u|_{B_j}$ can not be very small.

\begin{lemma}
	\label{rho_estim}
	Suppose that domain $\Omega$ has weak Bessel property, that is $\Cap2(B_j, \Omega^{(c)} \setminus B_j) \le\tilde C_B^2(\Omega) < +\infty$ for any $j=1,2,\dots$ \emph{(}see proposition \ref{necessary_cap}\emph{)}. Then, for any $j=1,2,\dots$, the function $u$ defined in the above satisfies the estimate 
	$$
	u|_{B_j} \ge \frac{\exp(-2\pi \tilde C_B^2(\Omega))}2\cdot {\dist(B_j, \partial\mathbb D)}. 
	$$
\end{lemma}

\noindent {\bf Proof.} 
Let $\Gamma\colon [0,T]\to \mathbb C, \, T<+\infty,$ be a  polygonal chain joining a point from $B_j$ with $\partial\mathbb D$. We have to prove that 
\begin{equation}
\label{estim:geodesic_lower}
\lambda_1\left(\Gamma^{-1}(\Omega)\right) \ge \frac{\exp(-2\pi \tilde C_B^2(\Omega))}2 \cdot \dist(B_j, \partial \mathbb D).
\end{equation}
Let us make some technical assumptions.

We may prove the desired estimate assuming that number of holes in  $\Omega$ is finite. Indeed, supposing that $j=1$, let us erase holes $B_{\nu+1}, B_{\nu+2}, \dots$ and put $\Omega_\nu:=\mathbb D\setminus \bigcup_{j'=1}^\nu B_{j'}$, here $\nu=1,2,\dots$. For any such $\nu$ and for any $j'=1,2,\dots, \nu$ estimate $\Cap2(B_{j'}, \Omega_\nu^{(c)}\setminus B_{j'}) \le \tilde C_B^2(\Omega)$ stays true. If we prove inequalities~(\ref{estim:geodesic_lower}) for domain  $\Omega=\Omega_\nu$ and for all $\nu$ then this inequality for the initial domain  $\Omega$ will be obtained by passing to a limit over $\nu\to\infty$ since $\lambda_1\left(\Gamma^{-1}(\Omega_1)\right)\le T <+\infty$.

We may assume that $\Gamma \cap B_{j'}$ is connected for any ${j'}=1,2,\dots$. Indeed, consider the first (with respect to the motion of polygonal chain $\Gamma$) set $B_{j'}$ which $\Gamma$ intersects twice;  the arc of curve  $\Gamma$ from its first entrance into $B_{j'}$ until its last exit from $B_{j'}$ can be replaced by a naturally parametrized polygonal chain in $B_{j'}$ (recall that we assume that boundary $\partial B_{j'}$ is smooth). Such a replacement does not increase $L(\Gamma)$. Further, if the obtained curve intersects some hole twice, then perform the same procedure with this curve,  and so on. The number of holes in $\Omega$ is finite according to the assumption we made above, then we will make only a finite number of changes with our curve $\Gamma$ and hence it will remain a polygonal chain with finite number of vertices.

Further, if $(t_1, t_2)$ is a maximal interval in $[0,T]$ such that $\Gamma\left((t_1,t_2)\right)\subset \Omega$ then replace arc $\Gamma|_{(t_1, 
t_2)}$ in our curve by line segment $[\Gamma(t_1), \Gamma(t_2)]$. Such a replacement does not increase length $L(\Gamma)$. Make such replacements for all connected components of $\Gamma^{-1}(\Omega)\subset [0,T]$. If necessary, reparametrize the obtained curve naturally. If new curve intersects some hole $B_{j'}$ by a non-connected set then again replace arc of curve from its first entrance into $B_{j'}$ until its last exit from $B_{j'}$ by a  polygonal chain lying in~$B_{j'}$. Thus we may  assume that \emph{$\Gamma([0,T])\cap \Omega$ consists of finite number of line segments whereas all the sets $\Gamma \cap B_{j'}$ are connected, $j'=1,2,\dots$}.

If polygonal chain $\Gamma$ intersects some hole $B_{j'}$ by at least two points then let us assume that all the interior points of arc $\Gamma\cap B_{j'}$ lie in $\Innt B_{j'}$. This can be achieved by an appropriate correction of our chain.

We may also assume that $\Gamma$ is a simple polygonal chain because self-intersections increase $L(\Gamma)$.

Finally, we make one more assumption. Polygonal chain $\Gamma$ ends in the point $\zeta\in\partial\mathbb D$. Attach to $\Gamma$ a line segment of length $5$ starting in $\zeta$ in order to have the end of curve $\Gamma$ to be placed on the distance $6$ from the origin. To the other hand, we may assume that polygonal chain $\Gamma$ starts in a point $w\in \partial B_{j}$. Choose $z\in \Innt B_{j}, \, z\neq w,$ such that line segment $[z, w)$ lies in $\Innt B_{j}$ entirely, and attach this segment to $\Gamma$. Now we assume that $\Gamma$ starts in a point $z\in \Innt B_j$.

\medskip

\begin{figure} 
\hspace{-0.85cm}
	{\includegraphics[scale=0.91]
		{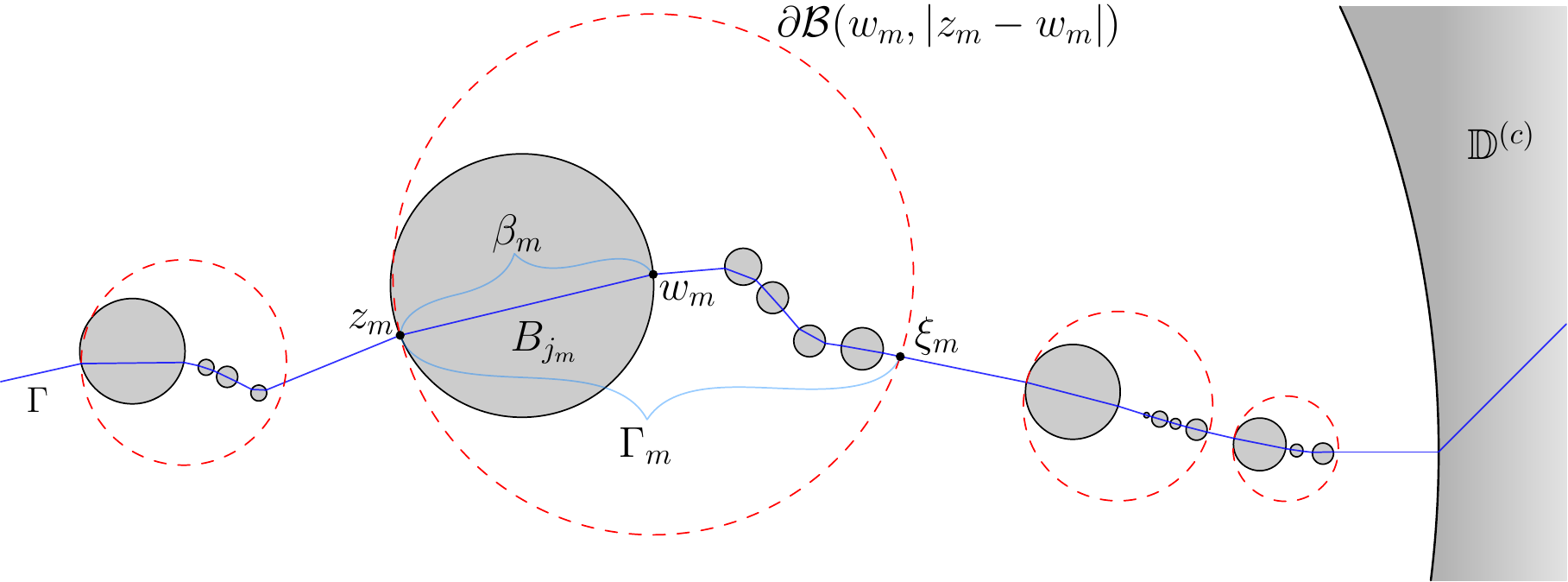}}
	\caption{Subdivision of curve $\Gamma$}
	\label{fig:chain}
\end{figure}

The goal of the argument below is as follows: we construct the covering of set 
$\Gamma\cap \Omega^{(c)}$ by arcs $\Gamma_m \, (m=0,1,\dots)$ of $\Gamma$ in 
such a way that in any of arcs $\Gamma_m$ one of arcs $B_{j'} \cap \Gamma, \, 
j'=1,2,\dots,$ will have a significant weight (see fig.~\ref{fig:chain}). This 
arc will be denoted by  $\beta_m$, its beginning (according to the motion of 
chain $\Gamma$)  by $z_m$, and its end  by $w_m$. Then weight of set 
$(\Gamma_m\cap \Omega^{(c)})\setminus \beta_m$ on arc $\Gamma_m$ will not be 
large (from capacity considerations) while  set $\Gamma_m\cap \Omega$ will, to 
the opposite, be large. This will lead us to the desired estimate. The 
construction will be worked out by a step-by-step process.

\begin{lemma}
	\label{lemma:chain_construct}
	In the settings of lemma~\ref{rho_estim} and assumptions made in the beginning of its proof there exist a number $M \in \mathbb N$, a  sequence of points $z_0, w_0,\xi_0, z_1, w_1,\xi_1, \dots, z_M, w_M, \xi_M$ on curve $\Gamma$, closed arcs $\Gamma_0, \beta_0, \Gamma_1, \beta_1, \dots,\Gamma_M, \beta_M$ on this curve and a sequence of distinct indices $j_0, j_1, \dots, j_M\in\mathbb N$ such that:
	\begin{enumerate}
		\item 
		sequence of points $z_0, w_0,\xi_0, z_1, w_1,\xi_1 \dots, z_M, w_M, 
		\xi_M$ is ordered in the direction of motion of oriented curve $\Gamma$ 
		\emph{(}that is, in order of increase of its parameter\emph{)};
		\item for $m=0,1,\dots, M$ arc $\beta_m$ is the arc of polygonal chain $\Gamma$ starting in $z_m$ and ending in $w_m$; any such arc $\beta_m$ lies entirely in hole $B_{j_m}$; point $w_m$ is the point of exit of $\Gamma$ from $B_{j_m}$;
		\item for $m=0,1,\dots, M$ arc $\Gamma_m$ starts in $z_m$ and ends in 
		$\xi_m$; moreover, $|\xi_m-w_m|=|z_m-w_m|$; arc $\Gamma_m\setminus 
		\beta_m$ does not intersect $B_{j_m}$ \emph{(}but may intersect other 
		holes\emph{)};
		\item distinct arcs $\Gamma_m, \, m=0,1,\dots,M,$ do not intersect by 
		inner points; $\Gamma\cap \left(\mathbb D\setminus \Omega\right) \subset\bigcup_{m=0}^M \Gamma_m$.
	\end{enumerate}
\end{lemma}

\noindent {\bf Proof.}
Let us organize a stepwise process. Curve $\Gamma$ starts at a point in $\Innt B_j$ ($j$ is the index of hole for which we are proving estimate~(\ref{estim:geodesic_lower})). Put $j_0=j$, take the beginning of curve $\Gamma$ for $z_0$,  for $w_0$ we take the point of exit of $\Gamma$ from $B_j$, and $\beta_0$ be the arc of $\Gamma$ from $z_0$ to $w_0$. 

Let $m=0,1,\dots$. Suppose that index $j_m$, arc $\beta_m$ and points $z_m, 
w_m$ are already chosen. Consider disk $\bar{\mathcal B}(w_m, |z_m-w_m|)$. 
Starting from point $w_m$ and moving by curve $\Gamma$ in the direction of 
increase of its parameter, we necessarily escape disk $\bar{\mathcal 
B}(w_m, |z_m-w_m|)$ since the end of polygonal chain $\Gamma$ is situated on 
the distance $6$ from the origin according to our assumptions on $\Gamma$. For 
$\xi_m$ we take the point of first exit of chain $\Gamma$ from $\bar{\mathcal 
B}(w_m, |z_m-w_m|)$; arc $\Gamma_m$ is defined as arc on $\Gamma$ from $z_m$ 
to  $\xi_m$. 

Recall that $\Gamma$ enters and exits any hole $B_{j'}$ no more than once. 
Let's check the following cases.

\begin{enumerate}
	\item 
	$\xi_m\in\Innt B_{j'}$ for some $j'$. Then put $z_{m+1}=\xi_m,\, j_{m+1}=j'$. For $w_{m+1}$ we take point of exit of chain $\Gamma$ from $B_{j'}$. Arc $\beta_{m+1}$ is then defined as a closed arc of chain $\Gamma$ from $z_{m+1}$ until $w_{m+1}$.
	
	\item Point $\xi_m$ lies in  $\Omega$ or at the boundary of one of the 	holes, but arc of chain $\Gamma$ after~$\xi_m$ does not lie in $\Omega\cup\mathbb D^{(c)}$ 	entirely. Let~$B_{j_{m+1}}$ be first of the holes $B_{j'}$ which~$\Gamma$ 	intersects after $\xi_m$. Let us define point $z_{m+1}$ as the point of 	entrance of $\Gamma$ into~$B_{j_{m+1}}$. (If  $\xi_m$ is a point of 	entrance of curve $\Gamma$ into some hole $B_{j'}$ then put $j_{m+1}:=j', 	\, z_{m+1}=\xi_m$.) Let $w_{m+1}$ be the point of exit of $\Gamma$ from 	$B_{j_{m+1}}$. Arc $\beta_{m+1}$ is defined as closed arc of chain $\Gamma$ 	from $z_{m+1}$ to $w_{m+1}$.
	
	\item Chain $\Gamma$ after point $\xi_m$ lies in $\Omega\cup\mathbb D^{(c)}$ entirely. Then we stop our process by putting $M=m$.
	
	\item $\xi_m \in \mathbb D^{(c)}$. In this case we also stop our process and set $M=m$.
\end{enumerate}

Our process must stop, since the number of holes in domain $\Omega$ is finite according to the assumptions which we made in the beginning of proof of lemma \ref{rho_estim}. Accounting the same  assumptions easily allows us to check all the properties of the points and arcs constructed. $\blacksquare$

\medskip

In the setting of construction built up in lemma~\ref{lemma:chain_construct}, the following estimate is true:
\begin{lemma}
	\label{capacity_estim_lemma}
	$
	\mathcal H^1(\Gamma_m \cap \Omega) \ge e^{-2\pi \tilde C_B^2(\Omega)} \cdot |z_m-w_m|
	$
	for any $m=0,1,\dots,M$.
\end{lemma}

\noindent {\bf Proof.} 
Let $E_m$ be the union of holes $B_{j'}$ intersected by $\Gamma_m$ and distinct 
from $B_{j_m}$. If $\xi_m \in \mathbb D^{(c)}$ then add $\mathbb D^{(c)}$ to 
this union. By monotonicity of capacity,
\begin{equation*}
\Cap2(B_{j_m}, E_m)\le \Cap2\left(B_{j_m}, \mathbb D^{(c)} \cup \bigcup\limits_{j\neq j_m} B_j\right)\le \tilde C_B^2(\Omega).
\end{equation*}
This condition implies the existence of such a function $f\in W^{1,2}_{\loc}(\mathbb C)$ that $f=0$ on $B_{j_m}$, $f=1$ on~$E_m$ almost everywhere, but $\int\limits_{\mathbb C} |\nabla f|^2 \dl \le 2\tilde C_B^2(\Omega)$.

Put 
$$
R_m=\left\{r\in (0, |z_m-w_m|)\colon \partial \mathcal B(w_m, r) \cap \Innt E_m \neq \varnothing\right\}.
$$
We may assume that $f$ is absolutely continuous on circle $\partial \mathcal 
B(w_m, r)$ for almost all $r \in R_m$. For all such $r$ circle $\partial 
\mathcal B(w_m, r)$ intersects both sets $\Innt E_m$ and  $\Innt B_{j_m}$ 
(since it must contain at least one inner point of arc $\beta_m$, whereas all 
such points lie in $\Innt B_{j_m}$ according to the assumptions made in the 
beginning of proof of lemma~\ref{rho_estim}). Hence for almost every $r\in R_m$ 
function $f$ takes values $0$ (on $B_{j_m}$) and $1$ (on $E_m$) at circle 
$\partial \mathcal B(w_m, r)$, which gives 
$$
\int\limits_{\partial\mathcal  B(w_m,r)} |\nabla f|^2\,d\mathcal H^1 \ge \frac{1}{\pi r}
$$
by Cauchy-Schwartz inequality. From this we conclude that 
$$
\int\limits_{R_m} \frac{dr}{r} \le 2\pi \tilde C_B^2(\Omega).
$$
This estimate implies that 
$$
\frac{\mathcal H^1(R_m)}{|z_m-w_m|} \le 1-e^{-2\pi \tilde C_B^2(\Omega)}.
$$
Indeed,  the set $X_0=[|w_m-z_m|-\mathcal H^1(X_0), |w_m-z_m|]$ gives minimum to integral $\int\limits_X \frac{dr}{r}$ over all sets with the same length.

We thus proved that 
\begin{equation}
\label{iskl_neq}
\mathcal H^1([0, |w_m-z_m|]\setminus R_m) \ge e^{-2\pi \tilde C_B^2(\Omega)} \cdot |w_m-z_m|.
\end{equation}
Recall that $\Gamma_m$ is a polygonal chain with finite number of vertices. 
Consider set $A_m=\{z\in {\mathcal B}(w_m, |z_m-w_m|) \colon |z-w_m| \notin 
\clos R_m\}$,  this set is a union of finite number of annuli, sum of their 
widths, according to~(\ref{iskl_neq}), is not less than $e^{-2\pi \tilde 
C_B^2(\Omega)} |z_m - w_m|$, and at the same time $A_m \subset \Omega$. Arc 
$\Gamma_m \setminus B_{j_m}$ starts in the common center of these annuli and 
escapes the greatest of them. Thus $\mathcal H^1(\Gamma_m \cap \Omega) \ge  
e^{-2\pi \tilde C_B^2(\Omega)} |z_m - w_m|$, and the desired estimate is 
established.
$\blacksquare$

\medskip

Now finish proof of lemma~\ref{rho_estim}. If point $\xi_M$ constructed in lemma \ref{lemma:chain_construct} lies in $\mathbb D^{(c)}$ then shrink curve $\Gamma$ by cutting from $\Gamma$ its arc after point $\xi_M$.  Now we assume that $\Gamma \setminus \bigcup\limits_{m=0}^M \Gamma_m\subset\Omega$.

Curve $\Gamma$ joins point on $B_j$ with a  point on $\mathbb D^{(c)}$. Hence, 
$$
\mathcal H^1\left(\Gamma \setminus \bigcup\limits_{m=0}^M \Gamma_m\right)+\sum\limits_{m=0}^M \diam(\Gamma_m) \ge \dist(B_j, \partial\mathbb D).
$$
At the same time, $\Gamma \setminus \bigcup\limits_{m=0}^M \Gamma_m\subset 
\Omega$ by lemma~\ref{lemma:chain_construct}. Also, by lemma 
\ref{capacity_estim_lemma}, 
$$
\mathcal H^1(\Gamma_m\cap \Omega) \ge \frac{\diam(\Gamma_m)}{2e^{2\pi \tilde C_B^2(\Omega)}}.
$$
Accounting disjointness of curves $\Gamma_m$,  $m=0,1,\dots,M,$ we get $\mathcal H^1(\Gamma\cap \Omega) \ge \dfrac{\dist(B_j, \partial\mathbb D)}{2e^{2\pi \tilde C_B^2(\Omega)}}$. Lemma~\ref{rho_estim} is proved.
$\blacksquare$

\medskip

Now we finish proof of theorem~\ref{th:uniform_local_finiteness}. Function $u$ 
was defined above by means of metric  $\rho$. Also, $u$ is constant on any of 
the sets $B_j$, and estimate 
$$
u|_{B_j} \ge \frac{\dist(B_j, \partial\mathbb D)}{2e^{2\pi \tilde C_B^2(\Omega)}},
$$
is true, whereas boundary values of $u$ on $\partial\mathbb D$ are zero. Moreover, $|\nabla u| \le 1$ almost everywhere.
By the second assertion of proposition~\ref{criterii_adm}, we have 
$$
\pi \ge \|u\|^2_{\sob} \ge C_I^{-2}(\Omega) \cdot \sum\limits_{j=1}^\infty (u|_{B_j})^2 \ge  \frac1{4 C_I^2(\Omega) e^{4\pi \tilde C_B^2(\Omega)}} \cdot \sum\limits_{j=1}^\infty \dist(B_j, \partial\mathbb D)^2,
$$
Moreover, weak Bessel property implies that quantities $\diam_H B_j$ are bounded from the above by a constant depending only on  $\tilde C_B(\Omega)$. This means, in particular, that quantities
$$
\frac{\max\{\dist(z, \partial\mathbb D)\colon z \in B_j\}}{\min\{\dist(z, \partial\mathbb D)\colon z \in B_j\}}
$$
are bounded by a constant depending only on $\tilde C_B(\Omega)$. Taking in account that 
$$
\sum\limits_{j=1}^\infty \dist(B_j, \partial\mathbb D)^2 \le 4\pi C_I^2(\Omega) e^{4\pi \tilde C_B^2(\Omega)},
$$
we conclude that disk $\mathcal B(0,1/2)$ can intersect only a finite number of holes $B_j$, the number of which can be estimated from the above only through  $C_I(\Omega)$ and $\tilde C_B(\Omega)$. Theorem is proved. $\blacksquare$

\medskip

\noindent {\bf Remark 1.} We again obtained the convergence of the series 
$\sum\limits_{j=1}^\infty \dist(B_j, \partial\mathbb D)^2$, now in weaker 
assumptions than in theorem~\ref{patch_theorem}. Non-sharpness of the exponent 
$2$ involved in this sum will be discussed in subsection 
\ref{subsection:Blaschke}.

\medskip 

\noindent {\bf Remark 2.} Bessel property (or weak Bessel property) alone or 
interpolation property alone do not imply uniform local finiteness property of $\Omega$.

To construct a domain $\Omega$ with $C_B(\Omega)$ not large but with big $N(\Omega)$, pick large $M\in\mathbb N$. For $m=0,1,\dots, M$ put $B_m:=\bar{\mathcal B}\left(\frac m{2M}, \frac 1{10M}\right)$. These holes are $3/2$-strongly separated, then Bessel constant of  domain $\Omega_M:=\mathbb D\setminus\bigcup_{m=0}^M B_m$ does not exceed some absolute constant (proposition~\ref{predl:strong_sep_Bessel_suff}). At the same time,  $N(\Omega) \ge M$. Acting with conformal shifts of such domains with $M\to\infty$ like in example~\ref{example:inverse} in subsection~\ref{section:examples}, one can construct a domain $\Omega$ for which $C_B(\Omega) <+\infty$ but $N(\Omega) =+\infty$.

Now construct a domain $\Omega$ with $C_I(\Omega)$ not large but large $N(\Omega)$. Again, pick large $M\in\mathbb N$ and also a small $\delta>0$. Further, take an increasing sequence of real numbers $x_1, x_2, \dots, x_M\in(-1,1)$ such that for disks $B_m:=\bar{\mathcal B}\left(x_m, \frac1{2M}\right), \, m=1,2,\dots, M,$ the equalities  $\dist(B_M, \partial\mathbb D) = \delta$,  $\dist(B_{m}, B_{m+1})=\delta$ are held for  $m=1,2,\dots, M-1$. To estimate interpolation constant of domain $\Omega_M := \mathbb D\setminus\bigcup_{m=1}^M B_m$ from the above, we apply proposition \ref{criterii_adm}. Let  $u\in\mathcal Adm(\Omega_m)$ be an admissible function. Then, according to proposition~\ref{predl:capacity_diam_low_estim}, we have inequalities $|(u|_{B_M})| \le  C(M, \delta)\cdot \|u\|_{\sob}$, \, $|(u|_{B_m})-(u|_{B_{m+1}})| \le  C(M, \delta)\cdot \|u\|_{\sob}$ for $m=1,2,\dots, M-1$, where quantity $C(M,\delta)\to 0$ if $\delta\to 0$ with $M$ fixed. Consecutive application of these estimates for $m=M-1, M-2, \dots, 2, 1$ 
lets us to conclude that inequality $\sum\limits_{j=1}^\infty (u|_{B_j})^2\le \|u\|_{\sob}^2$ will be true if for fixed $M$ number $\delta$ is small enough; then $C_I(\Omega_m) \le 1$. To the other hand, if $M$ is large  and $\delta$ is 
small then  we have $N(\Omega_M) \ge M/2$. To construct a  domain $\Omega$ for which $C_I(\Omega) <+\infty$ but $N(\Omega) =+\infty$, one may act by conformal shifts of obtained domains $\Omega_M$ for $M\to\infty$, as in example~\ref{example:inverse} in subsection~\ref{section:examples}.

\medskip 

\noindent {\bf Remark 3.} In the construction of function $u$ in the proof of theorem~\ref{th:uniform_local_finiteness} we were solving the following extremal problem: \emph{to maximize $\sum_{j=1}^\infty (u|_{B_j})^2$ if $u\in W^{1,2}(\mathbb D)$, $u|_{\partial\mathbb D}=0$, $|\nabla u|\le \varphi$ almost everywhere}. This was done for $\varphi = \mathds 1_{\Omega}$. One may also solve such a problem in the most generality. First, under very general assumptions there exists a function $u$ giving maximum to $u(z)$ simultaneously in (almost) all points $z\in\mathbb D$. Second,  function $\varphi\colon\mathbb D\to [0,+\infty)$ can be more or less arbitrary. In other words, one can, under pointwise upper estimate on $|\nabla u|$, maximize any norm of $u$ monotonically depending on its values. Moreover, such a function $u$ can be found as explicitly as we found one in the proof of theorem~\ref{th:uniform_local_finiteness}.

Metric $\rho$ then may be defined in an analogous way: if $z_1, z_2 \in \clos{\mathbb D}$ then 
\begin{equation}
\label{eq:rho_general}
\rho(z_1, z_2) := \inf \int_{\Gamma}\varphi\,d\mathcal H^1,
\end{equation}
where $\inf$ is taken over all curves $\Gamma$ joining $z_1$ with $z_2$. Further, one may put $u(z):=\dist_\rho(z,\partial\mathbb D), \, z\in\mathbb D$. Nevertheless,   problems with measurability of function under  integral in 
(\ref{eq:rho_general}) may arise now. Moreover, function $u$ may collapse. For example, for $x+iy\in\mathbb D$ put $\varphi(x+iy):=0$ if $x$ or  $y$ is rational and  $\varphi(x+iy):=1$ otherwise; then for metric $\rho$ defined by 
(\ref{eq:rho_general}) we have $\rho(z_1, z_2)=0$ for all $z_1, z_2 \in \clos{\mathbb D}$, and $u\equiv 0$ in $\clos{\mathbb D}$.

The solution of such difficulties is to replace $\inf$ in~(\ref{eq:rho_general}) by an \emph{essential infimum}, that is, by an infimum up to \emph{negligible}  family of curves. Such a  smallness should be understood in the sense of module (or extremal length) of the exceptional family, see~\cite{Fuglede57},~\cite{Fuglede60}. Fortunately, we did not need such a technique: in our case function $\varphi$ is not very non-smooth since it coincides with an indicator of a locally-smooth set $\Omega$. And difficulties related to the possibility of non-smoothness of curves $\Gamma$ in~(\ref{eq:rho_general}) were overcome by working with finite-verticed polygonal chains.

\section{Criteria under uniform local finiteness condition}

Our goal is to give a metric criterion of complete interpolation property. Questions on description of domains possessing only Bessel or only interpolation properties seem to be unreachable. At the same time, these questions can be answered if we, in addition, impose uniform local finiteness property on the family of holes $B_j$. We already proved (in theorem  \ref{th:uniform_local_finiteness})  that this is the case when $\Omega$ has complete interpolation property.

\subsection{Bessel property criterion}

Recall that notion of weak separatedness of holes was introduced in definition~\ref{def:separatedness}. Propositions~\ref{predl:capacity_diam_low_estim} and~\ref{necessary_cap} imply that weak (or strong) Bessel property implies weak separatedness of holes. In this subsection we will prove the opposite under condition of uniform local finiteness.

\begin{theorem}
	\label{th:bessel_sufficient}
	Let $\Omega$ be regular domain, $B_j$ be holes in $\Omega$. Suppose that $\sup\limits_{j\in\mathbb N} \diam_H(B_j) < +\infty$, $N(\Omega) < +\infty$ \emph{(}that is, any disk of the form $\mathcal B_H(z, 1)$ intersects no more than $N(\Omega)$ of holes $B_j$\emph{)} and that holes $B_j$ are $\eps$-weakly separated for some $\eps>0$. Then $\Omega$ possesses Bessel property and Bessel constant $C_B(\Omega)$ can be estimated from the above only through $\sup\limits_{j\in\mathbb N} \diam_H(B_j), N(\Omega)$ and $\eps$.
\end{theorem}

First we prove the following
\begin{predl} 
\label{predl:weak_bessel_suff}
Hypothesis of theorem~\ref{th:bessel_sufficient} provides weak Bessel property of $\Omega$. Moreover, constant $\tilde C_B(\Omega)$ can be estimated from the above only in terms of constants from this theorem.
\end{predl}

\noindent{\bf Proof.}
Indeed, according to proposition \ref{necessary_cap}, we have to get upper-estimate for capacity $\Cap2(B_j, \Omega^{(c)}\setminus B_j)$,  $j=1,2,\dots$. Recall that condition $\sup\limits_{j\in\mathbb N} \diam_H(B_j) < +\infty$ means that quantities $\dfrac{\dist(B_j, \partial \mathbb D)}{\diam (B_j)}$ are separated from zero by some constant $c>0$.  Let $\{B_{j'}\}_{j'\in J}$ be the set of holes $B_{j'}\,(j'\neq j)$ for which $\dist(B_j, B_{j'}) \le (c/2)\cdot \diam(B_j)$.  Let~$U_j$ be  $\left({c\cdot \diam(B_j)}/2\right)\mbox{-neighbourhood}$ of $B_j$. Quantity  $\diam_H U_j$ is finite and can be estimated from the above via $c$, hence, due 
to uniform local finiteness, the set of indices~$J$ is finite and number of its elements can be estimated through the constants from the hypothesis. Capacity $\Cap2\left(B_j,\mathbb D^{(c)}\cup \bigcup\limits_{j'\notin J,\, j'\neq j} 
B_{j'}\right)$ is no more than $\Cap2(B_j,  U_j^{(c)})$. As it was proved in proposition~\ref{predl:capacity_upper}, the last capacity is bounded from the above by a value depending only on  $c$. Moreover, for any $j'\in J$ the  
capacity $\Cap2(B_j, B_{j'})$ is bounded from the above by a value depending only on the constant~$\eps$ of weak separatedness of holes (corollary \ref{sled:weak_separ_cap}). The estimate on $\card J$ and semiadditivity  of 
capacity $\Cap2(\cdot, \cdot)$ by the second argument now imply the desired upper estimate for $\Cap2(B_j, \Omega^{(c)}\setminus B_j)$.~$\blacksquare$

\medskip

\noindent {\bf Proof of theorem~\ref{th:bessel_sufficient}.} By the proposition~\ref{predl:weak_bessel_suff} just proved, the quantity $\tilde C_B(\Omega)$ is finite and can be estimated through constants from the hypothesis of theorem. But, by theorem~\ref{th:weak_strong_bessel}, $C_B(\Omega) \le \sqrt 2 \cdot C_B(\Omega)$, what concludes the proof.
$\blacksquare$

\medskip

Theorem~\ref{th:bessel_sufficient} can be also proved without application of theorem~\ref{th:weak_strong_bessel} and in terms of only metric characteristics of domain~$\Omega$. It is interesting that a \emph{discrete structure of partial order} at the set of holes arises naturally in such a proof.

For $a>0$ denote by $U_a(B_j)$ the closed $a$-neighbourhood of set~$B_j$. Put $\Lambda_j:=\diam B_j$.

The next lemma, in fact, gives a description of domains satisfying hypothesis of theorem~\ref{th:bessel_sufficient}.

\begin{wrapfigure}[17]{r}{0.55\textwidth}
	\centering
{\includegraphics[scale=1.1]
	{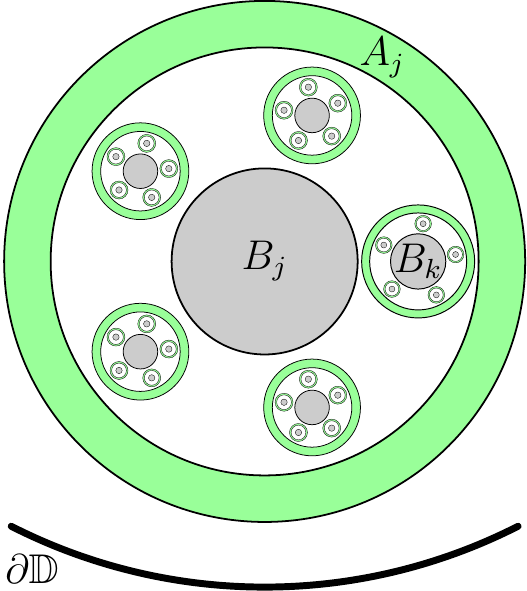}}
\caption{Partial order on the set of holes}
	\label{fig:bessel_order}
\end{wrapfigure}

\begin{lemma}[see fig. \ref{fig:bessel_order}] 
	\label{lemma:annular_construction_conditions}
	In the settings of theorem~\ref{th:bessel_sufficient} it is possible to define a strict partial order relation~$\succ$ on the set of holes $B_j$ and also to associate a set $A_j\subset\mathbb D$ to each hole $B_j$, $j=1,2,\dots$, such that:
\begin{enumerate}
 \item 
 For each $j=1,2,\dots$ set $A_j$ is of the form $U_{s_j}(B_j)\setminus U_{t_j}(B_j)$ for some $t_j, s_j \in (0,+\infty) \, (s_j>t_j)$. Also,  $s_j - t_j \ge c_2 \Lambda_j$, $s_j \le c_1 \Lambda_j$ where constants $c_1, c_2\in (0,+\infty)$ do not depend on $j$. 
 Moreover, $A_j \subset \Omega$.

 \item
 Overlapness multiplicity of sets $A_j$ does not exceed some constant $C_1<+\infty$.
 
 \item 
 $B_{j'} \prec B_{j} \, (j'\neq j)$ if and only if $U_{s_{j'}}(B_{j'}) \subset  U_{s_j}(B_j)$, and also if and only if $B_{j'}\subset U_{t_{j}}(B_{j})$.
 
 \item  For a fixed $j_0\in\mathbb N$ the number of indices $j$ for which $B_j \prec B_{j_0}$ does not exceed some constant $C$. In particular, lengths of chains in order $\succ$ are bounded from the above.  If $B_{j_1}, B_{j_2} \succ B_j$, then either $B_{j_1}\succ B_{j_2}$, or $B_{j_2}\succ B_{j_1}$.
  
\end{enumerate}

Constants $c_1, c_2, C, C_1$ depend only on constants from theorem~\ref{th:bessel_sufficient}.
\end{lemma} 

The proof is given in the Appendix. From this proof it is seen that we may force~$c_1$ to be arbitrarily small. Thence, relation $B_k \prec B_j$ for some $j,k=1,2,\dots$ means, roughly speaking, that $\diam B_k \ll \diam B_j$ and $\dist(B_j, B_k) \ll \diam B_j$.

\medskip

\noindent {\bf Remark.} The fourth assertion of lemma \ref{lemma:annular_construction_conditions} allows to endow the set of holes by the structure of \emph{oriented forest}. 

Namely, denote by $\mathcal M$ the set of indices of maximal elements in order $\succ$. Let  $j\notin \mathcal M$. According to the fourth assertion of lemma~\ref{lemma:annular_construction_conditions}, set $\{B_{j'}\colon B_{j'}\succ B_j\}$ is finite and is a chain in order $\succ$. Therefore, there exists a unique minimal element $B_{j''}$ in this set; denote $j''$ by $\mathcal P(j)$. Now we defined the mapping of "the nearest ancestor"{} $\mathcal P \colon \mathbb N \setminus \mathcal M \to \mathbb N$.

Let us define an oriented graph $\mathscr G$. Its vertices will be all the holes $B_j, \, j\in\mathbb N$. For each $j\in\mathbb N\setminus \mathcal M$, let us draw an edge in  graph $\mathscr G$ from $B_{\mathcal P(j)}$ to $B_j$. Let $\tilde{\mathscr G}$ be an non-oriented graph obtained by forgetting of directions of edges in $\mathscr G$. 
The fourth assertion of lemma~\ref{lemma:annular_construction_conditions} allows to conclude that graph $\tilde{\mathscr G}$ is a forest, that is, a countable disjoint union of trees. Also, numbers of vertices in connected components of $\tilde{\mathscr G}$ are bounded from the above. In any of such components, one vertex has zero incoming degree in graph $\mathscr G$, whereas all the other vertices in this component have incoming degree $1$.

\medskip

Suppose that all holes $B_j$ are disks. Then any set $A_j$, $j=1,2,\dots$, constructed in lemma \ref{lemma:annular_construction_conditions} is an annulus; $A_j$ winds exactly around those holes $B_{j'}$ for which $B_{j'}\preceq B_{j}$. By the first assertion of lemma \ref{lemma:annular_construction_conditions} annulus  $A_j$ is wide enough, and if $\omega\in\lo$ then 
$$
\int\limits_{A_j}|\omega|^2\dl \ge c\cdot\left(\sum\limits_{j'\colon B_{j'}\preceq B_j}\Per_{j'}\omega\right)^2,
$$
with some $c>0$ not depending on $j$ and  $\omega$. Consecutive application of this inequality starting from the minimal holes in order $\succ$ and up to maximal ones and accounting boundedness of lengths of chains in order $\succ$ and bounded overlapness multiplicity of sets  $A_j$, we can estimate $\|\Per\omega\|_{\ell^2}$ from the above through $\|\omega\|_{\lo}$.

A more strict derivation of theorem~\ref{th:bessel_sufficient} from lemma \ref{lemma:annular_construction_conditions} without assumption that all the holes are disks is given in the Appendix.

\subsection{Interpolation criterion}
\label{subsec:interp_criterion}

Now our goal is to establish necessary and sufficient condition for interpolation if uniform local finiteness of family of holes $B_j$ is held. For this aim we will need a graph.

Pick a number $S\in(0,+\infty)$. Define a non-oriented graph $G(\Omega, S)$ 
as follows. Vertices of graph $G(\Omega, S)$ will be all the connected components of 
$\Omega^{(c)}$, that is, all the holes~$B_j$,  $j=1,2,\dots,$ and also set 
$\mathbb D^{(c)}$. Join two holes $B_{j_1}$ and $B_{j_2}$ by an edge in the 
graph $G(\Omega, S)$ if $\dist(B_{j_1}, B_{j_2}) \le S\cdot \min\{\diam B_{j_1}, \diam 
B_{j_2}\}$. Further, join $B_j$ and~$\mathbb D^{(c)}$ by an edge if  
$\dist(B_j, \mathbb D^{(c)}) \le S\cdot \diam B_j$. If $S$ increases then the 
set of edges in this graph does.

We define \emph{distance} between two vertices in (non-oriented) graph $G$ as the number of \emph{edges} of the shortest path joining these vertices (thus, the  distance between two adjacent vertices is $1$); if there is no such a path, then the distance is set to be infinity. Denote by $\dist_G$ the distance in the graph. 

\emph{Capacity connectedness} property which we mentioned in the introduction is connectedness of graph $G(\Omega, S)$ for some $S<+\infty$ together with finiteness of its diameter in metric $\dist_{G(\Omega, S)}$; the latter condition will be often understood as estimate $\sup\limits_{j\in\mathbb N} \dist_{G(\Omega, S)} (B_j, \mathbb D^{(c)}) <+\infty$. 

\medskip

\noindent {\bf Remark.} Let $\Omega$ be a \emph{regular} domain, recall that in this case any set $B_j$, $j=1,2,\dots$, is a closure of a domain with smooth boundary. We may  define another graph $g(\Omega, s)$ where $s>0$. Its vertices will be, like in the case of $G(\Omega, S)$, all the connected components of $\Omega^{(c)}$. Join two such vertices $E_1$ and $E_2$ by an edge in the graph $g(\Omega, s)$ if $\Cap2(E_1, E_2) \ge s$. According to  proposition~\ref{predl:capacity_diam_low_estim} if graph $G(\Omega,S)$ is connected or its diameter admits an upper estimate, then, for some $s>0$ depending only on $S$, graph $g(\Omega,s)$ will be connected (or, respectively, its diameter will admit an upper-estimate). By corollary~\ref{sled:weak_separ_cap}, the opposite is true as well:  connectedness  (or boundedness of diameter) of graph $g(\Omega, s)$ implies connectedness (respectively, boundedness of diameter) of graph $G(\Omega, S)$ for some $S<+\infty$ depending only on $s$.

In this subsection we will mainly deal with graph $g(\Omega, s)$ using its conformal invariance. Graph $G(\Omega, S)$ will be appropriate for criteria in the case when holes $B_j$ have non-smooth boundaries (subsection \ref{subsection:nonsmooth}).

\begin{theorem}
\label{th:interp_sufficient}
	Let $\Omega$ be a regular domain and suppose that $N(\Omega)<+\infty$ \emph{(}that is, that family of holes in $\Omega$ is uniformly locally finite\emph{)}. 
	If $\Omega$ possesses interpolation property then there exist $S<+\infty$ and $M \in \mathbb N$ such that graph $G(\Omega, S)$ is connected and for any hole $B_j$ distance from $B_j$ to $\mathbb D^{(c)}$ in the graph $G(\Omega, S)$ is no more than $M$.
	
	Numbers $S$ and $M$ can be estimated from the above only through $C_I(\Omega)$ and $N(\Omega)$.
\end{theorem}

\noindent {\bf Remark.} As it will be seen from the argument below, we may, for appropriate $S$, take $M=N(\Omega)+1$.
	
\medskip	
	
\noindent {\bf Proof of theorem~\ref{th:interp_sufficient}.} The idea of the argument in the below can be briefly sketched  as follows. Suppose that  $r<1/4$ and let $\mathcal C$ be some connected component of graph $G(\Omega, S)$ such that all vertices from $\mathcal C$ lie in $\mathcal B(0,r)$. Assume also that annulus $\mathcal B(0,1/2)\setminus\mathcal B(0,r)$ lies in $\Omega$. Suppose that either $r$ is small, or the number of vertices in $\mathcal C$ is large. Function $u$ for which $u|_{\mathcal B(0,r)}=1$, $u|_{\mathbb C\setminus\mathcal B(0,1/2)}=0$, $\Delta u=0$ in $\mathcal B(0,1/2)\setminus\bar{\mathcal B}(0,r)$ is admissible for $\Omega$. Applying the second assertion of proposition~\ref{criterii_adm} to $u$, we conclude that the constant $C_I(\Omega)$ is big, what contradicts the hypothesis of the theorem. 

\medskip

Now we pass to the formal argument. According to remark before the theorem, it is sufficient to check the connectedness and boundedness of diameter of graph $g(\Omega,s)$ for some $s>0$ depending only on $N(\Omega)$ and $C_I(\Omega)$.

Let $N=N(\Omega)$. Put $\tilde s=C_I^{-2}(\Omega)/2N$. Choose $F=F(\tilde s)$ such that if $\Cap2(B_j, B_{j'}) \ge \tilde s$ then $\dist(B_j, B_{j'}) \le F\cdot \min\{\diam B_j, \diam B_{j'}\}$ (this can be done by corollary~\ref{sled:weak_separ_cap}); we may assume that $F \ge 1$. Pick a number $\delta>0$ such that  $3\delta N F < \tanh 1$ and also such that if $B_j\subset \mathcal B(0, 3 N F\delta)$ then $\Cap2\left(B_j,  (\mathcal B(0, \tanh 1))^{(c)}\right) \le C_I^{-2}(\Omega)/2$. Finally, pick a positive number $s \le \tilde s$ such that if $\diam B_j \ge \delta$ then $\Cap2(B_j, \mathbb D^{(c)}) \ge s$. By the construction, $s$ depends only on~$N$ and $C_I(\Omega)$. Let us show that graph $g(\Omega, s)$ is connected and its diameter does not exceed $2N+2$.

Take any hole $B_{j_1}$. We are going to prove that $\dist_{g(\Omega, s)}(B_{j_1},  \mathbb D^{(c)}) \le N+1$. Uniform local finiteness property, interpolation property and $g(\Omega, s)$ are invariant under M{\"o}bius transforms, hence we may assume that $0\in B_{j_1}$.  Let ${\mathcal C}$ be  connected component  of vertex $B_{j_1}$ in the graph $g(\Omega,  s)$ and  $\tilde {\mathcal C}$ be  connected component  of vertex $B_{j_1}$ in the graph $g(\Omega, \tilde s)$ (this is a subgraph in $g(\Omega, s)$). We have $\tilde{\mathcal C} \subset \mathcal C$. Let us prove the following: \emph{one of the vertices $B_j \in\tilde{\mathcal C}$ is adjacent to $\mathbb D^{(c)}$ in $g(\Omega, s)$}. Assume the contrary.

Let $B_j \in \tilde {\mathcal C}$.  Then $\diam B_{j}\le \delta$ (otherwise $\Cap2(B_j, \mathbb D^{(c)}) \ge s$ by the choice of $s$ and then vertices $B_j$ and $\mathbb D^{(c)}$ are adjacent in the graph $g(\Omega, s)$ what contradicts our assumption). Now we will work with graph $g(\Omega, \tilde s)$.

Let vertices $B_{j}, B_{j'}\in\tilde{\mathcal C}$ be adjacent  in the graph $g(\Omega, \tilde s)$, then $\Cap2(B_j, B_{j'}) \ge \tilde s$, hence 
$$
\dist (B_{j}, B_{j'}) \le F\cdot \min\{\diam B_{j},\diam B_{j'}\} \le \delta F
$$
(the first inequality is true by the choice of $F$, while the second -- by the 
diameters estimate proved above). If $\tilde{\mathcal C}$ contains a simple 
path of length $N+1$ consisting of edges of graph $g(\Omega, \tilde s)$   and 
starting in  $B_{j_1}$ then all  $B_j$'s from this path lie in $\mathcal 
B(0,(N+1)\delta+N\delta F)$ since diameters of sets $B_j$ from this path are no 
more than~$\delta$, while distances between adjacent sets do not exceed $\delta 
F$ (recall that $0\in B_{j_1}$). By the choice of $\delta$ we have 
$(N+1)\delta+N\delta F\le \tanh 1$. But by the definition of $N$ disk $\mathcal 
B(0, \tanh 1)$ can not intersect more than $N$ of holes $B_j$. Hence 
$\tilde{\mathcal C}$ can not contain simple paths starting in $B_{j_1}$, going 
by edges of $g(\Omega, \tilde s)$ and of length more than $N$. It means that 
any vertex $B_j \in \tilde{\mathcal C}$ can be joined in $g(\Omega, \tilde s)$ 
with $B_{j_1}$ by a path of length no more than $N$. Then, analogously to the 
distances estimates worked out above,  $B_j\subset \mathcal 
B(0,N\delta+(N-1)\delta F)\subset\mathcal B(0, 3\delta FN)\subset\mathcal B(0, 
\tanh 1)$. This implies that number of vertices in $\tilde{\mathcal C}$ is 
finite and no more than $N$. So, let $\tilde{\mathcal C}= 
\{B_{j_1},B_{j_2},\dots, B_{j_K}\}, \, K\le N$.

Let $J=\{j\in \mathbb N\colon j\neq j_k,\, k =1, 2, \dots, K\}$. Interpolation property of domain $\Omega$ implies that 
$$
\Cap2\left(B_{j_1}\cup B_{j_2}\cup \dots\cup B_{j_K}, \mathbb D^{(c)}\cup\bigcup\limits_{j\in J} B_j\right) \ge K\cdot C_I^{-2}(\Omega).
$$
Indeed, it is enough to apply the second assertion of proposition 
\ref{criterii_adm} to functions in $\inf$ for the definition of capacity in the 
left-hand side of the last inequality (see definition~\ref{def:capacity}). 
Semiadditivity of capacity (proposition~\ref{predl:capacity_simple}) allows us 
to conclude that there exists $k=1,2,\dots, K$ for which 
$$
\Cap2\left(B_{j_k}, \mathbb D^{(c)}\cup\bigcup\limits_{j\in J} B_j\right) \ge C_I^{-2}(\Omega).
$$
It has already been proved that $B_{j_k}\subset \mathcal B(0, 3 \delta F N)$. If $\{B_j\}_{j\in J_1}$ is the set of holes lying outside of $\mathcal B(0, \tanh 1)$ then 
$$
\Cap2\left(B_{j_k}, \mathbb D^{(c)}\cup\bigcup\limits_{j\in J_1} B_j\right) \le \frac{C_I^{-2}(\Omega)}{2N}
$$
by the choice of $\delta$. Thus, again by semiadditivity of capacity, 
$$
\Cap2\left(B_{j_k}, \bigcup\limits_{j\notin J_1\cup J} B_j\right) \ge \frac{C_I^{-2}(\Omega)}{2}.
$$
In the capacity from the last inequality, the second plate is formed by holes $B_j\notin\tilde{\mathcal C}$ intersecting with $\mathcal B(0, \tanh 1)$. The number of such holes is no more than $N$. Hence, for one of them, say, for $B_\alpha$, the following inequality should be true:
$$
\Cap2(B_{j_k}, B_\alpha) \ge \frac{C_I^{-2}(\Omega)}{2N}=\tilde s. 
$$
But, by the definition of graph $g(\Omega, \tilde s)$, in this case vertices $B_{j_k}$ and $B_\alpha$ are adjacent in the graph $g(\Omega, \tilde s)$, and hence $B_\alpha\in\tilde{\mathcal C}$. We get a contradiction. So, we proved that  graph $g(\Omega, \tilde s)$ contains  a path starting in $B_{j_1}$ and such that its last vertex is adjacent to $\mathbb D^{(c)}$ in graph $g(\Omega, s)$ (in particular, graph $g(\Omega, s)$ is connected). 

Suppose that $\dist_{g(\Omega, s)}(B_{j_1}, \mathbb D^{(c)}) > N+1$. Consider the shortest simple path in the graph $g(\Omega, \tilde s)$ from $B_{j_1}$ to some $B_\alpha$ adjacent to $\mathbb D^{(c)}$ in the graph $g(\Omega, s)$. By the assumption, the length of this path is not less than $N+2$. Consider its first $N+1$ vertices, denote them by $B_{j_1}, B_{j_2}, \dots, B_{j_{N+1}}$. If for some of them $\diam B_{j_k} \ge \delta\, (k=1,2,\dots,N+1)$ then in the graph $g(\Omega, s)$ vertex $B_{j_k}$ is adjacent to $\mathbb D^{(c)}$ by the choice of $s$. Then distance in  $g(\Omega, s)$ from $B_{j_1}$ to $\mathbb D^{(c)}$ is no more than $N+1$, but that is the desired. If for all holes  of path $(B_{j_1}, B_{j_2}, \dots, B_{j_{N+1}})$ inequality $\diam B_j < \delta$ is true then, by the choice of $F$, Euclidean distance between adjacent vertices of this path is no more than~$F\delta$ (since $B_{j_k}$ is adjacent to  $B_{j_{k+1}}$ in the graph $g(\Omega, \tilde s), \, k=1,2,\dots,N$). In this case $B_{j_k}\subset \mathcal B(0,3 N F\delta)\subset \mathcal B(0, \tanh 1)$ for all $k=1,2,\dots, N+1$. But disk $\mathcal B(0, \tanh 1)$ cannot intersect more than $N$ holes. We got a contradiction again.

So, under our choice of $s$ distance from any vertex of graph $g(\Omega, s)$ to vertex $\mathbb D^{(c)}$ does not exceed $N+1$ and diameter of this graph is no more than $2N+2$, the desired. Theorem is proved. $\blacksquare$

\medskip

Now we are going to prove the opposite: if diameter of graph $G(\Omega, S)$ is finite then~$\Omega$ possesses interpolation property. Our argument may be simplified if we add estimate $\sup\limits_{j\in\mathbb N}\diam_H(B_j)<+\infty$ to the hypothesis of the following

\begin{theorem}
\label{th:interp_tree}
	Let the family of holes in a regular domain $\Omega$ be uniformly locally finite. Suppose that, for some $S<+\infty$, graph $G(\Omega, S)$ is connected and that  $\dist_{G(\Omega, S)}(B_j, \mathbb D^{(c)}) \le M$ for any $B_j$ and for some $M \in \mathbb N$ not depending on $j$. Then domain $\Omega$ has interpolation property, moreover, constant $C_I(\Omega)$ can be estimated from the above in terms of only $S, M$ and $N(\Omega)$.
\end{theorem}

\begin{figure} 
	\centering
{\includegraphics[scale=0.85]
	{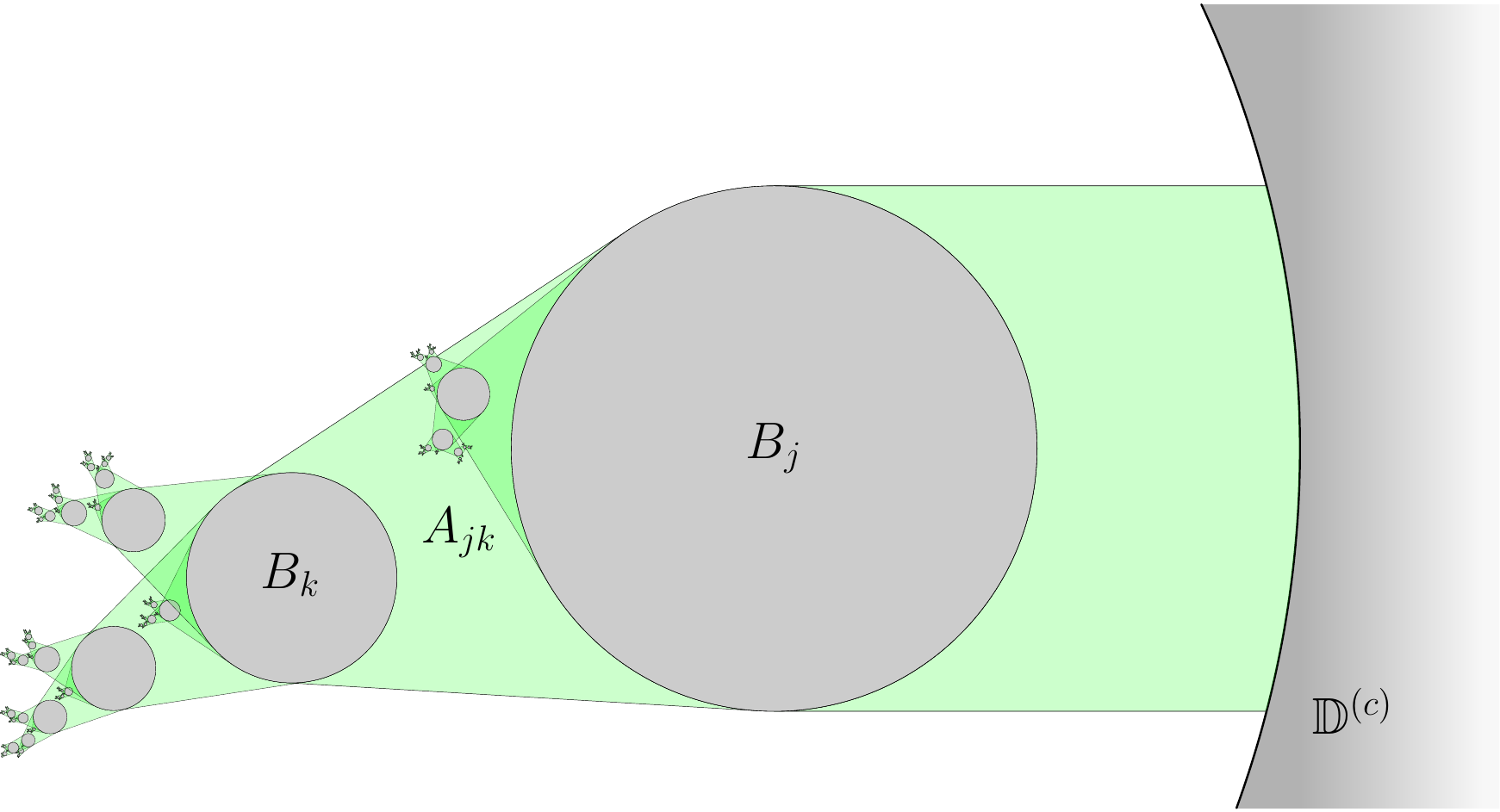}}
\caption{Spanning tree of graph $G(\Omega, s)$}
	\label{fig:tree}
\end{figure}

\noindent {\bf Proof.} 
 First let us sketch the proof. According to proposition~\ref{criterii_adm},  we have to prove that for any function $u\in\mao$ (that is, recall, a function from $\sob$ which is constant on any hole $B_j$) the following estimate is true:
$$
\int_{\mathbb D}^{} |\nabla u|^2\dl \ge C_I^{-2}(\Omega)\cdot \sum_{j=1}^{\infty}(u|_{B_j})^2,
$$
If $B_j$ and $B_k$ are joined by an edge in $G(\Omega, S)$,  then we will draw a "road"{} $A_{jk}$ between $B_j$ and $B_k$ (see fig.~\ref{fig:tree}) which is wide enough such that 
\begin{equation}
\label{eq:interp_explain}
\left|(u|_{B_j})-(u|_{B_k})\right|^2\le C\cdot\int_{A_{jk}}|\nabla u|^2\dl.
\end{equation}
(The constant  $C<+\infty$ depends only on the constants from the hypothesis of theorem but not on $j,k$ and $u$.) \emph{Almost} all $A_{jk}$ have bounded overlapness multiplicity. Moreover, the diameter of $G(\Omega, S)$ is finite. Starting from vertex $\mathbb D^{(c)}$, passing the graph in a breadth-first order and applying estimates~(\ref{eq:interp_explain}) consecutively we get the desired inequality.

If we argue in such a manner we will need, in particular, "roads"{} joining some holes~$B_j$ with $\mathbb D^{(c)}$. Unfortunately, overlapness multiplicity of such "roads"{} may turn to be unbounded. Thus we have to estimate  $u|_{B_j}$ for such holes $B_j$ separately. If a hole $B_j$ is adjacent to $\mathbb D^{(c)}$ in the graph $G(\Omega, S)$ then $\diam_H(B_j) \ge c$ where $c>0$ and depends only on $S$. But in case of such holes we may argue as in theorems~\ref{th:disks_interp_sufficient} and~\ref{th:arbitrary_interp_sufficiency_separated}. 

\medskip

Now we pass to the formal argument. During this proof constants $C,C_1, C_2, \dots \in(0, +\infty)$ depend only on constants $S, M$ and $N(\Omega)$ from hypothesis. 

According to remark before theorem~\ref{th:interp_sufficient}, there exists $s >0$ depending only on $S$ such that graph $g(\Omega, s)$ is connected and $\dist_{g(\Omega, s)}(B_j, \mathbb D^{(c)}) \le M$ for any $B_j$. We will be more convenient to work with this graph using conformal invariance of capacity.

Fix a function $u\in\mao$. Put $a_j := u|_{B_j},\, j=1,2,\dots$. We have to get estimate 
$$
\sum\limits_{j=1}^\infty a_j^2 \le C\cdot \|u\|_\sob^2
$$
with $C<+\infty$ depending only on constants from the hypothesis but not on $u$.

By corollary~\ref{sled:weak_separ_cap}, there exists a number $F=F(s)$ such that: if $B_j$ and $B_k$ are adjacent  in the graph $g(\Omega, s)$, that is $\Cap2(B_j, B_k) \ge s$, then $\dist(B_j, B_k) \le F\cdot \min\{\diam B_j , \diam B_k\}$. Pick $\delta\in(0,1/3]$ small enough. We will require the following. First, if $B_j$ is adjacent in $g(\Omega, s)$ with $\mathbb D^{(c)}$ then  $\diam_H(B_j) \ge \delta$. Second, if $B_j \subset \mathcal B(0,\delta)$ and $B_k\cap \mathcal B(0, 1/3)=\varnothing$ then $\Cap2(B_j , B_k) < s$ (and hence  $B_j$ and $B_k$ are not adjacent in the graph $g(\Omega, s)$). Third, $\delta (F+2) <1/3$. 

Denote by $J^\delta_1$ the set of indices $j=1,2,\dots$ for which $\diam_H(B_j) 
\ge \delta$. Next, denote by $J_2^\delta$ the set of indices $j\notin 
J^\delta_1$ for which $B_j$ is adjacent in $g(\Omega, s)$ with one of the 
vertices  $B_k$,  $k \in J^\delta_1$. (We could add in $g(\Omega, s)$ all the 
edges from vertices $B_j, \, j \in J^\delta_1,$ to $\mathbb D^{(c)}$ and 
arrange all holes by their  distances from $\mathbb D^{(c)}$ in the obtained 
graph. We, though, will not use all the advantages of such a trick.) 

Recall that if $Z$ is finite set then symbol $\card Z$ denotes the number of its elements.

\begin{lemma}
	\label{lemma:interp_suff_j2}
	Let  $\{E_j\}_{j\in J}$ be finite or countable family of compact connected disjoint sets in disk $\mathbb D$. We are assuming that any $E_j, \, j\in J,$  is a closure of a domain with {$C^\infty$-smooth} boundary. Suppose that $\diam_H (E_j) \ge \delta$ for all $j\in J$ where $\delta\in(0,1/3]$.
	
	Suppose also that we are given by a finite or countable set $X\subset \bigcup\limits_{j\in J} E_j$ for which $n_j:=\card(X\cap B_j) <+\infty$ for any $j\in J$. Assume also that set $X$ is \emph{uniformly locally finite} in the following sense: any hyperbolic disk of radius $1$ contains no more than $\nu$ points of set $X$ where $\nu\in\mathbb N$ is constant.
	
	Suppose that a function  $u\in \sob$ is constant almost everywhere on any set~$E_j$, $j\in J,$ and equals to some number $a_j \in\R$ on this set. Then
	$$
	\sum\limits_{j\in J} n_j a_j^2\le \tilde C_1\cdot  \|u\|_{\sob}^2,
	$$
	where $\tilde C_1$ depends only on $\nu$ and $\delta$.
\end{lemma}

\noindent {\bf Proof.}
Take any $z\in X$. Denote by $E(z)$ the set $E_j$ containing $z$. Put $K(z):=\mathcal B_H(z, \delta/2)$. The hypothesis of lemma implies that disks $K(z), \,z\in X,$ have overlapness multiplicity not exceeding $\nu$ (since $\delta \le 1/3$). Let us bring into consideration the \emph{hyperbolic area} 
$$
dA = \frac{4\dl}{(1-|z|^2)^2} \asymp \frac{\dl}{(1-|z|)^2}.
$$
This area is invariant under M{\"o}bius transforms. Let us prove that there exists a constant $\tilde C_2>0$ depending only on $\delta$ for which 
\begin{equation}
\label{estim:interp_suff_lemma1}
(u|_{E(z)})^2 \le \tilde C_2\cdot \left(\int\limits_{K(z)} u^2 \, dA +
\int\limits_{K(z)} |\nabla u|^2 \dl\right).
\end{equation}
The expressions at both sides of the last estimate are conformally invariant. 
Thus we may, by an application of an appropriate M{\"o}bius transform, assume 
that $z=0\in E(0)$. Then $K(z)=K(0)=\mathcal B(0, \tanh(\delta/2))$. 

Let $u|_{E(z)}\equiv a$. For $t\in [0,\tanh(\delta/2)]$ consider circles 
$\beta_t=\partial\mathcal{B}(0, t)$. For any such~$t$ we have $\beta_t\cap 
E(0)\neq \varnothing$ since $\diam_H E(0) \ge \delta$. Moreover, we may assume 
that for almost all $t\in [0,\tanh(\delta/2)]$ the function $u$ is absolutely 
continuous on $\beta_t$ and takes value $a$ on this circle.  Arguing as in the 
proof of theorem~\ref{th:arbitrary_interp_sufficiency_separated}, we obtain the 
estimate  
\begin{equation}
\label{eq:grad_mean}
2\pi t a^2 \le 2\int\limits_{\beta_t} u^2\,d\mathcal H^1+8\pi^2 t^2 \int\limits_{\beta_t}|\nabla u|^2\,d\mathcal H^1.
\end{equation}
If $\zeta=te^{i\theta}\in \beta_t, \, \theta \in [0,2\pi]$, then $(d\mathcal H^1|_{\beta_t})(\zeta) = t\, d\theta, \,\, dA(\zeta) = \dfrac{4t\,dt\,d\theta}{(1-t^2)^2}\asymp \dl(\zeta)$ on~$K(0)$ (the constant of comparability is absolute: recall that we assume that $\delta \le 1/3$). Taking this in account and integrating~(\ref{eq:grad_mean}) over  $t\in[0, \tanh(\delta/2)]$ we get~(\ref{estim:interp_suff_lemma1}) with some constant $\tilde C_2$ depending only on $\delta$.

Now sum up~(\ref{estim:interp_suff_lemma1}) over all $z\in X$ to obtain
$$
\sum\limits_{j\in J} n_j a_j^2 \le \nu \tilde C_2\cdot\left(\int\limits_{\mathbb D} u^2 \, dA +
\int\limits_{\mathbb D} |\nabla u|^2 \dl\right),
$$
since disks $K(z), \, z\in X,$ have overlapness multiplicity not exceeding $\nu$; in the left-hand side we have exactly this sum because $\card(E_j\cap X)=n_j$ for all $j\in J$. The assertion of  theorem~\ref{th:Hardy_neq} can be rewritten in the following way:
$$
\int\limits_{\mathbb D} u^2 \, dA \le \const\cdot \|u\|^2_{\sob}, ~~ u\in \sob,
$$
with absolute constant in the right-hand side. This and and previous inequalities imply the assertion of lemma. $\blacksquare$

\medskip

Now we proceed proving theorem~\ref{th:interp_sufficient}. Let us estimate  
$\sum\limits_{j\in J_1^\delta} a_j^2$. In each set  $B_j$,  $j\in 
J_1^\delta,$ pick a point $\zeta_j$,  put $X:=\{\zeta_j\}_{j\in J_1^\delta}$. 
Set $X$ is uniformly locally finite since a disk of the form $\mathcal B_H(z, 
1)$ can not intersect more than $N(\Omega)$ of holes $B_j$. By the definition 
of the set $J_1^\delta$, if $j \in J_1^\delta$ then $\diam_H(B_j) \ge \delta$. 
Application of lemma~\ref{lemma:interp_suff_j2} lets us to conclude that 
\begin{equation}
\label{eq:interp_j1}
\sum\limits_{j\in J_1^\delta}a_j^2 \le C_1\cdot \|u\|_{\sob}^2,
\end{equation}
where $C_1<+\infty$ depends only on $N(\Omega)$ and $\delta$, that is only on constants from the hypothesis of the theorem.

Recall that, by the choice of set of indices, if $j \notin J_1^\delta$ then $\diam_H B_j < \delta$.

\begin{lemma}
	\label{lemma:road_exist}
	Let $B_j, B_k$ be two holes adjacent in the graph $g(\Omega, s)$, that is,  
	$\Cap2(B_j, B_k) \ge s$. Suppose that $\diam_H B_j < \delta$ where constant 
	$\delta$ was chosen in the beginning of proof of theorem 
	\ref{th:interp_sufficient}; suppose also that $\diam_H B_j \le \diam_H 
	B_k$. 
	
	There exists a Borel set $A_{jk}\subset \mathbb D$ such that:
	\begin{enumerate}
		\item If $u\in\mao$  then $\left|(u|_{B_j}) - (u|_{B_k})\right|^2 \le C_2\cdot \int\limits_{A_{jk}} |\nabla u|^2 \dl.$
		\item  Inequalities $\dist_H(z, B_j) \le C_3, \,\dist_H(z, B_k) \le C_3$ are held for any point $z\in A_{jk}$. Moreover, sets $A_{jk}\cap B_j,\, A_{jk}\cap B_k$ are non-empty and even infinite.
	\end{enumerate}
	Constants $C_2, C_3<+\infty$ depend only on $\delta$ and $s$.
\end{lemma}	

\noindent {\bf Proof.} All the quantities acting in lemma are conformally invariant. Thus, by an application, if necessary, of a M{\"o}bius automorphism of disk $\mathbb D$, we may assume that $0\in B_j$. 
Euclidean  distance is no more than the hyperbolic one, therefore $\diam B_j \le \delta$.

By the choice of $F$ made in the beginning of the proof of theorem~\ref{th:interp_sufficient}, we have $\dist(B_j, B_k) \le F \cdot \diam B_j$. Thus there exists $z_1\in B_k\colon |z_1| \le (F+1) \cdot \diam B_j$.  By the choice of $\delta$ we have estimate $|z_1| \le 1/3$. Applying, if necessary, a rotation,  we may assume that $z_1 \in (0, 1)$.

Put $\Lambda=\diam B_j$. Then $\diam_H B_j \ge \Lambda$. Let us show that $\diam B_k \ge \Lambda/4$. Assume the contrary. Then, accounting that $z_1 \in B_k,\, |z_1| \le (F+1)\Lambda$, we conclude that $B_k \subset \mathcal B(0, (F+2) \Lambda)$. In this disk, by the choice of $\delta$, the density of metric $\dfrac{2|dz|\hspace{4pt}}{1-|z|^2}$ with respect to Euclidean metric is no more than $4$. Then, estimate $\diam B_k< \Lambda/4$  implies that  $\diam_H B_k < \Lambda$. But $\diam_H B_k \ge \diam_H B_j \ge \Lambda$. We got a contradiction, hence $\diam B_k \ge \Lambda/4$. 

For $t \in [0, \Lambda/8]$ construct a contour $\beta_t$ formed by two circles $\partial\mathcal B(0, t)$ and  $\partial\mathcal B(z_1, t)$ and also by two line segments $[0, z_1]+it$ and $[0, z_1]-it$. This contour is connected. Since $\diam B_j, \diam B_k \ge \Lambda/4$,  sets $B_j$ and $B_k$ intersect $\beta_t$ for all $t < \Lambda/8$ (because $0\in B_j,\, z_1 \in B_k$).

Let $u\in\mao, \, u|_{B_j} = a_j, \, u|_{B_k} = a_k$. Let us estimate  $|a_j - 
a_k|$. Repairing, if necessary, the function $u$, we can assume that for almost 
every $t\in [0, \Lambda/8]$ function~$u$ is absolutely continuous on  $\beta_t$ 
and takes values $a_j$ and $a_k$ on this contour. For all such~$t$ we have 
\begin{multline}
\label{estim:interp_suff2}
|a_j - a_k| \le \int\limits_{\beta_t} |\nabla u|\,d\mathcal H^1 \le 
\sqrt{4\pi t+2|z_1|}\cdot\left(\int\limits_{\beta_t} |\nabla u|^2\,d\mathcal H^1\right)^{1/2} \le \\ \le
\sqrt{\Lambda(\pi/2+2F+2)}\cdot\left(\int\limits_{\beta_t} |\nabla u|^2\,d\mathcal H^1\right)^{1/2}
\end{multline}
(we used the fact that length of $\beta_t$ is $4\pi t+2|z_1|$).

Define the set $A_{jk}:= \bigcup\limits_{t\in[0, \Lambda/8]} \beta_t$. Any 
point in $A_{jk}$ lies on $\beta_t$ for no more than three values of  $t$. 
Squaring~(\ref{estim:interp_suff2}) and integrating over $t\in[0, \Lambda/8]$ 
give us inequality
$$
\frac{\Lambda|a_j - a_k|^2}{8} \le 3\Lambda(2F+4)\cdot \int\limits_{A_{jk}} |\nabla u|^2\dl,
$$
which implies the first assertion of lemma with $C_2=48(F+2)$, this factor depends only on $F$, that is only on $s$. 

Now check the second assertion. By the construction, the distance from the 
origin to each of the point of contour $\beta_t$ is no more than 
$\Lambda(F+9/8)< \delta(F+2)$. Then, according to the assumptions made on the 
choice of  $\delta$, for any $z\in A_{jk}$ we have $|z| < 1/3$; but also $0\in 
B_j$, $z_1 \in B_k$,  $|z_1| < 1/3$. From this we have $\dist_H(z, B_j) \le 
C_3, \,\dist_H(z, B_k) \le C_3$ with an absolute constant $C_3<+\infty$. 
Finally, sets $A_{jk}\cap B_j,\,A_{jk}\cap B_k$ are infinite since~$\beta_t$ 
intersects $B_j$ and $B_k$ for any $t<\Lambda/8$. Lemma is proved.
$\blacksquare$

\medskip

Let us, for any pair of indices $j, k=1,2,\dots$ such that $j\not\in J_1^\delta$ or $k\notin J_1^\delta$, construct sets~$A_{jk}$ the existence of which was stated in lemma~\ref{lemma:road_exist}. (Application of this lemma is valid by the definition of set $J^\delta_1$.) \emph{Overlapness multiplicity of such sets $A_{jk}$ does not exceed $N(\Omega, C_3)^2$ where $C_3$ is the constant from lemma~\ref{lemma:road_exist}.} Indeed, if $z\in A_{jk}$ then, by the second assertion of this lemma, $\dist_H(z, B_j) < C_3, \, \dist_H(z, B_k) < C_3$. But disk $\mathcal B_H(z, C_3)$ can intersect no more than  $N(\Omega, C_3)$ of holes $B_j$. Therefore point $z$ cannot lie in more than $N(\Omega, C_3)^2$ sets of the form $A_{jk}$.

Now estimate $\sum\limits_{j\in J_2^\delta} (u|_{B_j})^2$ from the above through $\|u\|^2_{\sob}$. If $j\in J_2^\delta$ then there exists $k=k(j) \in J_1^\delta$ such that vertices $B_j$ and $B_{k(j)}$ are adjacent in the graph $g(\Omega, s)$. By the first assertion of lemma~\ref{lemma:road_exist} we have 
\begin{equation}
\label{estim:interp_suff3}
a_j^2 \le 2 a_{k(j)}^2 + 2C_2 \int\limits_{A_{jk(j)}}|\nabla u|^2 \dl.
\end{equation}
For each $j\in J_2^\delta$ put a point $\zeta_j$ on the set $A_{jk(j)} \cap B_{k(j)}$ such that all the marked points are distinct (this can be done by the second assertion of lemma~\ref{lemma:road_exist}); denote by $X$ the set of chosen points. 

Let us check the uniform local finiteness of set $X$. If some disk $\mathcal B_H(z_0, 1)$ contains point $\zeta_j, \, j\in J^\delta_2$,  then, by the second assertion of lemma~\ref{lemma:road_exist},  $\mathcal B_H(z_0, 1+C_3)\cap B_j\neq \varnothing$. Therefore $\mathcal B_H(z_0, 1)$ contains no more than $N(\Omega, C_3+1)$ of points $\zeta_j$, the desired.

Let $k \in J_1^\delta$, put $n_k:=\card(X \cap B_k)$ (this number is finite due to uniform local finiteness property of set $X$ and compactness of set $B_k$). Now sum up inequalities~(\ref{estim:interp_suff3}) over all $j \in J_2^\delta$: we have
$$
\sum\limits_{j\in J_2^\delta} a_j^2 \le 2\sum\limits_{k\in J_1^\delta} n_k a_k^2 + 2C_2 \sum\limits_{j\in J_2^\delta} \int_{A_{jk(j)}} |\nabla u|^2 \dl \le
2\sum\limits_{k\in J_1^\delta} n_k a_k^2 + 2C_2 N(\Omega, C_3)^2\cdot \|u\|_\sob^2 
$$
since overlapness multiplicity of sets $A_{jk(j)}, \, j\in J^\delta_2,$ is no more than $ N(\Omega, C_3)^2$. By lemma~\ref{lemma:interp_suff_j2} applied to the family $\{B_k\}_{k\in J^\delta_1}$ and the constructed set $X$, we have an estimate of the form $\sum\limits_{k\in J_1^\delta} n_k a_k^2\le C_4\cdot \|u\|_\sob^2$ where $C_4<+\infty$ depends only on $s$ and $N(\Omega)$. So, we have estimate 
\begin{equation}
\label{eq:interp_j2}
\sum\limits_{j\in J_2^\delta} a_j^2 \le C_5 \cdot \|u\|_\sob^2
\end{equation}
with a constant $C_5<+\infty$, depending only on constants from the hypothesis of theorem and not on function $u\in\mao$. (Recall that $s$ depends only on $S$.)

It remains to estimate $\sum\limits_{j\notin  J_1^\delta\cup J_2^\delta} a_j^2$. 
\emph{Degree of any vertex $B_j, \, j \notin J^\delta_1$, in graph $g(\Omega, s)$ is no more than $N(\Omega)$.} Indeed, consider such a vertex $B_j$. By application of a M{\"o}bius automorphism of disk $\mathbb D$ we may assume that $0\in B_j$. By the definition of set $J_1^\delta$, we have $\diam_H(B_j) \le \delta$ if $j\notin J^\delta_1$; then $\diam B_j \le \delta$. By the choice of $\delta$, if $B_k$ is adjacent with $B_j$ in the graph $g(\Omega, s)$ then $B_k \cap \mathcal B(0,1/3)\neq \varnothing$. The number of such sets does not exceed $N(\Omega)$, then the degree of vertex $B_j$ is no more than $N(\Omega)$.

\medskip

Note that degree of a vertex $B_j$ can be large if $\diam_H(B_j)$ is big, this can occur for $j\in J^\delta_1$. That is why we do estimates for $a_j, \, j\in J_2^\delta$ separately. If add estimate $\sup\limits_{j\in\mathbb N} \diam_H(B_j) < +\infty$ to the hypothesis of the theorem then we will be able to avoid application of lemma~\ref{lemma:interp_suff_j2} and estimate $a_j$ for all $j \notin J_1^\delta$ in the same manner. 

\medskip

For each $B_{j_0},\, j_0\notin J_1^\delta\cup J_2^\delta,$ there exists a path in the graph $g(\Omega, s)$ of length no more than $M$ (where $M \in \mathbb N$ is the constant from the hypothesis of theorem), connecting~$B_{j_0}$ with one of the vertices $B_l, \, l \in J_2^\delta$. Indeed, there exists a path with no more than  $M$ edges connecting $B_{j_0}$ and $\mathbb D^{(c)}$. The vertex in this path preceding its end lies in~$J_1^\delta$ (indeed, if $\Cap2(B_k, \mathbb D^{(c)}) \ge s$ then, by the choice of $\delta$, we have $\diam_H(B_k) \ge \delta$, and therefore $k \in J_1^\delta$). Hence, this path contains some vertices with indices from $J_1^\delta$ and therefore some vertices with indices from $J_2^\delta$. Taking the last of them as $B_l$ we obtain a path $(B_{j_0}, B_{j_1}, \dots, B_{j_m})$ in the graph $g(\Omega, s)$ such that $j_1, j_2, \dots, j_m \notin J_1^\delta$ but $j_m \in J^\delta_2$. At the same time, $m\le M$. Denote by $\mathcal M$ the set of all such paths, each of them corresponds to some set $B_{j_0}$,  $j_0 \notin J_1^\delta\cup J_2^\delta$.

Since $j_0, j_1, \dots, j_m \notin J_1^\delta$, sets $A_{j_k j_{k+1}}, \, k=0, 1, \dots, m-1,$ are defined. By construction from lemma~\ref{lemma:road_exist}, for all such $k$ we have estimate 
$$
a_{j_k}^2 \le 2 a_{j_{k+1}}^2 + 2 C_2\cdot \int\limits_{A_{j_k j_{k+1}}} |\nabla u|^2\dl.
$$
Applying this for $k=0,1,\dots, m-1$ we obtain an estimate of the form 
\begin{equation}
\label{estim:interp_suff4}
a_{j_0}^2 \le C_6\cdot\left(a_{j_m}^2+\sum\limits_{k=0}^{m-1}\int_{A_{j_k j_{k+1}}} |\nabla u|^2\dl\right),
\end{equation}
where $C_6 <+\infty$ depends only on $m$ and $C_1$. Since $m \le M$, then $C_6$ can be bounded from the above only in terms of constants from the hypothesis of the theorem. 

All paths from $\mathcal M$ have lengths no more than $M$ and go by vertices $B_j, \, j \notin J^\delta_1$. Degrees of such vertices do not exceed $N(\Omega)$. Thus there exists $L\in\mathbb N$ such that each edge in the graph $g(\Omega, s)$ belongs to no more than $L$ paths of class $\mathcal M$ and,  moreover, any vertex $B_j, \, j \in J^\delta_2,$ is the endpoint of no more than $L$ paths of class $\mathcal M$. Summation of estimates~(\ref{estim:interp_suff4}) over all $j_0 \notin J_1^\delta \cup J_2^\delta$ and taking into account that sets~$A_{jk}$ have overlapness multiplicity no more than $N(\Omega, C_3)^2$ allows us to conclude that 
$$
\sum\limits_{j \notin J_1^\delta \cup J_2^\delta}a_j^2 \le
C_6L\cdot \sum\limits_{k \in J_2^\delta}a_k^2 +
C_6L\cdot N(\Omega, C_3)^2\cdot \int\limits_{\mathbb D}|\nabla u|^2 \dl.
$$
Since the quantity $\sum\limits_{k \in J_2^\delta}a_k^2$ has already been estimated (inequality~(\ref{eq:interp_j2})), we obtain an estimate of the form 
\begin{equation}
\label{eq:interp_j3}
\sum\limits_{j \notin J_1^\delta \cup J_2^\delta}a_j^2 \le
C_7\cdot \|u\|_{\sob}^2,
\end{equation}
where $C_7 <+\infty$ depends only on constants from the hypothesis of the theorem. 
Gathering estimates~(\ref{eq:interp_j1}),~(\ref{eq:interp_j2}) and~(\ref{eq:interp_j3}), we conclude that if $u\in \mao$ then 
$$
\sum\limits_{j =1}^\infty a_j^2 \le C_8\cdot\|u\|_{\sob}^2
$$
with some constant $C_8 <+\infty$ depending only on constants from the hypothesis of the theorem. Proof is finished. $\blacksquare$

\medskip

\noindent {\bf Remark.} Interpolation property can fall if we erase some holes 
of domain $\Omega$. Consider, for example,  domain $\tilde\Omega_\delta$ with 
two holes that we constructed in the example~\ref{example:inverse} in 
subsection~\ref{section:examples}. Quantities $C_I(\tilde \Omega_\delta)$ are 
bounded uniformly. If we erase the greatest hole~$\tilde B_1$ in this domain 
then interpolation constant of domain $\mathbb D\setminus \tilde B_2^\delta$ 
will become large if $\delta$ is small. (This follows from theorem 
\ref{patch_theorem} and smallness of quantity $\diam_H(\tilde B_2 ^\delta)$.) 
Performing a construction by use of conformal shifts like in example 
\ref{example:inverse} in subsection~\ref{section:examples} one can obtain a 
domain $\Omega$ possessing interpolation property but such that this property 
will fail after erasing of a countable family of holes in $\Omega$.

In terms of graph $G(\Omega, S)$ this non-monotonicity by domain can be understood as follows: connectedness of a graph may fall if we remove some of its vertices. If we erase hole $B_j$ in domain at fig.~\ref{fig:tree} then vertex $B_k$ will become not connected with $\mathbb D^{(c)}$ in the graph $G(\Omega\cup B_j, S)$ if we, of course, do not increase $S$.

The non-monotonicity just mentioned can be also explained geometrically. Let $\Omega_1$ be a domain obtained by erasing of some holes $B_j, \, j\in J,$ in domain $\Omega$. Constant $C_I(\Omega_1)$ can be large whereas $C_I(\Omega)$ be small under the following  circumstances. Suppose that form $\omega\in L^{2,1}_c(\Omega_1)$ has given periods in $\Omega_1$ and minimizes the integral $\displaystyle\int\limits_{\Omega_1}|\omega|^2\dl$ under this condition. It may turn that the main mass of this integral falls into the union of the holes erased and then the integral 
$\displaystyle\int\limits_{\Omega}|\omega|^2\dl$
may turn to be not large.

Bessel property, to the opposite, is preserved if we erase in $\Omega$ some holes. Bessel constant does not increase after this removal. This can be seen immediately from the definition of Bessel property.

\subsection{Complete interpolation criterion}

Theorem~\ref{th:uniform_local_finiteness}, proposition~\ref{predl:weak_separated_bessel}, theorem~\ref{th:interp_sufficient}, theorem~\ref{th:bessel_sufficient} and theorem~\ref{th:interp_tree} which are all already proved immediately imply the following 

\begin{theorem}[complete interpolation criterion]
	\label{th:full_interp_criterion}
	Let $\Omega$ be a regular domain obtained by removing from disk $\mathbb D$  connected  holes $B_j, \, j=1,2,\dots,$ with smooth boundaries and not accumulating to inner points of  $\mathbb D$. Domain $\Omega$ possesses complete interpolation property \emph{(}that is, range of operator $\Per$ defined by $\Omega$ coincides with  $\ell^2$\emph{)} if and only the following conditions are held.
	\begin{enumerate}
		\item Family of holes $\{B_j\}_{j=1}^\infty$ is uniformly locally finite, that is, any disk of radius $1$ in the hyperbolic metric intersects no more than $N=N(\Omega)<+\infty$ of holes $B_j$.
		\item $\sup\limits_{j\in\mathbb N} \diam_H(B_j)<+\infty$ where $\diam_H$ is hyperbolic diameter.
		\item Holes $B_j$ are $\eps$-weakly separated with some $\eps>0$ that is: $$\dist(B_j, B_k) \ge \eps \cdot \min\{\diam B_j, \diam B_k\}$$ for all $j, k =1,2,\dots, \, j\neq k$.
		\item For some $S<+\infty$ graph $G(\Omega, S)$ defined at the beginning of subsection~\ref{subsec:interp_criterion} is connected and its diameter $M$ is finite.
	\end{enumerate}
	
	Moreover, if $\Omega$ has complete interpolation property then quantities $M, S, N(\Omega)$ and $\sup\limits_{j\in\mathbb N} \diam_H(B_j)$  can be estimated from the above whereas $\eps$ can be estimated  from the below  by some constants depending only on $C_I(\Omega)$ and $C_B(\Omega)$.

	Vice versa, if conditions 1 -- 4 are held then constants $C_I(\Omega)$ and $C_B(\Omega)$ can be estimated from the above only in terms of $\eps, M, S, N(\Omega)$ and $\sup\limits_{j\in\mathbb N} \diam_H(B_j)$.
\end{theorem}

\noindent {\bf Remark.} Author suspects that there exist abstract countable graphs whose connectedness is difficult to check mathematically or algorithmically. In our problem the situation is simplified, since we have condition of finiteness of diameter of $G(\Omega, S)$ in theorems~\ref{th:interp_sufficient} and \ref{th:interp_tree}. If also $\sup\limits_{j\in\mathbb N} \diam_H(B_j)<+\infty$ and holes in $\Omega$ are uniformly locally finite then degrees of vertices $B_j$ in $G(\Omega, S)$ can be estimated from the above. This allows to check (at least, algorithmically) the inequality $\dist_{G(\Omega, S)} (B_j, \mathbb D^{(c)})\le M$ for given $S, M$ and $j$, of course, if we know everything about metrics of holes.

\section{Comments}

%
%


\subsection{Blaschke condition for regular domain}

\label{subsection:Blaschke}

As a side result of the proofs of theorem~\ref{th:uniform_local_finiteness} on uniform local finiteness and theorem~\ref{patch_theorem} we obtained the estimate $\sum_{j=1}^\infty \dist(B_j, \partial\mathbb D)^2 <+\infty$. It is naturally to ask: \emph{for which exponents $\alpha \in \R$ the series $\sum_{j=1}^\infty \dist(B_j, \partial\mathbb D)^\alpha$ will necessarily converge if domain $\Omega$ has complete interpolation property?} The answer can be easily obtained only from uniform local finiteness property: \emph{for $\alpha >1$}. 

Indeed, suppose that domain $\Omega$ has uniform local finiteness property (by 
theorem~\ref{th:uniform_local_finiteness}, this is the case if $\Omega$ 
possesses complete interpolation property). In each hole~$B_j$ mark a point $z_j$ 
closest to circle $\partial\mathbb D$. For $n=1,2,\dots$ and $k=0,1,2,\dots, 
2^n-1$, consider dyadic sectors 
$$
Q_{n,k} = \{r e^{i\theta}\colon r \in (1-2^{-n+1}, 1-2^{-n}], \, \theta \in [2\pi k\cdot 2^{-n}, 2\pi (k+1) \cdot 2^{-n})\}.
$$
We have $\diam_H Q_{n,k} \le 20$. Therefore each $Q_{n,k}$ contains no more than $N(\Omega, 20)$ of marked points $z_j$. If $z_j\in Q_{n,k}$ then $1-|z_j| \le 2^{-n+1}$. Hence
\begin{multline*}
\sum_{j=1}^\infty \dist(B_j, \partial\mathbb D)^\alpha \le \sum_{j=1}^\infty(1-|z_j|)^\alpha 
\le \sum_{n=1}^\infty \sum_{k=0}^{2^n-1} \sum_{j\colon z_j \in Q_{n,k}}(1-|z_j|)^\alpha 
\le \\ \le \sum_{n=1}^\infty 2^n\cdot N(\Omega, 20)\cdot 2^{-\alpha(n-1)} < +\infty, 
\end{multline*}
if $\alpha >1$.

\emph{The series $\sum_{j=1}^\infty \dist(B_j, \partial\mathbb D)$ may not converge when domain $\Omega$ has complete interpolation property.} Indeed, in any dyadic sector $Q_{n,k}$ $(n=1,2,\dots; \, k=0,1,2,\dots, 2^n-1)$, place a round hole $B_{n,k}$ centered in $r e^{i\theta}$ where $r=1-3\cdot 2^{-n-1}, \, \theta = 2\pi(k+1/2)\cdot 2^{-n}$, set the radius of hole $B_{n,k}$ to be $2^{-n}/100$. It is easy to see that holes $B_{n,k}$ defined in such a way are strongly separated and their hyperbolic diameters are bounded from the above and separated from zero. Then we may apply theorem~\ref{th:predv_criterii}, and domain $\Omega$ with holes $B_{n,k}$ has complete interpolation property. At the same time, $\dist(B_{n,k}, \partial\mathbb D) \ge 2^{-n}$, therefore $\sum\limits_{n=1}^{\infty}\sum\limits_{k=0}^{2^n-1} \dist(B_{n,k}, \partial\mathbb D) \ge \sum_{n=1}^\infty 2^n \cdot 2^{-n} = +\infty$.

\subsection{Domains with non-smooth holes}

\label{subsection:nonsmooth}

Our investigation of  problem on interpolation by periods is, in general, based on reproducing kernels. In this way we time and again (in proofs of theorem~\ref{th:reprodicing_existence} and proposition~\ref{criterii_adm}) needed smoothness of these kernels up to the boundaries of holes in the domain. We would like to consider  problem on interpolation by periods also in domains with non-smooth holes\footnote{ Appropriateness of consideration of such domains, in particular, domains with line slits, was pointed by E.L. Korotyaev. }. Here we will not rely on reproducing kernels. Instead we will pass to limit by approximation of non-smooth holes by smooth ones. The statement of our problem does not use smoothness of holes; in the criteria we obtained we sometimes use capacity (which was defined by a way appropriate only for sets "fat"{} enough), but we can avoid this and rewrite the involved conditions through metrics. We have to check that our problem as well as the answers  obtained admit a pass to the limit if we approximate domain with arbitrary holes by domains with locally smooth boundaries in an appropriate way. The idea of such an approximation is rather clear; we will meet some minor technical difficulties only in lemma~\ref{lemma:domain_sequence} and proposition~\ref{predl:interp_necessary_approx}.

So, let $B_1, B_2, \dots$ be connected disjoint compact sets (continua) in disk $\mathbb D$. We will suppose that each $B_j$ does not separate plane $\mathbb C$ and consists of more than one point; also that sets $B_j, \, j=1,2,\dots,$ accumulate only to the boundary of disk $\mathbb D$. Put $\Omega := \mathbb D \setminus \bigcup_{j=1}^\infty B_j$. As in the above, we say that sets $B_j, \, j=1,2,\dots,$ are holes in domain $\Omega$, and such a domain $\Omega$ itself will be called a \emph{domain with arbitrary holes}. 

The space $\lo$ is defined in the same way as in the case of domains with smooth holes. It is easy to show that for each $j=1,2,\dots$ there exists a curve $\gamma_j \subset \Omega$ winding around $B_j$ once in the positive direction and not winding around the other holes in $\Omega$. Functional $\Per^{(\Omega)}_j(\omega) = \int_{\gamma_j}\omega, \, \omega \in \lo,$ is well-defined and continuous on $\lo$. Operator $\Per^{(\Omega)} \colon \lo \to \ell^2$ is defined in the same way as in the case of smooth holes. Bessel, weak Bessel, interpolation and complete interpolation properties and also constants $C_B(\Omega), \, \tilde C_B(\Omega)$ and $C_I(\Omega)$ for the domains of the above-mentioned type are introduced like in definition~\ref{def:interp_problem}. Proposition~\ref{predl:conformal_quasiconformal_invariance} stays true for domains with arbitrary holes. Corollary~\ref{sled:monotonicity_domain} on monotonicity by domain also remains true for such domains because this can be proved directly through forms and with no use of reproducing kernels.

Let $\Omega$ be a domain with arbitrary holes and  $B_1, B_2, \dots$ be bounded 
connected components of $\Omega^{(c)}$. Let $\Omega_m, \, m=1,2,\dots,$ be a 
sequence of regular domains and $B_{1,m},B_{2,m}, \dots$ be  holes in 
$\Omega_m$. Let us say that the sequence $\{\Omega_m\}_{m=1}^\infty$ 
\emph{increases to domain $\Omega$ nicely} if for each $j=1,2,\dots$ sets 
$B_{j,m}$ decrease to hole  $B_j$ when $m\to+\infty$.   

\begin{lemma}
	\label{lemma:domain_sequence}
	Let $\Omega$ be a domain with arbitrary holes. There exists a sequence of domains $\Omega_m, \, m=1,2,\dots,$ nicely increasing to domain $\Omega$ for which 
	\begin{enumerate}
		\item quantities $C_B(\Omega_m)$ decrease to  $C_B(\Omega)$;
		\item quantities $\tilde C_B(\Omega_m)$ decrease to $\tilde C_B(\Omega)$;
		\item quantities $C_I(\Omega_m)$ increase to $C_I(\Omega)$.
	\end{enumerate}
\end{lemma}

The proof is given in the Appendix. This is the only time in our work 
when we make an essential use of theory of conformal mappings. The third assertion of lemma~\ref{lemma:domain_sequence} is true for any nicely increasing sequence of domains. First and second limit passes may fall in the case of unsuccessful choice of such a sequence. It is enough to take  domain with conformal symmetry constructed in the example~\ref{example:conformal_symmerty} in subsection \ref{section:examples} and replace round holes in this domain by disks concentric with them and of small hyperbolic diameters one-by-one.

\begin{sled}
	If domain $\Omega$ with arbitrary holes has weak Bessel property then it has Bessel property, moreover, $C_B(\Omega) \le \sqrt 2 \tilde C_B(\Omega)$.
\end{sled}

\noindent {\bf Proof.}
It is enough to apply theorem~\ref{th:weak_strong_bessel} to the sequence of domains $\Omega_m$ constructed in lemma~\ref{lemma:domain_sequence} and pass to a limit in the estimate $C_B(\Omega_m) \le \sqrt 2 \tilde C_B(\Omega_m)$.
$\blacksquare$

\medskip

\noindent {\bf Remark.} If we enlarge any of domains constructed in lemma 
\ref{lemma:domain_sequence} such that the obtained sequence will still consist 
of regular domains and will increase to $\Omega$ monotonically then the 
assertion of lemma~\ref{lemma:domain_sequence} will be still true. This follows 
from the monotonicity of all the constants from this lemma.

\medskip

Definition~\ref{def:ULF} of uniform local finiteness property and constants $N(\Omega), N(\Omega, d)$ remains the same for domains with arbitrary holes.

\begin{predl}
	\label{predl:ulf_approx}
	If domain $\Omega$ with arbitrary holes possesses uniform local finiteness 
	property then the sequence of domains $\Omega_m, \,m=1,2,\dots,$ 
	constructed in lemma~\ref{lemma:domain_sequence} can be also subjected to 
	uniform local finiteness condition in order to have $N(\Omega_m, d) \le 
	N(\Omega,d+1)$ for all $m=1,2,\dots$ and $d\ge 0$.
\end{predl}

\noindent {\bf Proof.}
Indeed, for this it is enough to act in such a way that for all $m=1,2,\dots$ and $j=1,2,\dots$ set $B_{j,m}$ lies in $1$-hyperbolic neighbourhood of set $B_j$. But it is easy to achieve by shrinking holes in domains $\Omega_m$ preserving regularity of these domains and monotonicity of sequence formed by them.
$\blacksquare$

\begin{predl}
	If domain $\Omega$ with arbitrary holes  possesses complete interpolation property then it has uniform local finiteness property. Moreover, $N(\Omega)$ can be estimated from the above in terms of only $C_I(\Omega)$ and $C_B(\Omega)$.
\end{predl}

\noindent {\bf Proof.}
If $\Omega_m, \, m=1,2,\dots,$ is the sequence of domains constructed in lemma~\ref{lemma:domain_sequence}, then,  according to estimates from this lemma, 
there exists at least one $m_0=1,2,\dots$ such that domain $\Omega_{m_0}$ has 
complete interpolation property with $C_I(\Omega_{m_0}) \le C_I(\Omega), \, 
C_B(\Omega_{m_0}) \le 2 C_B(\Omega)$. By theorem 
\ref{th:uniform_local_finiteness}, the quantity $N(\Omega_{m_0})$ is finite and 
can be estimated only through $C_I(\Omega)$ and $C_B(\Omega)$. But so does 
$N(\Omega)$ since it is no more than $N(\Omega_{m_0})$.
$\blacksquare$
\medskip

Definition~\ref{def:separatedness} of $\eps$-weak separatedness of holes remains the same for domains with arbitrary holes.

\begin{predl}
If a domain $\Omega$ with arbitrary holes possesses Bessel property then holes in this domain are $\eps$-weakly separated with some $\eps>0$ depending only on $C_B(\Omega)$.	
\end{predl}

\noindent {\bf  Proof.} The sequence of domains $\Omega_m, \,m=1,2,\dots,$ 
constructed in lemma~\ref{lemma:domain_sequence}, by virtue of estimates from 
this lemma, contains a domain $\Omega_{m_0}$ for which $C_B(\Omega_{m_0}) < 
2C_B(\Omega)$. By proposition~\ref{predl:weak_separated_bessel},  holes in 
$\Omega_{m_0}$ are then $\eps$-weakly separated with some $\eps>0$ depending 
only on  $C_B(\Omega)$. But then the same will be true for $\Omega$ as well.
$\blacksquare$

\begin{predl}
	\label{predl:bessel_suff_approx}
	If domain $\Omega$ with arbitrary holes possesses uniform local finiteness property and  holes in this  domain  are $\eps$-weakly separated with $\eps>0$ and also $\sup\limits_{j\in\mathbb N} \diam_H(B_j) < +\infty$, then $\Omega$ has Bessel property. Moreover, the constant  $C_B(\Omega)$ can be estimated from the above in terms of only $\eps,\,\sup\limits_{j\in\mathbb N} \diam_H(B_j)$ and $N(\Omega)$.
\end{predl}

The proof is given in the Appendix.

Graph $G(\Omega, S)$ for $S\in(0,+\infty)$ and for domain $\Omega$ with arbitrary holes is defined in the same manner as in the beginning of subsection~\ref{subsec:interp_criterion}. This definition involves only metric characteristics of holes of domain but not capacities.

Suppose that a sequence of domains $\Omega_m, \,m=1,2,\dots,$ increases nicely 
to domain~$\Omega$ and  let $B_{1,m},B_{2,m},\dots$ be holes in $\Omega_m$. For 
a fixed $j_0=1,2,\dots$, let us identify corresponding vertices  $B_{j_0,m}$ in 
all the graphs of the form $G(\Omega_m, S)$ with different $S\in(0,+\infty)$ 
and $m=1,2,\dots$.  Moreover, we identify vertices $\mathbb D^{(c)}$ in all the graphs $G(\Omega_m, S)$.

\begin{predl}
	\label{predl:interp_necessary_approx}
	Let $\Omega$ be a domain with arbitrary holes and $N(\Omega)<+\infty$. If $\Omega$ possesses interpolation property, then there exist $S<+\infty$ and $M \in \mathbb N$ such that graph $G(\Omega, S)$ is connected and  $\dist_{G(\Omega, S)}(B_j, \mathbb D^{(c)}) \le M$  for any hole $B_j$.
	
	Numbers $S$ and $M$ can be estimated from the above in terms of only $C_I(\Omega)$ and~$N(\Omega)$.
\end{predl}

The proof is given in the Appendix. Here limit pass in the sequence of graphs $G(\Omega_m, S)$ requires a certain caution because edge-wise limit of a sequence of connected graphs may turn to be non-connected if degrees of vertices of these graphs may be infinite.

\begin{predl}
	Let $\Omega$ be a domain with arbitrary holes and $N(\Omega) <+\infty$. If for some $S<+\infty$ graph $G(\Omega, S)$ is connected and its diameter does not exceed some number $M\in\mathbb N$ then  $\Omega$ has interpolation property. Moreover, interpolation constant~$C_I(\Omega)$ can be estimated from the above only through $N(\Omega), S$  and $M$.
\end{predl}

\noindent {\bf Proof.}
Let $\Omega_m, \,m=1,2,\dots,$ be the sequence of domains constructed in proposition~\ref{predl:ulf_approx}. Identify vertices of graphs $G(\Omega_m, S)$ to ones of $G(\Omega, S)$ as before proposition~\ref{predl:interp_necessary_approx}. By the definition of these graphs, any edge from $G(\Omega, S)$ presents also in the graph $G(\Omega_m, S)$. Therefore, graphs $G(\Omega_m, S)$ are connected for all $m=1,2,\dots$ and their diameters are no more than $M$. Moreover, $N(\Omega_m) \le N(\Omega, 2)$ by proposition~\ref{predl:ulf_approx}. Hence, constants $C_I(\Omega_m)$ are finite and admit a simultaneous upper estimate only in terms of constants from the hypothesis. But   $C_I(\Omega_m)\xrightarrow{m\to\infty}C_I(\Omega)$, thus the constant~$C_I(\Omega)$ is also finite and admits the same estimate.
$\blacksquare$

\medskip

The results proved in this subsection immediately imply 

\begin{sled}
The assertion of theorem~\ref{th:full_interp_criterion} stays true for domains with arbitrary holes.
\end{sled}

\subsection{Some open questions}

\label{subsec:open_final}

\setcounter{paragraph}{0}

\paragraph{Higher dimensions.} A straightforward generalization of our results into the case of higher dimensions seems to be difficult. Suppose that we study the question on interpolation by periods for  differential forms of degree $k\in \mathbb N$ integrable with an exponent $p\in [1, +\infty]$ over an open set in  $\R^n, \, n=2,3,\dots$ (or at $n$-dimensional Riemann manifold). Period operator in this case is naturally to consider as an operator with range in  $\ell^p$. It is convenient to have scale-invariance of such a problem, the most of capacities estimates are based on this fact. This invariance will be the case if $kp=n$. 

To the other hand, theorem~\ref{th:reprodicing_existence} on explicit form of period reproducing kernels will be true in the same form if $k=n-1$; in this case duality passes a problem on forms to a problem on scalar functions. To have scale-invariance we then have to take $p=n/(n-1)$. This particular case seems to be degenerated to the author's taste. Nevertheless, even in this situation it is not possible to translate the arguments leading to the solution of the problem for $n=2$ to the higher dimension. This is because period reproducing kernels will still be gradients of harmonic functions. Writing down the Riesz basis conditions for linear combinations of such kernels leads us to \emph{estimates of harmonic gradients in $L^n$-norm} (since exponent $n$ is conjugate to $p=n/(n-1)$). Pass from Riesz basis conditions to ones of proposition~\ref{criterii_adm} was possible due to the fact that harmonic functions minimize $L^2$-norm of gradient under fixed Dirichlet boundary data. But this pass will not already be possible if we work with estimates of $L^n$-norms of such gradients when $n>2$. Therefore even in the case $k=n-1, \, p=n/(n-1)$ our arguments can not be generalized directly to the higher dimensions.

\paragraph{Riemann surfaces.} Author thinks that results of our study may be 
generalized into the case of Riemann surfaces. Then we should give up the 
simple definition of domain  $\Omega$ as a subset in $\mathbb C$ under which we 
were easily able to formulate metric conditions on holes $B_j$. Author suspects 
that the problem can be solved on Riemann surfaces as well and extremal lengths 
should stay instead of metrics.

\paragraph{Normed interpolation.} One more generalization of our problem is investigation of \emph{normed period} operator
$$
\omega \mapsto \left\{{\Per_j(\omega)}/{\|\Per_j\|_{\left(\lo\right)^*}}\right\}_{j=1}^\infty, ~~ \omega\in \lo.
$$
Namely, the problem on \emph{interpolation by normed periods} is to clarify when such an operator is bounded and has bounded right-inverse as an operator form  $\lo$ to~$\ell^2$. Such a problem differs from the studied problem on \emph{unnormed} interpolation by periods by the fact that in case of unnormed  interpolation we shall assume that norms $\|\Per_j\|_{(\lo)^*}$,  $j=1,2,\dots,$ are separated from zero and bounded from the above. The last condition (weak Bessel property) gives us a lot of useful information on geometry of domain $\Omega$.

Author has not paid a serious attention to the above-mentioned  problem on interpolation by normed periods. Nevertheless, let us point out  a case in which such a statement may occur to be relevant. Consider a nested sequence of domains $\Omega_1 \supset \Omega_2 \supset \dots$ convergent to Sierpi{\'n}ski triangle $K$. If constants $C_B(\Omega_j)$ and $C_I(\Omega_j)$ would stay bounded with $j\to\infty$ then we would be able to pass to a limit over domains  and get a good Hilbert (co)homology calculus on $K$. But, by virtue of theorem~\ref{th:uniform_local_finiteness} on uniform local finiteness, Bessel  or interpolation constants of domains approximating~$K$ should increase infinitely when domains tend to  $K$, that's because family of holes in~$K$ is not uniformly locally  finite. (Under a natural choice of domains $\Omega_j$ we have $C_B(\Omega_j)\to\infty$.) Therefore, theorem~\ref{th:uniform_local_finiteness}  on uniform local finiteness, that easily allowed us to describe complete interpolation property of a regular domain,  is also a negative result which makes it impossible to pass to a limit by domains with locally-smooth boundaries to domains with more complicated, that is, fractal configurations of holes. It may occur that such a limit pass will be possible exactly in the problem on normed interpolation by periods. For Sierpi{\'n}ski triangle, it seems to be possible to make direct calculations relying on the rich resource of its symmetries and in terms of only intrinsic geometry of this set.

\paragraph{Integer Riesz bases in Hilbert homologies.} Finally, let us state one more generalization of problem on interpolation  by periods -- a problem on \emph{integer Riesz basis in Hilbert homologies}. 

Let $\Omega$ be a countably-connected domain in $\mathbb C$. Spaces $H^1_{L^2}(\Omega)$ and $H_{1,L^2}(\Omega)$ of Hilbert cohomologies and homologies in  $\Omega$ are defined as in the remark in the end of subsection~\ref{subsection:Riesz}. These definitions make sense in the case when $\Omega$ is a \emph{countably-connected Riemann surface}.

The space $H_{1,L^2}(\Omega)$ is Hilbert and there are, of course, orthonormal bases and Riesz bases in this space. But, from geometrical viewpoint, the fact of their presence is too abstract. The question is are there "geometrically natural"{} Riesz bases in $H_{1,L^2}(\Omega)$.  For us, countable systems of \emph{integer homologies} in $\Omega$ seem to be a natural ones. Let us explain this notion. If $\beta_1, \beta_2, \dots, \beta_n$ are some closed smooth oriented curves in $\Omega$ while numbers $a_1, a_2, \dots, a_n\in \mathbb Z$, then a linear combination $a_1\beta_1+a_2\beta_2+\dots +a_n \beta_n$ gives a bounded linear functional on $H^1_{L^2}(\Omega)$, namely, a form $\omega$ defined up to an addition of an exact form is mapped to the number $\sum_{j=1}^n a_n \int_{\beta_n}\omega$. If set $\Omega$ is connected then, for integer coefficients $a_1, a_2, \dots, a_n$, a linear combination $a_1\beta_1+a_2\beta_2+\dots +a_n \beta_n$ equals in $H_{1,L^2}(\Omega)$  to a homology given by some closed curve $\beta\in \Omega$. Therefore integer $L^2\mbox{-homologies}$ in $\Omega$ can be defined as functionals on $H^1_{L^2}(\Omega)$ which can be represented as $\omega \mapsto \int_{\beta}\omega, \, \omega \in H^1_{L^2}(\Omega),$ where $\beta$ is some smooth closed curve in $\Omega$.

Consider the following property of a domain or a Riemann surface $\Omega$:
\begin{itemize}
	\item[(\dag)] \emph{In the space of Hilbert homologies $H_{1,L^2}(\Omega)$,   there exists a Riesz basis consisting of integer homologies.}
\end{itemize}
Proposition~\ref{predl:conformal_quasiconformal_invariance} implies that property (\dag) is preserved under application to $\Omega$ of a quasiconformal diffeomorphism. Note that the examples of inverse domains (examples~\ref{example:inverse} and~\ref{example:inversion_composition} in subsection~\ref{section:examples}) show that even one inversion can throw our problem out of the class of problems on interpolation by periods \emph{around holes}. That is, if we concern this problem as a conformally invariant one then its statement is attached to a conformal implementation of domain $\Omega$ as a subset of plane. Property (\dag) is free of such disadvantage, since there is no preferred homology system in this property. In other words, this property is geometrically invariant.

Now let us translate property (\dag) into analytical language. Let $\Omega$ be 
a regular domain (the class of such domains was introduced in subsection 
\ref{subsection:definitions}). Let, as usual, closed curves $\gamma_j\subset 
\Omega$ be winding around the corresponding holes $B_j, \, j=1,2,\dots$. 
Identify these curves with homologies from $H_{1,L^2}(\Omega)$ given by them  
(one may also work with kernels~$\kappa_j$ reproducing periods along curves 
$\gamma_j$). Property (\dag) of domain $\Omega$ is equivalent to the following 
one: \emph{there exists a Riesz basis in  $H_{1,L^2}(\Omega)$ whose elements 
are finite integer linear combinations of homologies  $\gamma_j$,  $j=1,2,\dots$.}

Let $H$ be an abstract Hilbert space and  $\{x_j\}_{j=1}^\infty$ be a countable system of vectors in $H$. Consider the following property of this system:
\begin{itemize}
	\item[(\ddag)] \emph{There exists a Riesz basis in $H$ whose elements are integer linear combinations of vectors $x_j, \, j=1,2,\dots$.}
\end{itemize}
There are other properties similar to (\ddag): one may require that elements of system~$\{x_j\}_{j=1}^\infty$  can be represented as integer linear combinations of elements of some  Riesz basis; or that matrix of transition from $\{x_j\}_{j=1}^\infty$ to some Riesz basis has integer values and has integer (not necessarily bounded) inverse matrix.

Author does not know any investigation in such an \emph{integer Riesz basis theory}. Note only the following fact: \emph{the system $\{je_j\}_{j=1}^\infty$ where $\{e_j\}_{j=1}^\infty$ is an orthonormal basis in~$H$, does not possess property \emph{(}\ddag\emph{)}}. Indeed, suppose that $\{f_k\}_{k=1}^\infty$ is a Riesz basis in~$H$ and that each vector $f_k, \, k=1,2,\dots,$ is obtained as a integer linear combination of elements of system $\{je_j\}_{j=1}^\infty$. For each $j=1,2,\dots$ there exists $k=1,2,\dots$ such that $\langle{f_k, e_j}\rangle \neq 0$ (otherwise, system $\{f_k\}_{k=1}^\infty$ is not complete in  $H$). This means that $je_j$ is involved in the linear combination defining $f_k$ with a non-zero integer coefficient. But system $\{je_j\}_{j=1}^\infty$ is orthogonal, therefore $\|f_k\|_H \ge \|je_j\|_H = j$. Hence, the system~$\{f_k\}_{k=1}^\infty$ contains elements with arbitrarily large norms and can not be a Riesz basis in~$H$.

Property (\dag) of a regular domain $\Omega$ is equivalent to property (\ddag) of systems of homologies $\{\gamma_j\}_{j=1}^\infty$ in the space $H_{1,L^2}(\Omega)$ (or to property (\ddag) of system of period reproducing kernels $\{\kappa_j\}_{j=1}^\infty$ in the space $H_{L^2}^1(\Omega)$). \emph{There exist regular domains not having property \emph{(}\dag\emph{)}.} To check this, one may  construct a domain $\Omega$ such that, first, norms $\|\kappa_j\|_{(\lo)^*}$ tend to infinity rapidly enough and, second, scalar products $\langle{\kappa_j, \kappa_{j'}}\rangle_{\lo}$ are small with respect to norms of kernels. Fulfilment  of these conditions can be achieved by cutting from disk $\mathbb D$ of round holes $B_j$ consecutively such that quantities $\diam_H B_j$ tend to infinity rapidly enough and such that holes $B_j$ are far enough from each other in hyperbolic metric. One can prove that a domain constructed in such a manner will not possess  property (\dag): system of kernels reproducing periods around holes in such a domain is, in general, similar to the system $\{je_j\}_{j=1}^\infty$ considered in the above. At the same time, complete interpolation property of a regular domain provides its property (\dag): integer Riesz basis in  $H_{1,L^2}(\Omega)$ is given by curves $\gamma_j$ linked to holes~$B_j$. Therefore, \emph{property \emph{(}\dag\emph{)} of existence of integer Riesz basis in Hilbert homologies is a non-trivial quasiconformal invariant of countably-connected Riemann surfaces}.

\appendix

\section{Appendix: some technical proofs}

\label{section:appendix}

\subsection{Proofs from section~\ref{section:statement_simple}}

\noindent {\bf Proof of proposition~\ref{predl:conformal_quasiconformal_invariance}.} Pick $\tilde\omega\in {{L^{2,1}_c}(\tilde\Omega)}$. 
Pull-back $\varphi^\sharp\tilde\omega$ of $\omega$ by mapping $\varphi$ is closed in $\Omega$, also $\int_{\gamma_j}\varphi^\sharp\tilde\omega=\int_{\varphi(\gamma_j)}\tilde\omega, \, j=1,2,\dots$.

Let us compare $\|\varphi^\sharp \tilde\omega\|_{\lo}$ and $\|\tilde\omega\|_{L^{2,1}_c(\tilde{\Omega})}$. Let $z\in\Omega$. By the definition of quasiconformality of mapping $\varphi$ we have $|D\varphi(z)|^2 \le K \cdot J(z)$, from which 
$$|(\varphi^\sharp\tilde\omega)(z)|^2 \le |D\varphi(z)|^2\cdot|\tilde\omega(\varphi(z))|^2\le K\cdot J(z)\cdot |\tilde\omega(\varphi(z))|^2.$$
Now by change of variable we get $\int_\Omega |\varphi^\sharp\tilde\omega|^2\dl \le K\cdot\int_{\tilde\Omega} |\tilde\omega|^2\dl$. This inequality is true if $\varphi$ is $K$-quasiconformal  and turns to equality with $K=1$ if $\varphi$ is conformal. One may push-forward forms from $\Omega$ to $\tilde{\Omega}$ and the analogous estimate will be held. Therefore, $\varphi^\sharp\colon L^{2,1}_c(\tilde{\Omega}) \to \lo$  is a bounded linear operator and $\|\varphi^\sharp\|_{L^{2,1}_c(\tilde{\Omega})\to \lo}\le \sqrt{K}$,  $\|(\varphi^\sharp)^{-1}\|_{\lo\to L^{2,1}_c(\tilde{\Omega})}\le \sqrt{K}$ while period sequences of forms $\tilde\omega\in  L^{2,1}(\tilde{\Omega})$ and $\varphi^\sharp\tilde\omega\in \lo$ along given systems of curves coincide. This implies the equalities from the first assertion and the estimates from the second assertion.
$\blacksquare$

\medskip

\noindent {\bf Proof of proposition~\ref{predl:minimizer_coexact}.} Let us prove the second assertion. Let form $\omega\in \lo$ be exact. For $t \in \R$, we have $\Per\left(\omega_0(a)+t\omega\right) = a$; minimality of norm of $\omega_0(a)$ implies that 
$$\left.\dfrac{ d\|\omega_0(a)+t\,\omega\|^2_{\lo}}{dt}\right|_{t=0}=0.$$
Differentiation according to~(\ref{eq:scalar_define}) gives us 
\begin{equation}
\label{eq:norm_diff}
\int\limits_\Omega \omega\wedge*\omega_0(a)  =0,
\end{equation}
which implies the second assertion.

Let $\eta\in W^{1,2}(\mathbb D)$ be a function. Then $d\eta|_{\Omega}\in \lo$. Since any such a function~$\eta$ can be approximated in $W^{1,2}(\mathbb D)$ by smooth functions, then  $\Per_j(d\eta)=0$ for all $j=1,2,\dots$. Therefore, by~(\ref{eq:norm_diff}),  
\begin{equation}
\label{eq:norm_diff1}
\int\limits_\Omega d\eta\wedge*\omega_0(a)  =0,
\end{equation}
for any function $\eta\in W^{1,2}(\mathbb D)$. Since the last relation is true 
for all $\eta \in C_0^\infty(\Omega)$, then form $*\omega_0$ is closed. 
Further, in subsection~\ref{subsection:reproducing_kernels} we will, with no 
use of the result now proving,   construct for $j=1,2,\dots$  functions 
$\mathfrak v_j\in\sob$ such that $\langle \omega, -(*d\mathfrak v_j)\rangle_\lo 
= \Per_j \omega$ for all $\omega\in \lo$. In our case $*\omega_0\in \lo$. 
Substituting $\eta=\mathfrak v_j$ to~(\ref{eq:norm_diff1}) gives  
$\Per_j(*\omega_0(a))=0$ for all $j=1,2,\dots$. Then form $*\omega_0(a)$ is 
exact, that is, form $\omega_0(a)$ is co-exact. $\blacksquare$

\medskip

\noindent {\bf Proof of proposition~\ref{predl:Bergman}.}
Let $\omega = \omega_x dx+\omega_y dy \in \Harm^{2,1}(\Omega)$. Consider a function $f_\omega = \omega_x - i \omega_y$. Harmonicity of form $\omega$ implies analyticity of function $f_\omega$, thence $f_\omega \in \mathscr A(\Omega)$. By a computation we ensure that $\int_{\gamma_j} f_\omega(\zeta)\,d\zeta = \Per_j(\omega) + i \, \Per_j(*\omega)$.

Suppose that domain $\Omega$ possesses interpolation property for forms. Let $a^{(1)},\,a^{(2)}\in \ell^2$. Consider a harmonic form $\omega = \omega_0(a^{(1)})-(*\omega_0(a^{(2)})) = \omega_x dx+\omega_y dy$. Then $\int_{\gamma_j} f_\omega = a^{(1)}_j+i a^{(2)}_j$. Since form $*\omega_0(a^{(2)})$ is exact then $\langle \omega_0(a^{(1)}), *\omega_0(a^{(2)})\rangle_{L^{2,1}(\Omega)}=0$. Therefore $$\|f\|_{\mathscr A(\Omega)}^2 = \|\omega_0(a^{(1)})\|_{\lo}^2 + \|\omega_0(a^{(2)})\|^2_\lo\le C_I^2(\Omega)\cdot \left(\|a^{(1)}\|_{\ell^2}^2+\|a^{(2)}\|_{\ell^2}^2\right).$$ From the arbitrariness  of the choice of sequences $a^{(1)}, \, a^{(2)} \in \ell^2$ it follows that $\Omega$ possesses interpolation property for Bergman functions as well, moreover, $C_{I, \mathscr A}(\Omega) \le C_I(\Omega)$. If, vice versa,  domain $\Omega$ possesses interpolation property in the Bergman space, then, for a given real sequence $a \in \ell^2$ there exists an analytic function $f \in \mathscr A(\Omega)$ with $\int_{\gamma_j} f = a_j$ and $\|f\|_{\mathscr A(\Omega)} \le C_{I, \mathscr A} (\Omega)\|a\|_{\ell^2}$. If $f = f_1 - i f_2, f_1, f_2 \colon \Omega \to \R$, then put $\omega =f_1 dx+ f_2 dy$. Analyticity of function $f$ implies that $\omega\in \lo$ and $\Per_j(\omega) = a_j, \|\omega\|_\lo = \|f\|_{\mathscr A(\Omega)}$. Arbitrariness of the choice of sequence $a$ now implies that $C_I(\Omega) \le C_{I, \mathscr A}(\Omega)$. 

Suppose that $\Omega$ possesses  Bessel property in the Bergman space. Let $\omega \in \lo$. Let us prove that $\|\omega\|_{\lo} \ge C_{B, \mathscr A}^{-1}(\Omega) \cdot \|\Per \omega\|_{\ell^2}$. By the remark to proposition~\ref{predl:minimizer_exists}, we may assume that $\omega$ is harmonic. But in this case the desired estimate follows from the analogous estimate for function $f_\omega \in \mathscr A(\Omega)$. 

Vice versa, let $\Omega$ possess  Bessel property in $\lo$. Take a function $f=f_1-i f_2 \in \mathscr A(\Omega)$,  $f_1, f_2 \colon \Omega\to \R$, and consider a form $\omega= f_1 dx + f_2 dy\in\Harm^{2,1}(\Omega)$. Let $a^{(1)}= \Per\omega,\, a^{(2)}= \Per(*\omega)$. Since sequences $a^{(1)}, \, a^{(2)}$ are realized as periods then there exist forms $\omega_0(a^{(1)}), \, \omega_0(a^{(2)})$, these forms are co-exact. Put $\tilde \omega := \omega - \omega_0(a^{(1)})+ *\omega_0(a^{(2)})$, then $\Per \tilde{\omega}=0, \, \Per(*\tilde{\omega})=0$. The first equality implies that $\langle \omega_0(a^{(1)}), \tilde\omega\rangle_{L^{2,1}(\Omega)}=0$ whereas the second provides   $\langle *\omega_0(a^{(2)}), \tilde\omega\rangle_{L^{2,1}(\Omega)}=\langle \omega_0(a^{(2)}), *\tilde\omega\rangle_{L^{2,1}(\Omega)}=0$. Therefore, due to Bessel property of~$\Omega$  for periods of forms,  
\begin{multline*}
\|f\|_{\mathscr A(\Omega)}^2 = \|\omega\|_{L^{2,1}(\Omega)}^2 = 
\|\omega_0(a^{(1)})\|_{L^{2,1}(\Omega)}^2+
\|\omega_0(a^{(2)})\|_{L^{2,1}(\Omega)}^2+
\|\omega\|_{L^{2,1}(\Omega)}^2\ge\\ \ge \|\omega_0(a^{(1)})\|_{L^{2,1}(\Omega)}^2+
\|\omega_0(a^{(2)})\|_{L^{2,1}(\Omega)}^2 \ge 
C_B^{-2}(\Omega)\cdot \left(\|a^{(1)}\|_{\ell^2}^2+\|a^{(2)}\|_{\ell^2}^2\right) =\\= C_B^{-2}(\Omega)\cdot \left(\|a^{(1)}+ia^{(2)}\|_{\ell^2}^2\right).
\end{multline*}
But  $\int_{\gamma_j} f = a^{(1)}_j + ia^{(2)}_j$. The estimate proved implies that domain $\Omega$ has Bessel property for functions from Bergman class.
$\blacksquare$

\subsection{Proofs from subsection~\ref{subsection:reproducing_kernels}}

\noindent {\bf Proof of theorem~\ref{th:reprodicing_existence}.} To prove the existence of function $\mathfrak v_j$ it is enough to consider closed affine subspace $Y\subset \overset{\circ}{W}{}^{1,2}(\mathbb D)$ consisting of functions $u\in\sob$ such that  $u=1$ almost everywhere on $B_j$ and $u=0$ almost everywhere on $B_{j'}$ for $j\neq j'$. Notice that $Y\neq \varnothing$ since holes $B_{j'}$ do not accumulate to $B_j$ (this easily allow to construct a function with the boundary values desired and a finite Dirichlet integral). Now, as a~$\mathfrak v_j$ we may take an element with a minimal norm in $Y$. Uniqueness of a \emph{minimizing} function follows from the fact that a sphere in a Hilbert space does not contain a segment. Equality $\Delta \mathfrak v_j=0$ in $\Omega$ is obtained by the usual variational argument.

The third assertion of  theorem follows from the second one. Let us prove the second assertion. Let a function $u_0\in \sob$ be  a solution of homogeneous boundary problem, that is, $u_0=0$ almost everywhere on each hole $B_j$,\, $j=1,2,\dots,$ and $\Delta u_0=0$ in $\Omega$. Then 
\begin{equation}
\label{eq:Laplace_gen}
\int\limits_\Omega \langle\nabla u_0, \nabla \eta\rangle \dl=0
\end{equation}
for any "test"{} function $\eta\in C^\infty_0(\Omega)$ (this is exactly the 
distributional form of relation $\Delta u_0=0$). Then, in order to prove that 
$u_0=0$ it is enough to approximate $u_0$ in 
$\overset{\circ}{W}{}^{1,2}(\mathbb D)$ by smooth functions compactly supported 
in $\Omega$ (then~(\ref{eq:Laplace_gen}) will imply that $\int_\Omega 
\langle\nabla u_0, \nabla u_0 \rangle \dl=0$ and hence $u_0=0$ almost 
everywhere in $\mathbb D$). To this end, taking a small $\delta>0$, first 
approximate $u_0$ by a function of the form $\psi_\delta(z)\cdot u_0(z)$ where 
$\psi_\delta\colon \mathbb D \to[0,1]$ is a smooth function, $\psi_\delta(z)=1$ 
for $|z|\in[0,1-2\delta],\, \psi_\delta(z) = 0$ for $|z|\in[1-\delta, 
1]$,\,$|\nabla\psi_\delta|\le 2/\delta$.  We have $\nabla(\psi_\delta u_0) - 
\nabla u_0 = (1-\psi_\delta) \nabla u_0+u_0\nabla \psi_\delta$, and 
\begin{multline*}
\|\psi_\delta u_0 - u_0\|_\sob^2 \le 
2\int\limits_{\mathbb D}\left(|\nabla u_0|^2 (1-\psi_\delta)^2+
 u_0^2 |\nabla \psi_\delta|^2\right)\dl\le \\ \le
2\int\limits_{|z|\ge 1-2\delta}|\nabla u_0|^2\dl+
8\int\limits_{|z|\ge 1-2\delta} \frac{u_0^2(z)\dl(z)}{(1-|z|)^2}.
\end{multline*}
by the choice of $\psi_\delta$. The first integral in the last expression tends 
to zero when $\delta \to 0$ since $u_0\in\sob$; the second integral also tends 
to zero when $\delta\to 0$ since, by theorem~\ref{th:Hardy_neq}, 
$$\int\limits_{\mathbb D} \frac{u_0^2(z)\dl}{(1-|z|)^2} \le \mathfrak c\cdot \|u_0\|_\sob^2 <+\infty.$$

So, functions $u_\delta=\psi_\delta u_0$ are supported in $\mathcal 
B(0,1-\delta)$ and tend to $u_0$ in $\sob$ with $\delta\to 0$. Consider a 
domain $\Omega_\delta\subset\mathbb D$ containing $\mathcal B(0,1-\delta)$ and 
all the sets $B_{j'}$, $j'=1,2,\dots,$ which intersect disk  $\mathcal 
B(0,1-\delta)$. Domain $\Omega_\delta$ can be taken in such a way that boundary 
$\partial\Omega_\delta$ is smooth. Now $u_\delta\in 
\overset{\circ}{W}{}^{1,2}(\Omega_\delta)$ and $u=0$ on all the holes $B_{j'}$, 
then~$u_\delta$ can be approximated in 
$\overset{\circ}{W}{}^{1,2}(\Omega_\delta)$ by functions of class 
$C^\infty_0(\Omega_\delta\cap\Omega)$. Here we use the fact that boundary 
$\partial\left(\Omega_\delta\cap\Omega\right)$ is smooth, that is true because 
all boundaries~$\partial B_{j'}$ are smooth. Hence, we can approximate 
$u_\delta$ in  $\sob$  by functions from  $C_0^\infty(\Omega)$, therefore $u_0$ 
admits such an approximation as well, which implies that  $u_0=0$. The second 
assertion is proved.

Now pass to the proof of the fourth assertion. First, check the case when 
domain  $\Omega$ has a smooth boundary (and, in  particular, a finite number of 
holes) and  form $\omega$, for which we want to prove~(\ref{eq:reproducing}) is 
smooth enough. Namely, let 
$$
\Omega = V\setminus \bigcup_{k=1}^K B_k,
$$
where $V\subset \mathbb C$ is a domain with boundary smooth enough whereas 
$B_k,\,k=1,2,\dots, K$, are disjoint compact subsets in $V$ also having smooth 
boundaries (the number of holes is finite). In this case it is well-known (see, 
e.g.,~\cite{LU}) that classic solution of Dirichlet problem posed in analogy to 
problem~(\ref{reproducing_eq}) with replacement of $\mathbb D$ by $V$, exists 
and is smooth up to the boundary $\partial\Omega$ (the problem is stated in the 
class $\overset{\circ}{W}{}^{1,2}(V)$ instead of  $\sob$). Let form $\omega\in 
\lo$ be smooth up to the boundary  $\partial \Omega$. Then 
\begin{equation}\label{reprod_compute}
\langle\omega, -(*d{\mathfrak v}_j)\rangle_\lo = \int\limits_\Omega \omega\wedge\left(-\left(**d{\mathfrak v}_j\right)\right) = \int\limits_\Omega \omega\wedge d{\mathfrak v}_j  = -\int\limits_\Omega d({\mathfrak v}_j\wedge\omega) = -\int\limits_{\partial\Omega} {\mathfrak v}_j\omega 
\end{equation}
(application of Stokes' theorem is valid due to sufficient smoothness of forms involved and of boundaries of sets). Further, ${\mathfrak v}_j=0$ on $\partial\Omega\setminus \partial B_j$ and ${\mathfrak v}_j=1$ on $\partial B_j$ which implies that the last integral in~(\ref{reprod_compute}) equals 
$$
-\int\limits_{\partial B_j} \omega = \Per_j\omega
$$
(again, due to smoothness of form $\omega$ up to $\partial\Omega$; boundary  
$\partial B_j$ is oriented clockwise to have accordance with the orientation of 
$\Omega$; but periods are calculated as integrals in counter-clockwise 
direction, by our definition). Therefore, in the case of a domain $\Omega$ with 
boundary smooth enough and a form $\omega$ smooth enough equality 
(\ref{eq:reproducing}) for $\kappa_j = -(*d\mathfrak v_j)$ is proved.

Now we pass the the general case. Fix $j$. Let function ${\mathfrak v}_j$ be defined for domain $\Omega$ as in the above, that is, by projection. We have to prove that 
\begin{equation}
\label{Per_eq}
\langle \omega, -(*d{\mathfrak v}_j)\rangle_\lo = \Per_j\omega
\end{equation}
for any form $\omega\in \lo$. Notice that $\langle dv, -(*d{\mathfrak 
v}_j)\rangle_\lo = 0$ for any function  $v\in C^\infty_0(\Omega)$: indeed, 
$$
\langle dv, -(*d{\mathfrak v}_j)\rangle_\lo = \int\limits_\Omega dv\wedge d{\mathfrak v}_j =
\int\limits_\Omega d(v\wedge d\mathfrak v_j)= 0
$$
due to compact supportedness of $v$. 
Therefore  $-(*d{\mathfrak v}_j) \bot {\mathcal F}_1(\Omega)$ (see~(\ref{Hodge2})) in $\lo$, 
moreover, $\Per_j|_{{\mathcal F}_1(\Omega)}\equiv0$. Thence, by proposition 
\ref{closed_orthogonal}, it is enough to check the equality~(\ref{Per_eq}) only 
for forms  $\omega\in \Harm^{2,1}(\Omega)$. These forms are smooth up to 
boundary in each strictly inner subdomain of  $\Omega$ having boundary smooth 
enough. Let us approximate $\Omega$ from the below by such domains.

For $m=1,2,\dots$ choose domains $\Omega_m\subset\mathbb D$ such that:
\begin{enumerate}
	\item $\clos\Omega_m \subset \Omega_{m+1}$, $\bigcup\limits_{m=1}^\infty =\Omega$;
	\item $\partial \Omega_m\in C^\infty$;
	\item set $\Omega_1$ (and so all $\Omega_m,\, m=2,3, \dots$)  separates the hole $B_j$ from $\infty\in \hat{\mathbb C}$. 
\end{enumerate}

For any $m=1,2,\dots$ the set $\mathbb C \setminus \Omega_m$ consists of only finite number of connected components. Denote by $B_j^{(m)}$ such a component including $B_j$. Let ${\mathfrak v}_j^{(m)}$ be the solution of problem~(\ref{reproducing_eq}) for set $\Omega_m$ (that is, $\Delta {\mathfrak v}_j^{(m)}=0$ in $\Omega_m$, ${\mathfrak v}^{(m)}_j=1$ on $B^{(m)}_j$, ${\mathfrak v}^{(m)}_j=0$ on $\partial \Omega_j\setminus \partial B_j^{(m)}$). Defining function ${\mathfrak v}^{(m)}_j$ to be $1$ on $B^{(m)}_j$ and to be $0$ on $\mathbb D\setminus \left(\Omega_m\cup B^{(m)}_j\right)$, we obtain a function of class $\overset{\circ}{W}{}^{1,2}(\mathbb D)$.

Function ${\mathfrak v}^{(m)}_j\colon \mathbb D\to \R$ can be described as a function minimizing Dirichlet integral $\|{\mathfrak v}\|^2_{\overset{\circ}{W}{}^{1,2}(\mathbb D)}$ under  conditions ${\mathfrak v}|_{B^{(m)}_j}=1$ almost everywhere, ${\mathfrak v}|_{\mathbb D\setminus (\Omega_m\cup B^{(m)}_j)}=0$ almost everywhere. If $m$ increases then these conditions weaken, thence
$$
\|{\mathfrak v}^{(m)}_j\|_{\overset{\circ}{W}{}^{1,2}(\mathbb D)}\le \|{\mathfrak v}^{(1)}_j\|_{\overset{\circ}{W}{}^{1,2}(\mathbb D)}.
$$
Thus sequence $\{{\mathfrak v}^{(m)}_j\}_{m=1}^\infty$ is bounded in $\overset{\circ}{W}{}^{1,2}(\mathbb D)$. Choose from $\{{\mathfrak v}^{(m)}_j\}_{j=1}^\infty$ a weakly convergent subsequence: namely, reserving the prior notation for the sequence, let us assume that ${\mathfrak v}^{(m)}_j \rightharpoondown {\mathfrak v}_0$ in $\overset{\circ}{W}{}^{1,2}(\mathbb D)$ while $\nabla {\mathfrak v}^{(m)}_j \rightharpoondown \nabla {\mathfrak v}_0$ in $L^2(\mathbb D)$ when $m\to\infty$, where~$\mathfrak v_0$ is some function from $\sob$.

Since ${\mathfrak v}^{(m)}_j = 1$ almost everywhere on $B_j$ and ${\mathfrak v}^{(m)}_j = 0$ almost everywhere on $B_{j'}$ for $j'\neq j$ (domains  $\Omega_m$ are less than $\Omega$), then  ${\mathfrak v}_0 = 1$ almost everywhere on $B_j$ and ${\mathfrak v}_0 = 0$ almost everywhere on $B_{j'}$ for $j'\neq j$ (embeddings $\overset{\circ}{W}{}^{1,2}(\mathbb D)$ into $L^2(B_{j})$ and into $L^2(B_{j'})$ are compact). Since $\bigcup\limits_{m=1}^\infty \Omega_m = \Omega$, we have $\Delta {\mathfrak v}_0=0$ in $\Omega$. The above-proved uniqueness of solution of Dirichlet problem~(\ref{reproducing_eq}) implies that ${\mathfrak v}_0={\mathfrak v}_j$.

Let us prove equality~(\ref{Per_eq}) for $\omega\in \Harm^{2,1}(\Omega)$.  Form $\omega$ is smooth up to the boundary in each of the sets $\Omega_m$, $m=1,2,\dots$. Thus, as it is already proved,  
$$
\Per_j(\omega) = \int\limits_{\Omega_m} \omega\wedge -(*d{\mathfrak v}^{(m)}_j).
$$
Define form $\omega$ to be zero on $\mathbb D\setminus \Omega$, then $\Per_j(\omega) = \int_{\mathbb D} \omega\wedge -(*d{\mathfrak v}^{(m)}_j)$ (functions ${\mathfrak v}^{(m)}_j$ are defined and locally constant on $\mathbb D\setminus \Omega_m$).
Passing to weak limit in $L^2(\mathbb D)$, we conclude that $\Per_j \omega = \int\limits_{\mathbb D} \omega\wedge -(*d{\mathfrak v}_j)$. Equality~(\ref{Per_eq}) is thus proved for harmonic forms, it was mentioned above that then this relation is true for all forms $\omega \in \lo$. Theorem is completely proved.
$\blacksquare$

\medskip

\noindent {\bf Proof of proposition~\ref{predl:negative_prod}}. 
Our goal is to apply Stokes' theorem to the integral $\int_\Omega \langle{\nabla \mathfrak v_{j}, \nabla \mathfrak v_{j'}}\rangle\dl$ for  functions $\mathfrak v_{j}$ and $\mathfrak v_{j'}$ that we found in theorem~\ref{th:reprodicing_existence}. To this end, take a positive number sequence $\{r_m\}_{m=1}^\infty$ increasing to $1$ such that domains $\Omega_m := \mathcal B(0, r_m) \cap \Omega$ have piecewise-smooth boundaries (this can be done by Sard's lemma). Moreover, we may assume that $B_j, B_{j'} \subset \Omega_m$. Functions $\mathfrak v_{j}, \mathfrak v_{j'}$ are smooth up to $\partial\Omega_m$ (since they are harmonic in $\Omega$, see, e.g.,~\cite{LU} about boundary smoothness of solutions of elliptic equations), moreover, $\Delta \mathfrak v_{j'}=0$ in $\Omega_m$, then
\begin{equation}
\label{eq:negative_prod_by_parts}
\int\limits_{\Omega_m} \langle{\nabla \mathfrak v_{j}, \nabla \mathfrak v_{j'}}\rangle\dl = \int\limits_{\partial \Omega_m} \mathfrak v_{j} \cdot \langle\nabla \mathfrak v_{j'}, \vec n\rangle\,d\mathcal H^1,
\end{equation}
where $\vec n$ is outward unit normal field on $\partial\Omega_m$ and $\mathcal 
H^1$ is the length measure. Set $\partial \Omega_m$ in the last integral is 
formed by $\Omega\cap \partial\mathcal B(0, r_m)$ and some pieces of boundary of 
$\partial \Omega$. The contribution of the last pieces into the right-hand side 
of equality~(\ref{eq:negative_prod_by_parts}) equals $\int\limits_{\partial 
B_j}  \langle\nabla \mathfrak v_{j'}, \vec n\rangle\,d\mathcal H^1$ since 
$\mathfrak v_j=0$ on $\partial B_k$ for $k\neq j$ and $\mathfrak v_{j} =1$ on $\partial 
B_j$. Moreover, 
$$
\left|\int\limits_{\,\Omega\cap \partial \mathcal B(0,r_n)} \mathfrak v_{j} \cdot \langle\nabla \mathfrak v_{j'}, \vec n\rangle\,d\mathcal H^1\right| \le \left(\int\limits_{\partial \mathcal B(0,r_n)} \mathfrak v_{j}^2\,d\mathcal H^1\right)^{1/2}\cdot \left(\int\limits_{\partial \mathcal B(0,r_n)} |\nabla \mathfrak v_{j'}|^2\,d\mathcal H^1\right)^{1/2}.
$$
The last quantity can be made arbitrarily small by the choice of sequence $r_n$ because $\mathfrak v_{j}, \mathfrak v_{j'}\in\sob$. Thence
$$
\int\limits_{\Omega_m} \langle{\nabla \mathfrak v_{j}, \nabla \mathfrak v_{j'}}\rangle\dl\xrightarrow{m\to\infty}\int\limits_{\partial B_j} \langle\nabla \mathfrak v_{j'}, \vec n\rangle\,d\mathcal H^1.
$$
But $\int\limits_{\Omega_m} \langle{\nabla \mathfrak v_{j}, \nabla \mathfrak v_{j'}}\rangle\dl\xrightarrow{m\to\infty}\langle\nabla \mathfrak v_{j}, \nabla \mathfrak v_{j'}\rangle_{L^2(\mathbb D)}=\langle\kappa_j, \kappa_{j'}\rangle_{\lo}$ by the explicit form of reproducing kernels $\kappa_{j}, \kappa_{j'}$. So,  $\langle\kappa_j, \kappa_{j'}\rangle_{\lo} = \int\limits_{\partial B_j}  \langle\nabla \mathfrak v_{j'}, \vec n\rangle\,d\mathcal H^1$. 

Function $\mathfrak v_{j'}$ was obtained as a solution of Dirichlet boundary problem for Laplace equation in domain $\Omega$ while its boundary values on $\partial\Omega$ equal $0$ or $1$; the existence was proved by a variational argument. From this, it is not hard to derive that $\mathfrak v_{j'} \ge 0$  in  $\Omega$. This implies that $\langle\nabla \mathfrak v_{j'}, \vec n\rangle \le 0$ on $\partial B_j$ (since $\mathfrak v_{j'}$ is smooth in $\Omega$ near $\partial B_j$ and equals zero on $\partial B_j$). To the other hand, equality $\langle\nabla \mathfrak v_{j'}, \vec n\rangle = 0$  can not be fulfilled everywhere on~$\partial B_j$ due to uniqueness of solution of Cauchy problem for Laplace equation. Therefore $\langle\kappa_j, \kappa_{j'}\rangle_{\lo} = \int\limits_{\partial B_j}  \langle\nabla \mathfrak v_{j'}, \vec n\rangle\,d\mathcal H^1<0$. Proof is finished.
$\blacksquare$

\subsection{The other proofs from section~\ref{section:prerequisites}}

\noindent {\bf Proof of proposition~\ref{predl:Hodge_F4_presice}.}
If $\omega \in \Harm^{2,1}(\Omega) \cap \left(L^{2,1}(\Omega)\ominus {\mathcal F}_5(\Omega)\right)$ then $\omega\in {\mathcal F}_4(\Omega)$ by the definition of set ${\mathcal F}_4(\Omega)$: indeed, form $\omega$ is harmonic  and is orthogonal to all the exact forms and hence to all the exact harmonic forms.

Now let $\omega \in {\mathcal F}_4(\Omega)$. In order to prove that $\omega \in \Harm^{2,1}(\Omega) \cap \left(L^{2,1}(\Omega)\ominus {\mathcal F}_5(\Omega)\right)$ one should check that $\langle\omega, \eta\rangle_{L^{2,1}(\Omega)}=0$ for any exact form $\eta\in L^{2,1}(\Omega)$. Apply to form~$\eta$ decomposition~(\ref{Hodge1}):
$$
\eta = \eta_1 + \eta_2 + \eta_3, ~~ \eta_1 \in {\mathcal F}_1(\Omega),\, \eta_2 \in {\mathcal F}_2(\Omega),\, \eta_3 \in \Harm^{2,1}(\Omega). 
$$
Proposition~\ref{closed_orthogonal} implies that $\eta_2=0$ (since form $\eta$ is exact and is then closed). Further,  $\langle\omega, \eta_1\rangle_{L^{2,1}(\Omega)}=0$ by the orthogonality of spaces ${\mathcal F}_1(\Omega)$ and $\Harm^{2,1}(\Omega)$. Finally,  $\Per_j(\eta_1)=0$ for all $j=1,2,\dots$ (because all the periods of a form of a kind $du$, $u\in C^\infty_0(\Omega)$ are zero and functional $\Per_j$ is continuous on $\lo$). To the other hand,  $\Per_j(\eta)=0$ by the exactness of form $\eta$, then $\Per_j(\eta_3)=0$ for any $j=1,2,\dots$, and hence form $\eta_3$ is exact. Condition $\omega\in {\mathcal F}_4(\Omega)$ implies that $\langle\omega, \eta_3\rangle_{L^{2,1}(\Omega)}=0$. Therefore, we proved that $\langle\omega,\eta\rangle_{L^{2,1}(\Omega)}=0$, the desired.
$\blacksquare$

\medskip

\noindent {\bf Proof of proposition~\ref{predl:4th_piece_proj}.} Let
$$
\omega =\omega_1 + \omega_2 + \omega_3 + \omega_4, ~~~ \omega_\alpha \in {\mathcal F}_\alpha(\Omega) \,\, (\alpha=1,2,3,4),
$$ 
be decomposition~(\ref{Hodge2}) of form $\omega$. 
By virtue of proposition~\ref{closed_orthogonal}, $\omega_2=0$. Further,  $\Per_j(\omega_1)=0$ and $\Per_j(\omega_3)=0$ for all $j=1,2, \dots$ (since form $\omega_3$ is smooth and exact whereas form $\omega_1$ can be approximated by smooth exact forms by the definition of space~${\mathcal F}_1(\Omega)$). Hence $\Per_j(\omega) = \Per_j(\omega_4)$. The second assertion follows from the first one. $\blacksquare$

\medskip

\noindent {\bf Proof of proposition~\ref{predl:capacity_simple}.}
To prove the first assertion it is enough to cut a function~$u$ from $\inf$ in~(\ref{eq:capacity_def1}) by levels $0$ and $1$. The second assertion is obvious. 

To prove the third assertion, take functions $u_{12}$ and $u_{13}$ from $\inf$ in definition~(\ref{eq:capacity_def1}) of capacities $\Cap2(E_1, E_2)$ and $\Cap2(E_1, E_3)$ for which $u_{12}(z), u_{13}(z)\in[0,1]$ for almost all $z\in\mathbb C$. Function $u=\max \{u_{12}, u_{13}\}$ belongs to $W^{1,2}_{\loc}(\mathbb C)$ and can be substituted into $\inf$ in~(\ref{eq:capacity_def1}) defining  $\Cap2(E_1, E_2\cup E_3)$. We have $|\nabla u|\le \max\{|\nabla u_{12}|, |\nabla u_{13}|\}$, then $\int_{\mathbb C} |\nabla u|^2 \dl \le \int_{\mathbb C} |\nabla u_{12}|^2 \dl +\int_{\mathbb C} |\nabla u_{13}|^2 \dl$. It remains to pass to $\inf$ by $u_{12}$ and~$u_{13}$.~$\blacksquare$

\medskip

\noindent {\bf Proof of proposition~\ref{predl:capacity_invariance}.} The second assertion can be proved by a simple push-forward of a functions $u$ from the definitions of capacities. Let us prove the first assertion. Take a function $\tilde u$ from the definition of capacity $\Cap2(\varphi(E_1), \varphi(E_2))$. One may choose a cutting function $\psi\colon \mathbb C\to [0,1], \, \psi|_{\mathbb D}=1,\,\psi\in C_0^\infty(\mathbb C)$ such that norm $\|\nabla \tilde u - \nabla(\psi\cdot \tilde u)\|_{L^2(\mathbb C)}$ will be arbitrarily small. Now it remains to substitute the function $u=(\psi\tilde u) \circ \varphi$ to the $\inf$ defining $\Cap2(E_1, E_2)$, make a change of variables and pass to $\inf$ over $\tilde u$. $\blacksquare$

\medskip

\noindent {\bf Proof of proposition \ref{predl:reverse_counterex}.} 
Let $n=1,2,\dots,\,Ox_1x_2\dots x_n$ be some fixed orthogonal coordinate system in~$\R^n$ and  $(e_1, e_2, \dots, e_n)$  be the corresponding orthonormal basis. Further, denote by~$\ell_n^2$ the space $\R^n$ endowed with Euclidean norm, and by $\ell_n^\infty$ we denote the space $\R^n$ with the norm $\|(x_1, x_2, \dots, x_n)\|_{\ell_n^\infty}:=\max\limits_{j=1,2,\dots, n}|x_j|, \, x_1, x_2, \dots, x_n \in \R$. It is enough to construct, for any $n$, such  vectors $\xi_{1n}, \xi_{2n}, \dots, \xi_{nn} \in \ell^2_n$  that: 
\begin{gather*}
\label{eq:counterex11}
\|\xi_{jn}\|_{\ell_n^2} \le 2, ~ j=1,2,\dots, n; \\
\label{eq:counterex12}
\left\|\sum_{j=1}^n a_j \xi_{jn}\right\|_{\ell_n^2}\ge  \|\{a_j\}_{j=1}^n\|_{\ell^2_n} \mbox{ for any vector } \{a_j\}_{j=1}^n \in \ell_n^2;\\
\label{eq:counterex13}
\sup\left\{\left\|\sum_{j=1}^n a_j \xi_{jn}\right\|_{\ell_n^2}\colon \{a_j\}_{j=1}^n\in \ell_n^2, \, \|\{a_j\}_{j=1}^n\|_{\ell^2_n} \le 1\right\} \ge \sqrt{n}. 
\end{gather*}
In this case the desired system $\{\xi_j\}_{j=1}^\infty$ in $\ell^2$ can be easily constructed by identification of $\ell^2$ with $\bigoplus\limits_{n=1}^\infty \ell_n^2$.

Define $n\times n$ matrix $T_n$: for $j=1,2,\dots,n$ its  $j$-th  \emph{row} is formed by components of vector $\xi_{jn}$ in the orthonormal basis $(e_1, e_2, \dots, e_n)$ in  $\ell_n^2$. Now conditions on vectors $\xi_{1n}, \xi_{2n}, \dots, \xi_{nn}$ can be reformulated in the following way:
\begin{gather*}
\|T_n\|_{\ell^2_n\to \ell^\infty_n} \le 2;\\
\mbox{matrix } T_n \mbox{ is invertible and } \|T_n^{-1}\|_{\ell^2_n\to \ell^2_n} \le 1;\\
\|T_n\|_{\ell^2_n\to\ell^2_n} \ge \sqrt{n}. 
\end{gather*}
Put $T_ne_n:=e_1+e_2+\dots+e_n$. Choose an orthonormal basis $(f_1, f_2, \dots, f_{n-1})$ of the linear subspace $\{x\in \R^n\colon \langle{x, e_1+e_2+\dots+e_n}\rangle_{\ell_n^2}=0\}\subset \ell_n^2.$ Put $T_ne_j := f_j$ for $j=1,2,\dots, n-1$. Let $B$ be the unit ball in  $\ell_n^2$ centered in the origin. It is easy to check that $T_n(B)$ is an ellipsoid in  $\ell_n^2$ containing $B$ and lying in cube $\{(x_1, x_2, \dots, x_n)\colon |x_j| \le {2}, \, j=1,2,\dots,n\}$. This implies first two conditions on matrix~$T_n$. The last condition is fulfilled as well since $\|T_ne_n\|_{\ell^2_n}/\|e_n\|_{\ell_n^2} =\sqrt{n} \xrightarrow{n\to\infty}\infty$. Our construction is finished.
$\blacksquare$

\medskip

\noindent {\bf Proof of proposition~\ref{predl:capacity_diam_low_estim}.}
Put $T=\frac12\min\{\diam E_1, \diam E_2\}, \, \Lambda = \dist(E_1, E_2)$. From the closedness of sets $E_1, E_2$ and compactness of one of them it follows that there exist points $z_1 \in E_1, z_2 \in E_2$ such that $|z_1-z_2|=\Lambda$. By an orthogonal transform of plane we may assume that $z_1=0, \, z_2 = \Lambda$. 

By the definition of capacity,  there exists a function $u \in 
W^{1,2}_{\loc}(\mathbb C)$ such that $u \ge 1$ almost everywhere on~$E_1$, $u 
\le 0$ almost everywhere on $E_2$, $\int\limits_{\mathbb C}|\nabla u|^2\dl \le 
2 \Cap2(E_1, E_2)$. Let $t \in [0, T)$. Consider a contour  $\beta_t$ formed by 
two circles $\partial\mathcal B(0, t)$ and  $\partial\mathcal B(\Lambda, t)$ 
and also by two line segments $[0, \Lambda]+it$ and $[0, \Lambda]-it$.  
Repairing, if necessary, the function $u$ on the set of zero measure, we may 
assume that $u$ is absolutely continuous on contour $\beta_t$ for almost all $t 
\in (0, T)$.  Circle $\partial \mathcal B(0, t)$ must have non-empty 
intersection with $\Innt E_1$, whereas $\partial \mathcal B(\Lambda, t)$ --  
with  $\Innt E_2$, since $0\in E_1, \, \Lambda\in E_2$, diameters of sets $E_1$ 
and $E_2$ are more than $2t$, and these sets are connected. Therefore $\mathcal 
H^1(\beta_t\cap E_1),\mathcal H^1(\beta_t\cap E_2)>0$. Hence, for almost all 
$t\in(0,T)$ function $u$ must take values not less than $1$ and no more than 
$0$ on the contour $\beta_t$. The length of this contour $\beta_t$ is $4\pi 
t+2\Lambda$ and this contour is connected. Then, by Cauchy-Schwartz inequality, 
$$
\int\limits_{\beta_t}|\nabla u|^2 \, d\mathcal H^1 \ge \frac1{4\pi t+2\Lambda}
$$
for almost all $t\in(0, T)$. Each point $z\in\mathbb C$ belongs to contour $\beta_t$ for no more than three values of $t$. Integration of the estimate obtained by $t\in (0,T)$ gives 
$$
3\int\limits_{\mathbb C} |\nabla u|^2 \dl \ge \int\limits_0^{T} \frac{dt}{4\pi t+2\Lambda}=\frac1{4\pi}\log\left(1+\pi\cdot \frac{2T}\Lambda\right).
$$
The desired estimate now follows from the inequality $\int_{\mathbb C}|\nabla u|^2\dl \le 2 \Cap2(E_1, E_2)$.~$\blacksquare$

\medskip

\noindent {\bf Proof of proposition~\ref{predl:weak_separated_bessel}.}
Let us prove the first assertion. Let $j=1,2,\dots$. By the first assertion of proposition~\ref{necessary_cap}, we have $\Cap2(B_j, \mathbb D^{(c)}) \le \tilde C_B^2(\Omega)$. By proposition~\ref{predl:capacity_diam_low_estim} then $\min\{\diam B_j, \diam \mathbb D^{(c)})\} \le C_1\cdot \dist(B_j, \mathbb D^{(c)})$, where  $C_1$ depends only on $\tilde C_B(\Omega)$. In other words, $\dfrac{\diam B_j}{\dist(B_j, \mathbb D^{(c)})} \le C_1$. But then, by the second assertion of proposition~\ref{sled:diamh_metric}, we have $\diam_H(B_j)\le C_2$, where  $C_2$ depends only on $C_1$, that is, only on $\tilde C_B(\Omega)$. The second assertion is proved analogously. $\blacksquare$

\medskip

\noindent {\bf Proof of proposition~\ref{predl:strong_sep_Bessel}.} Let $j, k=1,2,\dots, \, j\neq k$. Let $\diam B_j \ge \diam B_k$. According to the second assertion of proposition~\ref{predl:weak_separated_bessel},  $\dist(B_j, B_k) \ge \eps_0 \cdot \diam B_k$, where $\eps_0 >0$ depends only on $\tilde C_B(\Omega)$. Let us prove that $\dist(B_j, B_k) \ge \eps \cdot \diam B_j$ with some $\eps$ depending only on the constants from the hypothesis. From $\diam_H(B_k) \ge r$ it follows that $\dist(B_k, \partial\mathbb D) \le c_1 \cdot \diam B_k$, where $c_1$ depends only on $r$ (the first assertion of proposition~\ref{sled:diamh_metric}). By the second assertion of the same proposition and also by the first assertion of proposition~\ref{predl:weak_separated_bessel} we have $\dist(B_j, \partial\mathbb D) \ge c_2 \cdot \diam B_j$ where $c_2>0$ depends only on $\tilde C_B(\Omega)$. Suppose that $\dist(B_j, B_k) \le \dfrac{c_2}{2} \cdot \diam B_j$. Then 
\begin{multline*}
c_2 \diam B_j \le \dist(B_j, \partial\mathbb D) \le \dist(B_j, B_k)+\diam B_k + \dist(B_k, \partial\mathbb D) \le \\ \le
\frac{c_2}{2} \cdot \diam B_j+(1+c_1)\diam B_k,
\end{multline*}
from which
$$
\diam B_j \le \frac{2(1+c_1)}{c_2}\cdot \diam B_k.
$$
Estimate for $\diam_H(B_k)$ now gives 
$$
\dist(B_j, B_k) \ge \eps_0\cdot \diam B_k\ge \frac{c_2\eps_0}{2(1+c_1)}\cdot \diam B_j.
$$
So, we may take $\min\left\{\dfrac{c_2}2, \dfrac{c_2\eps_0}{2(1+c_1)}\right\}$ as $\eps$, this quantity is positive and depends only on $r$ and $\tilde C_B(\Omega)$. $\blacksquare$

\medskip

\noindent {\bf Proof of proposition~\ref{predl:capacity_upper}.} Consider a  function 
$$
u(z) := \frac{\dist(z, \mathbb C \setminus U)}{\Lambda\eps},  ~~ z\in\mathbb C.
$$
We have inequalities $u\ge 1$ on $E$, $u=0$ on $\mathbb C\setminus U$. Further,  function $u$ is Lipschitz, moreover,
$$
|\nabla u| \le \frac{1}{\Lambda\eps}
$$
almost everywhere; $\nabla u=0$ outside of $U$. The area of $U$ is no more than $\pi(1+2\eps)^2 \Lambda^2$ since its diameter is no more than $(1+2\eps)\Lambda$. Hence
$$
\int\limits_{\R^2}|\nabla u|^2\dl \le \frac{\pi (1+2\eps)^2}{\eps^2}.
$$
Substitution of function $u$ into  $\inf$ from the definition of capacity leads us to the result.~$\blacksquare$

\medskip
\noindent {\bf Proof of proposition~\ref{predl:strong_sep_Bessel_suff}.} Condition $\sup\limits_{j\in\mathbb N} \diam_H(B_j)< +\infty$ implies that $\dist(B_j, \partial\mathbb D) \ge \eps_1 \cdot \diam B_j$ for each $j=1,2,\dots$ where $\eps_1>0$ and depends only on $\sup\limits_{j\in\mathbb N} \diam_H(B_j)$. Put $\eps_0 :=\min\{\eps_1/3, \eps/3\}$, this constant is positive and depends only on $\eps$ and $\sup\limits_{j\in\mathbb N} \diam_H(B_j)$. Put $\Lambda_j=\diam B_j, \, j=1,2,\dots$. According to definition of strong separatedness, sets $U_{\eps_0\Lambda_j}(B_j)$ are pairwise disjoint, moreover, $U_{\eps_0\Lambda_j}(B_j)\subset \mathbb D$  for all $j=1,2,\dots$.

Let $u_j, \, j=1,2,\dots,$ be the functions constructed in the proof of proposition~\ref{predl:capacity_upper} for set $B_j$ and complement to its neighbourhood $U_{\eps_0\Lambda_j}(B_j)$, that is, $u_j=1$ almost everywhere on $B_j$,  $u_j=0$ almost everywhere out of  $U_{\eps_0\Lambda_j}(B_j)$,  $\|u_j\|_\sob\le C$, here~$C$ depends only on $\eps_0$ and thus only on the constants from the hypothesis.

To prove Bessel property of domain $\Omega$, let us check the hypothesis of the first assertion of proposition~\ref{criterii_adm}. Let $a=\{a_j\}_{j=1}^\infty \in \ell^2$. Put $u_a:=\sum_{j=1}^\infty a_j u_j$. For finite sums of such a form we have $\|u_a\|_\sob^2 \le C^2 \cdot \|a\|_{\ell^2}$ (sets $\supp u_j$ are pairwise disjoint), hence series $\sum_{j=1}^\infty a_j u_j$ is convergent in $\sob$ for any sequence $a=\{a_j\}_{j=1}^\infty\in \ell^2$ and $\|u_a\|_\sob\le C\cdot\|a\|_{\ell^2}$. Sum $u_a$ of this series is an admissible function since each $u_j$ is admissible. Also $u_a|_{B_j}\equiv a_j$ almost everywhere. Thence, the hypothesis of proposition~\ref{criterii_adm} is checked and therefore $\Omega$ possesses  Bessel property with $C_B(\Omega) \le C$.~$\blacksquare$

\medskip

\subsection{Metric proof of theorem~\ref{th:bessel_sufficient}}



\noindent {\bf Proof of lemma \ref{lemma:annular_construction_conditions}.} Let us start the choice of constant $c_1$. Since $\sup\limits_{j\in\mathbb N} \diam_H(B_j) < +\infty$, it follows that there exists a constant $c_h>0$ such that $U_{c_h \Lambda_j}(B_j)\subset\mathbb D$ for all $j=1,2,\dots$ (see proposition~\ref{sled:diamh_metric}). We may assume that $c_1 < c_h/2$. It is easy to show that in this case $\dist_H(z, B_j) \le 2$ for any point $z\in U_{c_1 \Lambda_j} (B_j)$; thence $\diam_H U_{c_1 \Lambda_j} (B_j) \le \diam_H B_j+4$.  We will also assume that $c_1<\eps/2$ where $\eps$ is the constant of weak separatedness of holes $B_{j}$. The choice of $c_1$ will be refined in the below.

Uniform local finiteness property implies the existence of such $M\in \mathbb N$ that any set of hyperbolic diameter no more than $4+\sup\limits_{j\in\mathbb N} \diam_H(B_j)$ intersects no more than  $M$ of holes $B_j$ ($M=N(\Omega, 4+\sup\limits_{j\in\mathbb N} \diam_H B_j)$).

\begin{lemma} 
	\label{lemma:wide_neigh}
	In the hypothesis of theorem~\ref{th:bessel_sufficient} there exists a constant $C_2=C_2(\eps, M)<\infty$ such that for any $j=1,2,\dots$ and any positive $c_1 < \min\{c_h/2, \eps/2\}$ there exists an open interval $I\subset (0, c_1 \Lambda_j)$ having length not less than ${c_1 \Lambda_j}/{C_2(\eps, M)}$ for which 
	$$
	\{z\in\R^2 \colon \dist(z, B_j) \in I\} \subset \Omega. 
	$$
\end{lemma}

\noindent {\bf Proof.} Fix $j$. Let $B_{j_1}, B_{j_2},\dots, B_{j_m}$ be holes intersecting $U_{c_1 \Lambda_j}(B_j)$ (excepting the hole~$B_j$). Their number $m$ does not exceed $M$ by the choice of $c_1 < c_h/2$. Let $[t_k, s_k]$ be the image of set $B_{j_k}$ under mapping $\dist(\cdot, B_j)$  ($k=1,2,\dots, m$). This is indeed a segment since $B_{j_k}$ is a continuum. 

By the choice of $c_1$, $\dist(B_j, B_{j_k})< \eps\cdot \diam B_j$ for all $k=1,2,\dots, m$. But  holes in  domain $\Omega$ are $\eps$-weakly separated. Thence,   $\dist(B_j, B_{j_k})\ge \eps\diam B_{j_k}$. This implies that 
\begin{equation}
\label{eq:wide_neigh_lemma_separ}
\frac{s_k-t_k}{t_k} \le \frac1\eps.
\end{equation}
We are going to show that the last estimate and boundedness of $m$ imply that $[0, c_1 \Lambda_j]$ contains 
an interval $I$ with $\lambda_1(I) \ge c_2 \Lambda_j$ which does not intersect $[t_k, s_k]$; here $c_2>0$ is not very small. Take $c_2<c_1$ and suppose the contrary. We may assume that $t_1 \le t_2 \le \dots \le t_m$. If $t_1 > c_2 \Lambda_j$ then we may take interval of the form $I=[\delta , \delta + c_2 \Lambda_j]$ with $\delta >0$ small enough. Therefore $t_1 \le c_2 \Lambda_j$ then from~(\ref{eq:wide_neigh_lemma_separ}) we have $s_1 \le (1+1/\eps) c_2 \Lambda_j$. Let $k=2,3, \dots,m$. If $t_k > c_2 \Lambda_j+\max\{s_1, s_2, \dots, s_{k-1}\}$ then we again may take $I=\max\{s_1, s_2, \dots, s_{k-1}\}+\delta + [0,c_2\Lambda_j]$ with small positive  $\delta$. Thus
\begin{gather*}
t_k \le \max\{s_1, s_2, \dots, s_{k-1}\}+c_2 \Lambda_j,\\
s_k \le \left(1+\frac1\eps\right)t_k.
\end{gather*}
Consecutive  application of these estimates  for $k=2, 3, \dots,m$ gives us inequality of the form $s_k \le C(\eps, k)\cdot c_2 \Lambda_j$ where number $C(\eps, k)\in(0,+\infty)$ depends only on $k$  and $\eps$. Put  $C_2(\eps, M):=2\max\limits_{k\le M} C(\eps, k)+3$. Take $c_2 := c_1/C_2(\eps, M)$. If, on $[0, c_1\Lambda_j]$, between intervals  $[t_k, s_k]$ and also to the left of them there is no unoccupied interval of length~$c_2 \Lambda_j$, then  $s_k < c_1\Lambda_j/2$ for all $k=1,2,\dots, m$. Now we may take an interval of length $c_2 \Lambda_j$ in $[c_1\Lambda_j/2, c_1\Lambda_j)$ as $I$. $\blacksquare$

\medskip

Now we continue the proof of lemma \ref{lemma:annular_construction_conditions}. Let us refine the choice of $c_1$: in addition to the above-stated restrictions let us assume that $c_1\cdot \left(1+C_2(\eps, M)\right) < \eps$ where  $\eps$ is the constant of  weak separatedness of holes whereas number $C_2(\eps, M)$ was constructed in lemma~\ref{lemma:wide_neigh}. Put $c_2 := c_1/C_2(\eps, M)$. We have, in particular, 
\begin{equation}
\label{neq:dist_order_estim0}
c_1^2/c_2+c_1 < \eps.
\end{equation}

Denote by $I=[t_j, s_j]$ the interval constructed in the previous lemma for hole~$B_j$, $j=1,2,\dots$. Put $A_j:=U_{s_j}(B_j)\setminus U_{t_j}(B_j)$. Sets~$A_j$ satisfy all the requirements from the first assertion of (inclusion $A_j\subset \Omega$ is true by the choice of  $[t_j,s_j]$). Let us check the second assertion. Suppose that some point $z \in A_{j_1}\cap  A_{j_2}\cap\dots\cap A_{j_m}$ for some indices $j_1, j_2, \dots, j_m$. Let us estimate $m$. As it was mentioned before lemma~\ref{lemma:wide_neigh}, then $\dist_H(z, B_{j_k}) \le 2$ for all $k=1,2,\dots,m$. Consider a ball in the hyperbolic metric of radius $2$ centered in $z$. This ball intersects all the sets $B_{j_k}, \, k=1,2,\dots, m$. Hence,  $m\le N(\Omega,2)$ by the definition of constant $N(\Omega,2)$. We have checked the second assertion of lemma.

\medskip

Note that set $A_j$ can, in general, turn to be non-connected even if $B_j$ is connected and simply connected. (It is enough to take the difference of two neighbourhoods of set~$B_j$ looking like an arc of a circle of some radius $r$, the length of this arc be close to~$2\pi r$.) Even if set $U_{s_k}(B_j)\setminus U_{t_k}(B_j)$ is connected, it may turn not to be homeomorphic to an annulus. Our construction will simplify if holes are disks: then~$A_j$ are annuli. However, a mutual placement even of round annuli may occur to be complicated. Namely, it may turn that $A_j\cap A_k \neq \varnothing$ for some $j,k=1,2,\dots$, but $U_{s_j}(B_j) \not\subset U_{s_k}(B_k)$ and $U_{s_k}(B_k) \not\subset U_{s_j}(B_j)$. We have to show that our choice of $c_1$ excludes  such a possibility.

\medskip

The relation $\prec$ of strict partial order on the set of holes $B_j$ is defined as follows: $B_{j'}\prec B_j$ if and only if $j'\neq j$ and $U_{s_j'}(B_{j'})\subset U_{s_j}(B_{j})$. By the definition, this relation is transitive and anti-reflexive.

Let us check the remaining part of the third assertion of lemma. Let $B_{j'}\prec B_j$. Then $j'\neq j$ and $B_{j'}\subset U_{s_{j'}}(B_{j'})\subset U_{s_{j}}(B_{j})$ by the definition of order $\prec$. Now check that if $B_{j'}\subset U_{s_{j}}(B_{j})$ and~$j'\neq j$, then $B_{j'}\prec B_j$. We have to prove that $U_{s_{j'}}(B_{j'})\subset U_{s_{j}}(B_{j})$. 
Set $A_j=U_{s_j}(B_j)\setminus U_{t_j}(B_j)$ lies in $\Omega$ and thus does not intersect hole $B_{j'}$. Thence $B_{j'}\subset U_{t_j}(B_j)$, in particular, $\dist(B_j, B_{j'})\le t_j$. By the choice of $t_j$ and by $\eps$-weak separatedness of holes we have inequality
$$c_1 \Lambda_j \ge t_j \ge \dist(B_j, B_{j'})\ge\min\{\eps \Lambda_j, \eps \Lambda_{j'}\}.$$
By the choice of $c_1$ we have $c_1 \Lambda_j < \eps \Lambda_j$, which implies $\Lambda_{j'}< \Lambda_j$ and $c_1 \Lambda_j \ge \eps \Lambda_{j'}$, that is, $\Lambda_{j'}\le c_1 \Lambda_j/\eps$. By the choice of $s_{j'}$ we have $s_{j'}\le c_1 \Lambda_{j'}\le c_1^2 \Lambda_j/\eps$. By the choice of~$c_1$, the last quantity is strictly less than $c_2 \Lambda_j$. Therefore $s_{j'}+t_j< c_2 \Lambda_j +t_j \le s_j$.  We conclude that 
$U_{s_{j'}}(B_{j'})\subset U_{s_{j'}+t_j} (B_j)\subset U_{s_{j}}(B_j)$. 
The third assertion of lemma~\ref{lemma:annular_construction_conditions} is proved.

Now let's prove the fourth assertion in lemma~\ref{lemma:annular_construction_conditions}. Let us, for a fixed  $j=1,2,\dots$, estimate the number of holes $B_{j'}$ for which $B_{j'} \prec B_j$. For such $j'$ we have $B_{j'}\subset U_{t_j}(B_{j})$. But, as it was noticed above, $\diam_H( U_{t_j}(B_{j}))\le 4+\sup\limits_{j\in\mathbb N}\diam_H(B_j)$. Thus there exist no more than $N(\Omega,4+\sup\limits_{j\in\mathbb N}\diam_H(B_j))$ indices $j'$ for which $B_{j'}\prec B_j$. This implies that lengths of chains in order $\prec$ are finite and can be estimated from the above only through constants from the hypothesis of theorem.

Now we check that if $B_{j_1}, B_{j_2} \succ B_j$ then either $B_{j_1}\succ B_{j_2}$, or $B_{j_2}\succ B_{j_1}$. (In other words, $\{B_{j'}\colon B_{j'} \succ B_j\}$ is a chain in order $\succ$.) We may assume that $j_1=1,\,j_2=2$. Since $B_j \subset U_{t_{1}}(B_{1}) \cap U_{t_{2}}(B_{2})$ (by the third assertion of lemma), then  $U_{t_{1}}(B_{1}) \cap U_{t_{2}}(B_{2})\neq\varnothing$. Therefore,  we have, in particular,
\begin{equation}
\label{neq:dist_order_estim1}
\dist(B_1, B_2) \le t_1+t_2\le c_1(\Lambda_1+\Lambda_2).
\end{equation} 
Without limiting of generalization we may assume that $\Lambda_2 \le \Lambda_1$. From the condition of $\eps$-weak separatedness for sets $B_1$  and $B_2$ we have inequality
$\dist(B_1, B_2) \ge \eps \Lambda_2$ (because $c_1 < \eps/2$). Taking~(\ref{neq:dist_order_estim1}) in account, we obtain $\eps \Lambda_2 \le c_1(\Lambda_1+\Lambda_2)$, thence
\begin{equation}
\label{neq:dist_order_estim2}
\Lambda_1\ge \frac{\eps-c_1}{c_1}\,\Lambda_2.
\end{equation}

If $B_2\not\prec B_1$ then, by the third assertion of lemma,  $B_2\not\subset U_{t_1}(B_1)$. Distance between sets $U_{t_1}(B_1)$ and $\R^2 \setminus U_{s_1}(B_1)$ is positive whereas $B_2$ is connected and does not intersect $A_1=U_{s_1}(B_1) \setminus U_{t_1}(B_1)$. Therefore $\dist(B_2, B_1) \ge s_1\ge t_1+c_2 \Lambda_1$ by the choice of numbers $t_1, s_1$. Comparing this and~(\ref{neq:dist_order_estim1}), we have $c_2 \Lambda_1 \le t_2 \le c_1 \Lambda_2$ (the last inequality is true by the choice of $t_2$), from which $\Lambda_1 \le c_1 \Lambda_2/c_2$. This inequality and~(\ref{neq:dist_order_estim2}) imply 
$$
\frac{\eps-c_1}{c_1} \le \frac{c_1}{c_2},
$$
or $\eps \le c_1^2/c_2+c_1$. But this contradicts the inequality~(\ref{neq:dist_order_estim0}) established during the choice of $c_1$. We thus checked all the properties necessary for our construction, lemma~\ref{lemma:annular_construction_conditions} is proved.
$\blacksquare$

\medskip

\noindent {\bf Derivation of theorem~\ref{th:bessel_sufficient} from lemma \ref{lemma:annular_construction_conditions}.} Let $a_1, a_2, \dots \in \R$, $\sum_{j=1}^\infty a_j^2 <+\infty$. We are going to construct a function $u\in\mao$ for which $u|_{B_j}=a_j$, $|u\|_{\sob}^2 \le   \tilde C\cdot \sum_{j=1}^\infty a_j^2$ where constant $\tilde C<+\infty$ can be estimated from the above only through constants from the hypothesis of the theorem. It is sufficient to work out the construction only for finitely supported sequences $\{a_j\}_{j=1}^\infty$. According to proposition~\ref{criterii_adm}, this is enough to establish Bessel property of domain~$\Omega$.

Apply lemma~\ref{lemma:annular_construction_conditions} already proved. The set of holes $B_j$ is now endowed with a partial order  $\prec$. Recall that $\mathcal M$ is the set of indices of maximal elements in order~$\prec$ and "nearest ancestor"{} mapping $\mathcal P \colon \mathbb N \setminus \mathcal M \to \mathbb N$ was defined in the remark after lemma~\ref{lemma:annular_construction_conditions}.

Let $j=1,2,\dots$. Consider sets $U_{t_j}(B_j)$ and $\R^2\setminus U_{s_j}(B_j)$. Diameter of the first of them is no more than $(1+2c_1)\Lambda_j$ whereas distance between them is not less than~$c_2 \Lambda_j$. By proposition~\ref{predl:capacity_upper} there exists a function $u_j\in\sob$ such that 
$u_j=0$ in $\R^2\setminus U_{s_j}(B_j)$, $u_j=1$ in $U_{t_j}(B_j)$, $\|u_j\|_\sob \le C_3$, where constant $C_3$ depends only on numbers $c_1$  and $c_2$ worked out in lemma~\ref{lemma:annular_construction_conditions} and depending only on constants from the hypothesis of theorem. (Such a function was constructed explicitly in the proof of proposition~\ref{predl:capacity_upper}.)

Put 
\begin{equation}
\label{eq:bessel_existence}
u(z) := \sum\limits_{j\in \mathcal M} a_j u_j(z)+\sum\limits_{j\not\in \mathcal M} \left(a_j-a_{\mathcal P(j)}\right)u_j(z), ~~ z\in\mathbb D.
\end{equation}
By the fourth assertion of lemma~\ref{lemma:annular_construction_conditions}, for a fixed $j_0$ there exists only finite number of indices $j$ for which $\mathcal P(j)=j_0$. The finite supportedness of the sequence $\{a_j\}_{j=1}^\infty$ now implies that the sum in the right-hand side of~(\ref{eq:bessel_existence})  is finite. Moreover, $u\in\sob$ since all the $u_j\in\sob$. 

Let us check that $u|_{B_{j}}=a_{j}$. When $u_{j'}\neq 0$ on $B_j$? Function $u_{j'}$ equals $1$ on  $U_{t_{j'}}(B_{j'})$ and equals zero out of $U_{s_{j'}}(B_{j'})$. Connectedness of $B_j$ and estimate $$\dist\left(U_{t_{j'}}(B_{j'}), U_{s_{j'}}(B_{j'})\right)\ge s_{j'}-t_{j'}>0,$$ imply that $u_{j'} \equiv 0$ on $B_j$ if $B_{j'}\not\succ B_j$, and also $u_{j'} \equiv 1$ on $B_j$ if $B_{j'}\succ B_j$ (by the third assertion of lemma~\ref{lemma:annular_construction_conditions}). The set of holes $B_{j'}$ for which $B_{j'}\succ B_j$ is finite and forms a chain (by the fourth assertion of lemma~\ref{lemma:annular_construction_conditions}). Let it consist of holes $B_{j_1}\succ B_{j_2}\succ \dots\succ B_{j_m} \, (m\in\mathbb N)$ whereas $B_{j_m}\succ B_j$, $j_1\in\mathcal M$, $\mathcal P(j_{k+1})=j_k$  $(k=1,2,\dots, m-1)$. Now if $z\in B_j$ then 
$$
u(z)=a_{j_1}+(a_{j_2}-a_{j_1})+(a_{j_3}-a_{j_2})+\dots+(a_{j_{m}}-a_{j_{m-1}})+(a_{j}-a_{j_m})=a_j
$$
by the construction of function $u$. The desired equality is established.

Let's estimate $\|u\|_\sob$. By the construction,  $\nabla u_j\neq 0$ only on set $A_j$. Overlapness multiplicity of sets $A_j$ is no more than $C_1$ (the second assertion of lemma~\ref{lemma:annular_construction_conditions}; $C_1$ depends only on constants from the hypothesis of theorem). For  $z\in\mathbb D$ we have 
$$
\nabla u(z) = \sum\limits_{j\in \mathcal M} a_j \nabla  u_j(z)+\sum\limits_{j\not\in \mathcal M} \left(a_j-a_{\mathcal P(j)}\right)\nabla u_j(z).
$$
For fixed $z$ the last sum contains no more than $C_1$ non-zero terms. Thence
\begin{multline*}
|\nabla u(z)|^2 
\le C_1\cdot\left(\sum\limits_{j\in \mathcal M} a_j^2 \,|\nabla  u_j(z)|^2+\sum\limits_{j\not\in \mathcal M} \left(a_j-a_{\mathcal P(j)}\right)^2|\nabla u_j(z)|^2\right)\le\\ \le
 C_1\cdot\left(\sum\limits_{j\in \mathcal M} a_j^2 \,|\nabla  u_j(z)|^2+2\sum\limits_{j\not\in \mathcal M} \left(a_j^2+a^2_{\mathcal P(j)}\right)|\nabla u_j(z)|^2\right).
\end{multline*}
Integration by $\dl$ gives 
\begin{multline*}
\hspace{-4mm}
\int\limits_{\mathbb D}|\nabla u|^2\dl\le 
\sum\limits_{j\in\mathcal M} C_1 a_j^2 \|u_j\|_\sob^2+ 
\sum\limits_{j\notin\mathcal M}\left(2C_1 a_j^2 \|u_j\|_\sob^2
+2C_1 a^2_{\mathcal P(j)}\|u_j\|_\sob^2\right)\le\\ \le
C_1 C_3\cdot\sum\limits_{j\in\mathcal M}   a_j^2 +
2C_1C_3\cdot \sum\limits_{j\not\in\mathcal M} \left(a_j^2+ a^2_{\mathcal P(j)}\right) .
\end{multline*}
Since for any hole $B_j$ there exists no more than $C$ of holes $B_{j'}$ for which $\mathcal P(j')=j$, then the last sum does not exceed 
$$
(2C_1C_3+2C_1C_3C)\sum\limits_{j=1}^\infty a_j^2.
$$
The estimates proved imply that 
\begin{equation}
\label{est:norm_estimate_bessel_construct}
\|u\|^2_{\sob}\le C_4\cdot\sum_{j=1}^\infty a_j^2,
\end{equation}
where $C_4$ depends only on constants from the hypothesis of the theorem. So, we constructed, for any (finitely supported) sequence $\{a_j\}_{j=1}^\infty$, a function $u\in\mao$ such that $u|_{B_j}\equiv a_j$ and a norm estimate~(\ref{est:norm_estimate_bessel_construct}) is held. By proposition~\ref{criterii_adm}, this is sufficient for domain $\Omega$ to possess Bessel property. Also, $C_B(\Omega) \le \sqrt{C_4}$. Proof is complete.
$\blacksquare$

\subsection{Proofs from subsection~\ref{subsection:nonsmooth}}

\noindent {\bf Proof of lemma~\ref{lemma:domain_sequence}.}
The monotonicity required in all the three assertions is immediate from corollary~\ref{sled:monotonicity_domain}. 

To establish, for any nicely increasing sequence of domains, the limit pass from the third assertion, take any  sequence $a\in \ell^2$. For each  $m=1,2,\dots$ there exists a form $\omega_m \in L^{2,1}_c(\Omega_m)$ with $\Per^{(\Omega_m)}\omega_m = a$, $\|\omega_m\|_{L^{2,1}_c(\Omega_m)}\le C_I(\Omega_m) \cdot \|a\|_{\ell^2}$. Define $\omega_m$ to be zero in $\mathbb D \setminus \Omega_m$. Taking from $\{\omega_m\}_{m=1}^\infty$ a  subsequence weakly* convergent in  $L^2(\mathbb D)$, we obtain as a  weak limit a form $\omega \in \lo$ as a weak limit, for which $\Per^{(\Omega)}\omega = a, \, \|\omega\|_{\lo} \le \|a\|_{\ell^2}\cdot \varlimsup\limits_{m\to\infty}C_I(\Omega_m)$. This implies that $C_I(\Omega) \le \varlimsup\limits_{m\to\infty}C_I(\Omega_m)$, accounting the monotonicity of interpolation constants gives us the third assertion.

The limit pass for Bessel  and weak Bessel  constants requires additional considerations. 

\begin{lemma}
	\label{lemma:qc_construct}
	Let $\Omega$ be a domain with arbitrary holes, $B_1, B_2, \dots$ be bounded connected components of $\Omega^{(c)}$ and  $\lambda > 1$. There exist sets $B_1', B_1'', B_2', B_2'', \dots\subset\mathbb D$ such that:
	\begin{enumerate}
		\item for any $j=1,2,\dots$ sets $B_j', B_j''$ are simply connected closures of domains with smooth boundaries; $B_j \subset \Innt B_j'\subset B_j'\subset \Innt B_j'' \subset B_j''$; sets $B_j''$ are pairwise disjoint and accumulate only to $\partial\mathbb D$ when $j\to\infty$;
		\item there exists a quasiconformal diffeomorphism $\varphi$ of domain $\Omega$ onto domain $\Omega' := \mathbb D \setminus \bigcup_{j=1}^\infty B_j'$ for which $K(\varphi) \le \lambda^2$; moreover, $\,\varphi(z)=z$ if $z \in \Omega\setminus \bigcup_{j=1}^\infty B_j''$.
	\end{enumerate} 
\end{lemma}

\noindent {\bf Proof.} 
We do this by a step-by-step process. Suppose that, on $(j-1)$-th step ($j=1,2,\dots$), we have already constructed sets $B_1', B_1'', B_2', B_2'',\dots, B_{j-1}', B_{j-1}''$, each of them does not intersect $B_j$. Let us construct sets $B_j'$ and $B_j''$. Let $\psi_j$ be a conformal bijection of set $(\mathbb C\cup\{\infty\}) \setminus B_j$ onto $(\mathbb C\cup\{\infty\}) \setminus \clos{\mathbb D}$. Such a  mapping exists because set $B_j$ is connected, does not separate plane and consists of more than one point (see~\cite{Goluzin}). We may assume that $\psi_j(\infty)=\infty$.

Holes $B_{j+1}, B_{j+2}, \dots$ do not accumulate to $B_1$. Thence openness of mapping $\psi_j$ allows to derive easily that  set $\psi_j\left(\mathbb D \setminus\left(B_1''\cup B_2'' \cup \dots \cup B_{j-1}''\cup B_j \cup B_{j+1}\cup\dots\right)\right)$ contains an annulus of a kind $\bar{\mathcal B}(0,r_j) \setminus \clos{\mathbb D}$ with some $r_j>1$. 

Pick $\rho_j\in(1,r_j)$. Let us assume that $\rho_j$ is close to $1$ enough such that there exists a quasiconformal diffeomorphism $\tilde \varphi_j$ which maps $(\clos\mathbb D)^{(c)}$ onto $\mathbb C\setminus \bar{\mathcal B}(0,\rho_j)$, for which $K(\tilde \varphi_j) \le \lambda^2$ and $\tilde \varphi_j(z)=z$ for $|z| \ge r_j$. 

Sets $\psi_j^{-1}(\partial\mathcal B(0,\rho_j))$ and  $\psi_j^{-1}(\partial\mathcal B(0,r_j))$ are smooth Jordan curves on the plane and even in disk $\mathbb D$. By Jordan theorem, each of them encloses  a closed bounded   set, take the first of these sets as $B_j'$, and the second -- as $B_j''$. Since $\psi_j(\infty)=\infty$ whereas circle $\partial\mathcal B(0,r_j)$ separates $\infty$ from $\partial\mathcal B(0,\rho_j)$ then $B_j'\subset B_j''$ (and thus $B_j'\subset \Innt B_j''$ ). Circle $\partial\mathcal B(0,r_j)$ does not separate sets $\psi_j(B_k'')$ for $k=1,2,\dots,j-1,$ and sets $\psi_j(B_l)$ for $l=j+1,j+2,\dots$ from $\infty$; hence $B_j''$ does not intersect $B_k''$ and  $B_l$ for such $k$ and~$l$. The last circumstance allows to proceed the process of construction of sets $B_l''$. Set $\psi_j^{-1}(\mathcal B(0,\rho_j)\setminus \clos\mathbb D)$ can not be simply connected, this implies that $B_j \subset \Innt B_j'$.

We also have to provide the last assertion of the first statement of lemma, that is to act in such a way that sets $B_j''$ do accumulate only to $\partial\mathbb D$. For that, it is enough to do our construction such that set $B_j''$ lies in $1$-neighbourhood of set $B_j$ in the hyperbolic metric. Since $B_j\subset \mathbb D$, then it is enough to have $B_j''$ lying in $\eps$-neighbourhood of set $B_j$ in Euclidean metric for some $\eps>0$ depending on $B_j$; it is sufficient to check the  last condition only for curve $\psi_j^{-1}(\partial\mathcal B(0,r_j))$. The openness of mapping $\psi_j$ allows to prove easily that it can be done by shrinking, if necessary, of radius $r_j$ (before choosing~$\rho_j$).

Performing such steps for all $j=1,2,\dots$ gives us sets required in the first assertion of lemma. Let us construct the mapping $\varphi$ from the second assertion. To this end, put $\varphi(z)=z$ if $z$ does not belong to the union of constructed sets $B_j''$, $j=1,2,\dots$. If $z\in B_j''\setminus B_j$ for some $j$ then put $\varphi(z) := \left(\psi_j^{-1}\circ \tilde \varphi_j\circ \psi_j\right)(z)$. Then $\varphi$ is a smooth mapping without critical points in the whole set $\Omega$ by the construction of mappings $\tilde\varphi_j$. From this and also from conformality of mappings $\psi_j$ we have estimate $K(\varphi) \le \lambda^2$. Furthermore, annulus $\mathcal B(0,r_j)\setminus\bar{\mathcal B}(0,\rho_j)$ consists of points separated from $\infty$ by the curve $\partial\mathcal B(0,r_j)$, but not by the curve $\partial\bar{\mathcal B}(0,\rho_j)$. Thence $\psi^{-1}\left(\mathcal B(0,r_j)\setminus\bar{\mathcal B}(0,\rho_j)\right) = (\Innt B_j'')\setminus B_j'$. Therefore $\varphi$ is a diffeomorphism of $\Omega$ onto $\mathbb D \setminus \bigcup_{j=1}^\infty B_j'$.  The second assertion of lemma is proved. Proof is complete. $\blacksquare$

\medskip

Now finish the proof of lemma~\ref{lemma:domain_sequence}. For $m=1,2,\dots$ apply lemma~\ref{lemma:qc_construct} to domain $\Omega$ and $\lambda = 1+1/m$. Let $B_{j,m}',B_{j,m}'', \, j=1,2,\dots,$ be holes with smooth boundaries constructed in this lemma, and $\Omega_m' = \mathbb D \setminus \bigcup_{j=1}^\infty B_{j,m}'$, denote by $\varphi_m$ the corresponding mappings from the second assertion of this lemma.  Mapping~$\varphi_m$ is a diffeomorphism of domain $\Omega$ onto $\Omega_m'$. \emph{For any $j=1,2,\dots$, mapping~$\varphi_m$ carries a curve in $\Omega$ winding around hole $B_j$ onto a curve in domain  $\Omega_m'$ winding around~$B_{j,m}'$.} Indeed, by the first assertion of lemma~\ref{lemma:qc_construct}, we may take a contour lying in $\mathbb D\setminus \bigcup_{j=1}^\infty B_{j,m}''$ as such a curve, but on this contour mapping $\varphi_m$ is identical. It was mentioned that proposition~\ref{predl:conformal_quasiconformal_invariance} stays true for domains with arbitrary holes. Therefore  $C_B(\Omega_m') \le (1+1/m)\cdot C_B(\Omega), \, \tilde C_B(\Omega_m') \le (1+1/m)\cdot \tilde C_B(\Omega)$. So,  $\varlimsup\limits_{m\to\infty} C_B(\Omega_m') \le C_B(\Omega),\, \varlimsup\limits_{m\to\infty} \tilde C_B(\Omega_m') \le \tilde C_B(\Omega)$.

It remains to provide monotonicity of the constructed sequence of domains. By monotonicity of Bessel, weak Bessel  and interpolation constants by domains, this can be done by shrinking of constructed  domains $\Omega_m'$ in such an way that the obtained sequence of domains will decrease. But this condition can be provided by a consequential choice of domains $\Omega_1'', \Omega_2'', \dots$ such that for every $m=1,2,\dots$ domain $\Omega_m$ will stay regular, lie inside $\Omega_m'$ and inside of already constructed $\Omega_1'', \Omega_2'', \dots, \Omega_{m-1}''$ and still contain $\Omega$. Lemma \ref{lemma:domain_sequence} is proved. $\blacksquare$

\medskip

\noindent {\bf Proof of proposition~\ref{predl:bessel_suff_approx}.} Holes in domain $\Omega$ accumulate only to $\partial\mathbb D$. Thence, for each $j=1,2,\dots$ there exists $\delta>0$ such that  $\delta$-neighbourhood of hole $B_j$ in Euclidean metric intersects only finite number of other holes. This observation allows easily to construct a sequence of open sets $U_j, \, j=1,2,\dots,$ such that $B_j \subset U_j$ for all~$j$ and sets $\clos U_j$ are $\eps/2$-weakly separated. Also, it is easy to provide the estimate $\sup\limits_{j\in\mathbb N} \diam_H(U_j) < \sup\limits_{j\in\mathbb N} \diam_H(B_j)+1$ for all $j=1,2,\dots$.

For each $j=1,2,\dots$ choose a  set $B_j'$ with smooth boundary such that $B_j\subset B_j'\subset U_j$. Domain $\Omega':=\mathbb D\setminus \bigcup_{j=1}^\infty B_j'$ is regular. By the construction, holes in  $\Omega'_m$ are $\eps/2$-weakly separated, and  $\sup\limits_{j\in\mathbb N} \diam_H(B_{j,m}) \le \sup\limits_{j\in\mathbb N} \diam_H(B_j)+1$. One may also assume that $B_j'$ lies in  $1$-neighbourhood of set $B_j$ in hyperbolic metric. Then $N(\Omega_m) \le N(\Omega,2)$. Theorem~\ref{th:bessel_sufficient} implies that the constant $C_B(\Omega')$ is finite and can be estimated from the above only through $\eps,\,\sup\limits_{j\in\mathbb N} \diam_H(B_j)$ and $N(\Omega)$. But $C_B(\Omega)$ does not exceed $C_B(\Omega')$ and thus is finite and admits an estimate through the constants from the hypothesis. $\blacksquare$

\medskip

\noindent {\bf Proof of proposition~\ref{predl:interp_necessary_approx}.} Let $S'$ and $M'$  be the constants whose existence is provided by theorem~\ref{th:interp_sufficient} for domains $\Omega'$ such that $C_I(\Omega') \le C_I(\Omega)$ and $N(\Omega') \le N(\Omega, 2)$.

Let $\alpha>0$ be such that if $E\subset \mathbb D$,  $\diam_H E \le \alpha$, then $\dist(E, \partial\mathbb D) \ge 2S'\cdot \diam E$. Denote by $J$ the set of indices of holes $B_j$ in domain $\Omega$ for which $\diam_H(B_j)\le \alpha/2$. Choose number $S\ge S'$ such that if $E\subset \mathbb D, \, \diam_H E \ge \alpha/2$ then $\dist(E, \partial\mathbb D) \le S \cdot \diam E$. Put $M=M'+1$. Numbers $S$ and $M$ depend only on constants from the hypothesis of our proposition. Let us prove our assertion with such $S$ and $M$.

Choose a sequence of regular domains $\Omega_m, \, m=1,2,\dots,$ nicely increasing to $\Omega$ for which estimates from lemma \ref{lemma:domain_sequence} are true. These domains can be chosen such that if~$B_{j,m}$ is one of the holes in $\Omega_m$ and  $B_j$ is the hole in $\Omega$ corresponding to $B_{j,m}$, then $\diam_H(B_{j,m})\le {2\diam_H(B_j)}$ and $B_{j,m}$ lies in  $1$-neighbourhood of $B_j$ in hyperbolic metric. Then $N(\Omega_m) \le N(\Omega,2)$. Moreover, $C_I(\Omega_m) \le C_I(\Omega)$ by monotonicity of interpolation constant. Therefore, for each $m=1,2,\dots$  graph $G(\Omega_m, S')$ is connected, and $\dist_{G(\Omega_m, S')}(B_{j,m}, \mathbb D^{(c)}) \le M$ for all $j\in\mathbb N$. Identify the corresponding vertices of these graphs as it was mentioned before proposition \ref{predl:interp_necessary_approx}; if some edge presents in all the graphs $G(\Omega_m, S')$ then it presents in the graph $G(\Omega,S')$ as well.

\emph{If $j_0\in J$ then the degree of vertex $B_{j_0,1}$ in the graph $G(\Omega_1, S')$ is finite.} Indeed, for any such $j_0$ we have $\diam_H(B_{j_0,1}) \le \alpha$ by the construction. If vertex $B_{j,1}$ is adjacent in the graph $G(\Omega_1, S')$ to vertex $B_{j_0,1}$, then $$\dist(B_{j_0,1}, B_{j,1}) \le S'\cdot \diam B_{j_0,1} \le \dist(B_{j_0,1}, \partial\mathbb D)/2.$$ But the Euclidean $\left(\dist(B_{j_0,1}, \partial\mathbb D)/2\right)$-neighbourhood of set $B_{j_0,1}$ can intersect only a finite number of other holes in $\Omega_1$, since the closure of this neighbourhood lies strictly inside of the disk $\mathbb D$ whereas domain $\Omega_1$ is regular.

Let $j_0\in J$.  Denote by $\mathcal T_m$ the set of all paths $\tau$ in the graph $G(\Omega_m, S')$ such that: 
\begin{enumerate}
\item $\tau$ has edge-length no more than $M'$;
\item $\tau$ starts in vertex $B_{j_0,m}$;
\item all the vertices of path $\tau$, excepting its end, are of the form $B_{j,m}$ with some $j\in J$;
\item the end of path $\tau$ is either $\mathbb D^{(c)}$, or one of the holes of kind $B_{j', m}$ with $j'\notin J$.
\end{enumerate}
By monotonicity of sequence $\{\Omega_m\}_{m=1}^\infty$, the set of edges of graph  $G(\Omega_m, S')$ decreases if~$m$ increases. Then $\mathcal T_m$ decreases with the growth of $m$. Further, degrees of vertices~$B_{j,1}$ in the graph $G(\Omega_1,S)$ are finite for $j\in J$; then set $\mathcal T_1$ is finite. Since distance from $B_{j_0,m}$ to $\mathbb D^{(c)}$ in the graph  $G(\Omega_m, S')$ is no more than $M'$ then sets $\mathcal T_m$ are non-empty for all $m=1,2,\dots$.  Therefore, there exists a path $\tau$ belonging to sets $\mathcal T_m$ for all $m=1,2,\dots$. From the definitions of graphs  $G(\Omega_m,S')$ and  $G(\Omega, S')$ it follows that all edges of path $\tau$ present in the graph $G(\Omega, S')$ as well, whereas length of path  $\tau$ does not exceed $M'$. The last vertex of path $\tau$ in the graph $G(\Omega, S')$ is either $\mathbb D^{(c)}$ or such a hole  $B_j$ that $\diam_H B_j \ge \alpha/2$. By the choice of $S$, in the second case such a vertex is adjacent to  $\mathbb D^{(c)}$ in the graph  $G(\Omega, S)$. But  $\tau$ is a path in the graph $G(\Omega, S)$ as well. Therefore, vertex $B_{j_0}$ can be joined in the graph  $G(\Omega, S)$ with vertex $\mathbb D^{(c)}$ by a path of edge-length no more than $M'+1$. The same is true also for vertices of the form $B_{j'}$ for $j'\notin J'$, since these vertices are adjacent in  $G(\Omega, S)$ to $\mathbb D^{(c)}$ by the choice of  $S$. Proof is complete. $\blacksquare$


\end{document}